\pdfoutput=1
\documentclass[10pt, twoside]{article}

\usepackage{amssymb,epsfig,amsmath,enumitem,mathrsfs}
\usepackage{amsthm}
\usepackage{mathtools}
\usepackage{dsfont}
\usepackage{amsfonts}
\usepackage{indentfirst}
\usepackage{natbib}
\usepackage{graphicx}
\usepackage{caption}
\usepackage{subcaption}
\usepackage{color}
\definecolor{darkred}{rgb}{0.5,0.2,0.2}
\usepackage[pdftitle={Hybrid marked point processes: characterisation, existence and uniqueness},pdfauthor={Maxime Morariu-Patrichi and Mikko S. Pakkanen},colorlinks=TRUE,allcolors=darkred]{hyperref}
\usepackage{eso-pic}
\usepackage{chngcntr}
\usepackage[top=1.25in, bottom=1.25in, left=1in, right=1in]{geometry}
\usepackage[utf8]{inputenc}
\usepackage{datetime}
\usepackage{apptools}
\usepackage{titlesec}
\AtAppendix{\counterwithin{marked_point_process}{subsection}}

\newdateformat{mydate}{\THEDAY\ \monthname[\THEMONTH] \THEYEAR}

\titleformat{\subsubsection}[runin]
  {\normalfont\normalsize\bfseries}{\thesubsubsection}{1em}{}

\makeatletter
\def\@endtheorem{\endtrivlist}
\makeatother

\graphicspath{{plots/}}

\usepackage{fancyhdr}
\fancypagestyle{firstpg}{%
	\fancyhf{}
	\fancyfoot[C]{\small\thepage}
	\fancyhead[CO]{}
	\fancyhead[CE]{}
	
}

\newcommand{\as}{\mbox{a.s.}} 
\newcommand{\aeverywhere}{\mbox{-a.e.}} 
\newcommand{\indicator}{\mathds{1}} 
\newcommand{\borelSetsOf}[1]{\mathcal{B}(#1)} 
\newcommand{\realLine}{\mathbb{R}}  
\newcommand{\realLineNonNegative}{\mathbb{R}_{\geq 0}} 
\newcommand{\realLineNonPositive}{\mathbb{R}_{\leq 0}} 
\newcommand{\realLineNegative}{\mathbb{R}_{<0}} 
\newcommand{\realLinePositive}{\mathbb{R}_{>0}} 
\newcommand{\integers}{\mathbb{N}} 
\newcommand{\anyCSMS}{\mathcal{U}}  
\newcommand{\anyCSMSBorel}{\borelSetsOf{\anyCSMS}}  
\newcommand{\integerValuedMeasures}[1]{\mathcal{N}^{\infty}_{#1}}  
\newcommand{\boundedlyFiniteMeasures}[1]{\mathcal{N}^{\#}_{#1}} 
\newcommand{\boundedlyFiniteMeasuresBorel}[1]{\mathcal{B}(\boundedlyFiniteMeasures{#1})} 
\newcommand{\boundedlyFiniteMeasuresMeasurableSpace}[1]{(\boundedlyFiniteMeasures{#1}, \boundedlyFiniteMeasuresBorel{#1} )} 
\newcommand{\stateSpace}{\mathscr{X}}  
\newcommand{\stateSpaceMeasure}{\mu_{\stateSpace}} 
\newcommand{\stateSpaceMeasurable}{(\stateSpace, \borelSetsOf{\stateSpace}, \stateSpaceMeasure )}  
\newcommand{\eventSpace}{\mathscr{E}}  
\newcommand{\eventSpaceMeasure}{\mu_{\eventSpace}} 
\newcommand{\eventSpaceMeasurable}{( \eventSpace, \borelSetsOf{\eventSpace}, \eventSpaceMeasure )} 
\newcommand{\markSpace}{\mathscr{M}}  
\newcommand{\markSpaceBorel}{\mathcal{B}(\markSpace)} 
\newcommand{\markSpaceMeasure}{\mu_{\markSpace}} 
\newcommand{\markSpaceMeasurable}{ ( \markSpace, \markSpaceBorel, \markSpaceMeasure ) } 
\newcommand{\theCSMS}{\realLine\times\markSpace} 
\newcommand{\csmsBorel}{\mathcal{B}(\realLine\times\markSpace)} 
\newcommand{\mppSpace}{\mathcal{N}_{\realLine\times\markSpace}^{\# g}}  
\newcommand{\mppSpaceBorel}{\mathcal{B}(\mppSpace)} 
\newcommand{\mppMeasureSpace}{( \mppSpace, \mppSpaceBorel)}  
\newcommand{\weakHash}{w^{\#}} 
\newcommand{\weakHashDistance}[2]{d^{\#}(#1,#2)} 
\newcommand{\weakHashDistancePlain}{d^{\#}} 
\newcommand{\realisationsSpace}{\Omega}  
\newcommand{\sigmaAlgebra}{\mathcal{F}} 
\newcommand{\prob}{\mathbb{P}} 
\newcommand{\measureSpace}{( \realisationsSpace, \sigmaAlgebra)} 
\newcommand{\probSpace}{( \realisationsSpace, \sigmaAlgebra , \prob )} 
\newcommand{\E}[1]{\mathbb{E}\left[ #1 \right]}  
\newcommand{\internalHistory}[2]{\mathcal{F}^{#1}_{#2}} 
\newcommand{\internalHistoryCollection}[1]{\mathbb{F}^{#1}}
\newcommand{\collection}[2]{(#1)_{#2}} 
\newcommand{\aHistory}[1]{\mathcal{F}_{#1}} 
\newcommand{\aHistoryCollection}{\mathbb{F}}
\newcommand{\predictableSigmaAlgebra}{\mathcal{F}^{p}} 
\newcommand{\intensity}{\lambda} 
\newcommand{\intensityAtTimeAtMarkWithIndex}[3]{\lambda_{#3}(#1,#2)} 
\newcommand{\intensityFunctional}{\psi} 
\newcommand{\intensityFunctionalAtMarkGivenPastWithIndex}[3]{\intensityFunctional^{#3}(#1 \,|\, #2)} 
\newcommand{\stateFunctional}{\phi} 
\newcommand{\stateFunctionalAtStateGivenEventGivenPast}[3]{\phi(#1 \,|\, #2,#3)} 
\newcommand{\eventFunctional}{\eta} 
\newcommand{\eventFunctionalAtEventGivenPast}[2]{\eventFunctional(#1 \,|\, #2)} 
\newcommand{\stateProcess}{X} 
\newcommand{\shift}[1]{\theta_{#1}} 
\newcommand{\restrictionToNegativeOfMeasure}[1]{#1^{<0}} 
\newcommand{\restrictionToNonPositiveOfMeasure}[1]{#1^{\leq 0}} 
\newcommand{\pointProcessAtTime}[2]{\shift{#2}\restrictionToNegativeOfMeasure{#1}}  
\newcommand{\restrictionToNonNegativeOfMeasure}[1]{#1^{\geq 0}} 
\newcommand{\restrictionToPositiveOfMeasure}[1]{#1^{>0}} 
\newcommand{\pointProcessAfterTime}[2]{\shift{#2}\restrictionToPositiveOfMeasure{#1}}  
\newcommand{\inducedProb}[1]{\mathcal{P}^{#1}} 
\newcommand{\realisationsSpaceInitialCondition}{\Omega_{\leq 0}} 
\newcommand{\sigmaAlgebraInitialCondition}{\mathcal{F}_{\leq 0}} 
\newcommand{\probInitialCondition}{\mathbb{P}_{\leq 0}} 
\newcommand{\probSpaceInitialCondition}{( \realisationsSpaceInitialCondition, \sigmaAlgebraInitialCondition, \probInitialCondition )} 
\newcommand{\realisationsSpaceFuture}{\Omega_{> 0}} 
\newcommand{\sigmaAlgebraFuture}{\mathcal{F}_{> 0}} 
\newcommand{\probFuture}{\mathbb{P}_{> 0}} 
\newcommand{\probSpaceFuture}{( \realisationsSpaceFuture, \sigmaAlgebraFuture, \probFuture )} 
\newcommand{\drivingPoisson}{M} 
\newcommand{\dominatingMarkedPointProcess}{\overline{N}} 
\newcommand{\dominatingIntensityFunctional}{\overline{\psi}} 
\newcommand{\dominatingIntensityFunctionalAtMarkGivenPastWithIndex}[3]{\dominatingIntensityFunctional^{#3}(#1 \,|\, #2)} 
\newcommand{\dominatingIntensity}{\overline{\lambda}} 
\newcommand{\dominatingIntensityAtTimeAtMarkWithIndex}[3]{\dominatingIntensity_{#3}(#1,#2)} 

\begin{document}

\thispagestyle{firstpg}
\begin{center}
\large \textbf{HYBRID MARKED POINT PROCESSES: \\ CHARACTERISATION, EXISTENCE AND UNIQUENESS} \\
	\vspace{1cm}
	\large  Maxime Morariu-Patrichi\footnote{Corresponding author. Department of Mathematics, Imperial College London, South Kensington Campus, London SW7 2AZ, UK. E-mail: \href{mailto:m.morariu-patrichi14@imperial.ac.uk}{\texttt{m.morariu-patrichi14@imperial.ac.uk}} URL: \url{http://www.maximemorariu.com}} \hspace{1cm} Mikko S. Pakkanen\footnote{Department of Mathematics, Imperial College London, South Kensington Campus, London SW7 2AZ, UK and CREATES, Aarhus University, Aarhus, Denmark. E-mail: \href{mailto:m.pakkanen@imperial.ac.uk}{\texttt{m.pakkanen@imperial.ac.uk}} URL: \url{http://www.mikkopakkanen.fi}} $\mbox{\hspace{0cm}}$ \\
	\vspace{0.5cm}
	\normalsize {\mydate\today}
	\vspace{1cm}
\end{center}

\begin{abstract}
We introduce a class of hybrid marked point processes, which encompasses and extends continuous-time Markov chains and Hawkes processes. While this flexible class amalgamates such existing processes, it also contains novel processes with complex dynamics. These processes are defined implicitly via their intensity and are endowed with a state process that interacts with past-dependent events. The key example we entertain is an extension of a Hawkes process, a state-dependent Hawkes process interacting with its state process. We show the existence and uniqueness of hybrid marked point processes under general assumptions, extending the results of \citet{Massoulie:1998aa:StabilityResults} on interacting point processes.
\end{abstract}

\bigskip

\bigskip

\numberwithin{equation}{section}

\section{Introduction} \label{sec:intro}

Let $N=(N_1,\ldots,N_d)$ be a $d$-dimensional counting process, meaning that $N_i(t)$, $i=1,\ldots,d$, is the number of events of type $i$ that have occurred until time $t$ \citep{bremaud:1981:MartingaleDynamics, daleyVereJonesVolume1, last1995marked, Sigman:1995aa, jacobsen2006point}.
A fundamental concept describing the dynamics of $N$ is the intensity process $\intensity=(\intensity_1,\ldots,\intensity_d)$. Loosely speaking, when it exists, the intensity $\intensity_i(t)$ of $N_i$ at time $t$ is such that
\begin{equation*}
	\E{N_i(t+dt)-N_i(t)\,|\,\internalHistory{N}{t}} \approx \intensity_i(t)dt,
\end{equation*}
where $\internalHistoryCollection{N}=\collection{\internalHistory{N}{t}}{t\geq 0}$ is the natural filtration generated by $N$. Intuitively, the above equation says that $\intensity_i(t)dt$ is the expected number of events of type $i$ in the infinitesimal time window $(t,t+dt]$, given what has happened so far.

Besides describing counting processes, intensity processes can actually be used to specify them implicitly. A prime example of that is the class of linear Hawkes processes \citep{Hawkes:1971aa, Hawkes:1974aaClusterRepresentation}, which are characterised by the intensities
\begin{equation} \label{eq:linear_hawkes}
	\intensity_i (t) = \nu_i + \sum_{j=1}^{d}\int_{[0,t)} k_{ji}(t-s)dN_j(s),\quad t\geq 0,\, i=1,\ldots,d,
\end{equation}
where each $\nu_i\in\realLinePositive$ is fixed and each $k_{ji}:\realLineNonNegative\rightarrow\realLineNonNegative$ is a non-negative function, usually called a \emph{kernel}, see~\citet{laub:2015:Hawkes} for an introduction. This class of processes allows for self-excitation and cross-excitation effects: the arrival at time $t$ of an event of type $j$ increases the intensity $\intensity_i(t+h)$ of events of type $i$ at time $t+h$ by an amount of $k_{ji}(h)$. A Hawkes process is thus a good candidate for a model of interactions between different types of events. Gaining popularity in the last decade, Hawkes processes have indeed found applications in many areas including earthquake modelling \citep{Ogata:1998aa, Turkyilmaz:2013aa, Fox:2016aa}, criminology \citep{Lewis:2012aa, Mohler:2013aa, Loeffler:2016:gunViolence}, social networks analysis \citep{Blundell:2012:Reciprocating:Relationships, Zhou:2013aa, Farajtabar:2017:fakeNews}, neurology \citep{Chornoboy:1988aa, chevalier:2015:microscopic, Gerhard:2017aa} and finance \citep{Bowsher:2007aa:HawkesModel, Bacry:2015aa, Jaisson:2016:roughLimits}.

In spite of their success and attractiveness, we notice that these Hawkes-process models account only for the dynamics of events and ignore the state of the underlying system they may influence. For instance, when applied to financial markets, Hawkes processes describe the arrival in time of buy and sell orders \citep{Large:2007aa:MeasureResiliency, bacry2016estimation, rambaldi:2016:volumeHawkes} but capture neither the asset price nor the supply and demand imbalance, which are impacted by the arriving orders. In fact, Hawkes processes can be contrasted with models based on continuous-time Markov chains where the focus is instead on a state process that represents the underlying system \citep{cont2010stochastic, cont2013:priceDynamics, huang:2015:QueueReactive, huang:2015:ergodicit}. However, because of the Markov property inherent to this second group of models, the arrival rates depend only on the current state and, thus, interactions like in general Hawkes processes are not possible. In effect, Hawkes processes and continuous-time Markov chains can have either an \emph{event} viewpoint or a \emph{state} viewpoint, respectively. Note that the dichotomy between Hawkes processes and Markovian models, and the need for something bridging the two, was already suggested by \citet[p.~1190--1191]{bacry2016estimation}. In a recent paper, \citet{Gonzales:2017:orderImpact} make a similar observation and propose a discrete-time Markov chain to model both the dependence on the most recent event and the current state of the system.

The original idea that started the present work is to address this gap by endowing the counting process $N$ with a state process $X$ that is fully coupled to $N$ in the following manner. On the one hand,  we wish to make the intensity dependent on both past events and states by changing \eqref{eq:linear_hawkes} to
\begin{equation*}
	\intensity_i (t) = \nu_i + \sum_{j=1}^{d}\int_{[0,t)} k_{ji}(t-s, X_s)dN_j(s),\quad t\geq 0,\, i=1,\ldots,d,
\end{equation*}
where the kernels now depend on the state process $X$. On the other hand, we imagine that each event in $N$ prompts a state change according to transition probabilities that depend on the type of the event. Naturally, we call this new process $(N, X)$ a \textit{state-dependent} Hawkes process.
The practical relevance and strong potential of this new model is demonstrated in \citet{Morariu:Pakkanen:2018}, where parametric estimation of the kernels from high-frequency financial data indeed reveals significant state dependence. 
In fact, this extension of Hawkes processes opens an avenue to novel models of high-frequency data that feature both excitation effects and a feedback loop between events and the state of the underlying system.

The first contribution of this paper is to turn this idea into a class of \textit{hybrid marked point processes} that encompasses and extends continuous-time Markov chains and Hawkes processes. To generalise state-dependent Hawkes processes to hybrid marked point processes, we view the process $(N,X)$ as a single marked point process on a product mark space and allow the intensity of events to be any measurable functional of past events and states. These new hybrid marked point processes are actually defined implicitly via their intensity that takes a specific product form. We prove that the dynamics generated by this product form completely characterise the class of hybrid marked point processes (Theorem  \ref{thm:implied_dynamics_characterisation}). Offering an \textit{event--state} viewpoint, this new class is well-suited to the joint modelling of events and the time evolution of the state of a system. While this general and flexible class provides a unifying framework for various existing processes, it also contains new processes with complex dynamics, as illustrated by state-dependent Hawkes processes.

The second contribution of this paper is to prove the strong existence and uniqueness of non-explosive state-dependent Hawkes processes and, more generally, hybrid marked point processes. In fact, by dispensing with a Lipschitz condition, we extend the results currently available in the literature for marked point processes defined via their intensity. It is known that a marked point process whose intensity $\intensity$ is expressed in terms of an intensity functional $\intensityFunctional$ can be formulated as a solution to a Poisson-driven stochastic differential equation (SDE) \citep{Massoulie:1998aa:StabilityResults, Bremaud:1996aa:StabilityNonLinearHawkes}. However, the existence and uniqueness results available in these works cannot be applied to hybrid marked point processes because their intensity functional may fail to satisfy the Lipschitz condition imposed therein.
We show that, under certain integrability or decay conditions, it is enough for $\psi$ to be dominated by either a Hawkes functional or an increasing function of the total number of past events in order to obtain the existence of a strong solution to the Poisson-driven SDE (Theorem \ref{thm:pathwise_existence}) and, in particular, the existence of hybrid marked point processes (Corollary \ref{cor:existence_hybrid}). The solution is constructed piece by piece along the time axis in a pathwise manner, taking advantage of the discrete nature of the driving Poisson random measure. A domination argument is then used to show non-explosiveness. In the context of multivariate point processes, a similar construction has already been considered in \citet{ccinlar2011probability}, while a similar domination argument is given in \citet{chevallier:2015:Hawkes_Generalised}. We combine the two in a more general setting (i.e., general mark space, initial conditions and intensity functional). We are also able to obtain strong and weak uniqueness without any specific assumptions (Theorems \ref{thm:uniqueness_pathwise} and \ref{thm:uniqueness_weak}) and, in particular, uniqueness of hybrid marked point processes (Corollary \ref{cor:uniqueness_hybrid}).

The paper is organised as follows. Section \ref{sec:framework_point_processes} introduces the framework that we use and presents the main results. Section \ref{sec:hybrid_marked_pp} proves the results concerning the dynamics of hybrid marked point processes. Section \ref{sec:pathwise_construction} proves the existence and uniqueness results. The \hyperref[sec:appendix]{Appendix} gathers technical results concerning the space $\boundedlyFiniteMeasures{\theCSMS}$ of boundedly finite integer-valued measures  and the enumeration representation of marked point processes.

\section{Framework and main results} \label{sec:framework_point_processes}

\subsection{A framework for point processes}

In the following, $\anyCSMS$ refers to a complete separable metric space and we denote by $\anyCSMSBorel$ its Borel $\sigma$-algebra. We reserve the notation $\markSpace$ for a complete separable metric space that represents the set of marks in the context of marked point processes. For most definitions, we follow closely \citet{daleyVereJonesVolume2} along with \citet{bremaud:1981:MartingaleDynamics}. The former reference will be especially used to introduce (marked) point processes while the latter is essential when defining the intensity process.

\subsubsection{Spaces of integer-valued measures.}

Let $\xi$ be a Borel measure on  $\anyCSMS$. We say that $\xi$ is boundedly finite if  $\xi(A)<\infty$ for every bounded Borel set $A\in\anyCSMSBorel$. We denote by $\integerValuedMeasures{\anyCSMS}$ the space of Borel measures on $\anyCSMS$ with values in $\integers\cup\{\infty\}$. We denote by $\boundedlyFiniteMeasures{\anyCSMS}$ the set of all $\xi\in\integerValuedMeasures{\anyCSMS}$ such that $\xi$ is boundedly finite. We denote by $\mppSpace$ the set of all $\xi\in\boundedlyFiniteMeasures{\realLine\times\markSpace}$ such that their ground measure $\xi_{g}(\cdot):=\xi(\cdot\times\markSpace)$ satisfies:
	\begin{enumerate}[label=\textup{(\roman*)}]
		\item $\xi_{g}\in\boundedlyFiniteMeasures{\realLine}$ \textup{;}
		\item $\xi_{g}(\{t\})=0 \mbox{ or } 1$ for all $t\in\realLine$ (we say that the ground measure is simple).
	\end{enumerate}
Observe that $\mppSpace\subset\boundedlyFiniteMeasures{\theCSMS}\subset\integerValuedMeasures{\theCSMS}$. The space $\integerValuedMeasures{\theCSMS}$ corresponds to the realisations of potentially explosive point processes, while the space $\boundedlyFiniteMeasures{\theCSMS}$ corresponds to the realisations of non-explosive point processes and contains all the realisations of potentially explosive marked point processes. Regarding the space $\mppSpace$,
each $\xi\in\mppSpace$ is a realisation of a non-explosive marked point process. When $\xi(\{(t,m)\})=1$ for some $t\in\realLine$ and $m\in\markSpace$, this should be interpreted as an event happening at time $t$ with characteristics $m$. The boundedly finite property of the ground measure ensures that, in any finite amount of time, only finitely many events can occur (i.e., the marked point point process is non-explosive). The simpleness constraint on the ground measure means that there cannot be two events at the same time.

The so-called $\weakHash$-distance $\weakHashDistancePlain$ (``weak-hash'') introduced by \citet[p.~403]{daleyVereJonesVolume1}, makes $\boundedlyFiniteMeasures{\anyCSMS}$ a complete separable metric space, see Theorem A2.6.III in \citet[p.~404]{daleyVereJonesVolume1}. The corresponding $\sigma$-algebra $\boundedlyFiniteMeasuresBorel{\anyCSMS}$  coincides with the one generated by all mappings $\xi\mapsto\xi(A)$, $\xi\in\boundedlyFiniteMeasures{\anyCSMS}$, $A\in\borelSetsOf{\anyCSMS}$.
Proposition A2.6.II of \citet[p.~403]{daleyVereJonesVolume1} characterises convergence in this topology, called the $\weakHash$-topology. These properties of the space $\boundedlyFiniteMeasures{\anyCSMS}$ play an important role in this work. Note that \citet{morariu:2017:weakHashMetric} clarifies the proofs of Proposition A2.6.II and Theorem A2.6.III of \citet{daleyVereJonesVolume1}. Indeed, the original proofs assume a certain function to be monotonic. As this does not seem to actually hold, \citet{morariu:2017:weakHashMetric} proposes alternative arguments where required. Besides, in the \hyperref[sec:appendix]{Appendix}, we show that $\mppSpace$ is indeed a Borel set of $\boundedlyFiniteMeasures{\theCSMS}$, see Lemma \ref{prop:mppSpace_is_measurable}.

Finally, for any $u\in\anyCSMS$, we denote by $\delta_u$ the Dirac measure at $u$.

\subsubsection{Non-explosive marked point processes.}

In the following, the notation $\probSpace$ refers to a probability space. For any $\sigma$-algebra $\mathcal{S}$, the trace of $A\in\mathcal{S}$ on $\mathcal{S}$ is defined by $A\cap\mathcal{S}:=\{A\cap S\,:\,S\in\mathcal{S}\}$.
\theoremstyle{definition}
\newtheorem{marked_point_process}[]{Definition}[section]
\begin{marked_point_process}[Non-explosive point process] \label{def:point_process}
A \textit{non-explosive point process} on  $\anyCSMS$ is a measurable mapping from $\measureSpace$ into $\boundedlyFiniteMeasuresMeasurableSpace{\anyCSMS}$.
\end{marked_point_process}
\newtheorem{marked_point_process2}[marked_point_process]{Definition}
\begin{marked_point_process2}[Non-explosive marked point process] \label{def:marked_point_process}
A \textit{non-explosive marked point process} $N$ on $\theCSMS$ is a non-explosive point process $N$ on $\realLine\times\markSpace$ such that $N(\omega)\in\mppSpace$ for all $\omega\in\realisationsSpace$.
\end{marked_point_process2}
\newtheorem{remark_marked_point_process2}[marked_point_process]{Remark}
\begin{remark_marked_point_process2} \label{rem:remark_marked_point_process2}
By applying Lemma 1.6 in \citet[p.~4]{kallenberg2006foundations}, we obtain that $\mppSpaceBorel=\mppSpace\cap\boundedlyFiniteMeasuresBorel{\theCSMS}$, where $\mppSpace$ is also equipped with the $\weakHash$-metric $\weakHashDistancePlain$. This implies that Definition \ref{def:marked_point_process} is equivalent to saying that a non-explosive marked point process is a measurable mapping from $\measureSpace$ into $\mppMeasureSpace$.
\end{remark_marked_point_process2}
Next, a non-explosive point process induces a probability measure on $\boundedlyFiniteMeasures{\anyCSMS}$.
\newtheorem{induced_probability}[marked_point_process]{Definition}
\begin{induced_probability}[Induced probability] \label{def:induced_probability}
Let $N$ be a non-explosive point process on $\anyCSMS$. We define the \textit{induced probability measure}	$\inducedProb{N}$ on the measurable space $\boundedlyFiniteMeasuresMeasurableSpace{\anyCSMS}$ through the relation
\begin{equation*}
	\inducedProb{N}(A):=\prob\left( N^{-1}(A) \right), \quad A\in\boundedlyFiniteMeasuresBorel{\anyCSMS}.
\end{equation*}
\end{induced_probability}

\subsubsection{Enumeration representation.}

It is common to define instead marked point processes on $\realLineNonNegative\times\markSpace$ as a sequence $(T_{n},M_{n})_{n\in\integers}$ of random variables in $(0,\infty]\times\markSpace$ such that $(T_{n})_{n\in\integers}$ is non-decreasing and $T_{n}<\infty$ implies $T_{n}<T_{n+1}$ \citep{jacod:1975:projection, bremaud:1981:MartingaleDynamics}. We will call such a sequence an enumeration.
Here, $T_{n}$ is to be interpreted as the time when the $n^{th}$ event occurs while $M_{n}$ describes the characteristics of that event. Moreover, $T_{n}<\infty$ with $T_{n+1}=\infty$ means that there are no more events after time $T_{n}$. Note that this definition allows for explosion in the sense that $\lim_{n\rightarrow\infty}T_{n}<\infty$ is possible with positive probability. This is why we stress the non-explosive character of marked point processes in Definition \ref{def:marked_point_process}. There is a one-to-one correspondence between non-explosive marked point processes on $\realLineNonNegative\times\markSpace$ and enumerations such that $\lim_{n\rightarrow\infty}T_{n}=\infty$ $\as$ We give a proof of this correspondence in Appendix \ref{subsec:enumeration_representation} for completeness.

\subsubsection{Poisson processes.}

Let $\nu$ be a boundedly finite measure on $(\anyCSMS,\anyCSMSBorel)$. We say that a non-explosive point process $N$ on $\anyCSMS$ is a Poisson process on $\anyCSMS$ with parameter measure $\nu$ if $N(A_1),\ldots,N(A_n)$ are mutually independent for all disjoint and bounded sets $A_1,\ldots,A_n\in\borelSetsOf{\anyCSMS}$, $n\in\integers$, and $N(A)$ follows a Poisson distribution with parameter $\nu(A)$ for all bounded sets $A\in\anyCSMSBorel$. Their existence can be verified using Theorem 9.2.X in \citet[p.~30]{daleyVereJonesVolume2}, see Example 9.2(b) on p.~31 therein.

\subsubsection{Pathwise integration.} \label{subsec:integration}

Let $N$ be a non-explosive point process on $\anyCSMS$.
Let $H:\realisationsSpace\times\anyCSMS\rightarrow\realLineNonNegative$ be an $\sigmaAlgebra\otimes\anyCSMSBorel$-measurable non-negative mapping. In particular, $H$ is an $\realLineNonNegative$-valued stochastic process on $\anyCSMS$.
One can define the integral of $H$ against $N$ in a pathwise fashion as
\begin{equation*}
	I(\omega):=\int_{\anyCSMS}H(\omega,u)N(\omega,du),\quad \omega\in\realisationsSpace.
\end{equation*}
Besides, by a monotone class argument, one can check that $\omega\mapsto I(\omega)$ is $\sigmaAlgebra$-measurable. In the special case where $N$ is actually a non-explosive marked point process on $\realLineNonNegative\times\markSpace$, the integral can be rewritten as
\begin{equation*}
	\iint_{\realLineNonNegative\times\markSpace} H(t,m)N(dt,dm) = \sum_{n\in\integers}H(T_{n},M_{n})\indicator_{\{T_{n}<\infty\}},\quad \as,
\end{equation*}
where $(T_{n},M_{n})_{n\in\integers}$ is the enumeration corresponding to $N$. For any $\xi\in\mppSpace$ and any $\tau\in\realLine$ such that $\xi(\{\tau\}\times\markSpace)>0$, we abuse the notation and define $\iint_{\{\tau\}\times\markSpace} m \xi(dt,dm)$ as the unique element $m\in\markSpace$ such that $\xi(\{\tau\}\times\{m\})=1$.

\subsubsection{Shifts, restrictions, histories and predictability.} \label{subsec:shift_histories}

For all $t\in\realLine$, define the shift operator $\shift{t}:\boundedlyFiniteMeasures{\realLine\times\anyCSMS}\rightarrow\boundedlyFiniteMeasures{\realLine\times\anyCSMS}$ by
	$\shift{t}\xi(A):=\xi(A+t)$,  $A\in\borelSetsOf{\realLine\times\anyCSMS}$,
where $A+t:=\{(s+t,u)\in\realLine\times\anyCSMS \,:\, (s,u)\in A\}$. Then, for any non-explosive point process $N$ on $\realLine\times\anyCSMS$, define $\shift{t}N$ through $(\shift{t}N)(\omega):=\shift{t}(N(\omega))$, $\omega\in\realisationsSpace$. It will be useful to show that $\shift{t}\xi$ is jointly continuous in $t$ and $\xi$ (Lemma \ref{prop:shift_jointly_continuous}).

Denote the restriction to the negative real line of any realisation $\xi\in\boundedlyFiniteMeasures{\realLine\times\anyCSMS}$ by $\restrictionToNegativeOfMeasure{\xi}$, which is defined by
	$\restrictionToNegativeOfMeasure{\xi}(A):=\xi(A\cap\realLineNegative\times\anyCSMS)$,  $A\in\borelSetsOf{\realLine\times\anyCSMS}$.
We can then define the restriction to the negative real line of any non-explosive point process $N$ on $\realLine\times\anyCSMS$ by $\restrictionToNegativeOfMeasure{N}(\omega):=\restrictionToNegativeOfMeasure{(N(\omega))}$, $\omega\in\realisationsSpace$.  Similarly, define the notations  $\restrictionToNonPositiveOfMeasure{\xi}(A):=\xi(A\cap\realLineNonPositive\times\anyCSMS)$, $\restrictionToNonPositiveOfMeasure{N}(\omega):=\restrictionToNonPositiveOfMeasure{N(\omega)}$, $\restrictionToNonNegativeOfMeasure{\xi}(A):=\xi(A\cap\realLineNonNegative\times\anyCSMS)$, $\restrictionToNonNegativeOfMeasure{N}(\omega):=\restrictionToNonNegativeOfMeasure{N(\omega)}$, $\restrictionToPositiveOfMeasure{\xi}(A):=\xi(A\cap\realLinePositive\times\anyCSMS)$ and $\restrictionToPositiveOfMeasure{N}(\omega):=\restrictionToPositiveOfMeasure{N(\omega)}$.

These notations will allow us to refer to the internal history of $N$. For instance, for any $t\in\realLine$, $\restrictionToNegativeOfMeasure{(\shift{t}N)}$ contains the history of the process up to time $t$, excluding time $t$. To lighten these notations, we will use the conventions $\pointProcessAtTime{\xi}{t}:= \restrictionToNegativeOfMeasure{(\shift{t}\xi)}$, $\restrictionToNonPositiveOfMeasure{\shift{t}\xi} := \restrictionToNonPositiveOfMeasure{(\shift{t}\xi)}$, $\restrictionToPositiveOfMeasure{\shift{t}\xi} := \restrictionToPositiveOfMeasure{(\shift{t}\xi)}$ and $\restrictionToNonNegativeOfMeasure{\shift{t}\xi} := \restrictionToNonNegativeOfMeasure{(\shift{t}\xi)}$. It will be useful to note that these restriction mappings are measurable (Lemma \ref{prop:restrictions_are_measurable}) and that $\pointProcessAtTime{\xi}{t}$ is left-continuous as a function of $t\in\realLine$ (Lemma \ref{prop:history_left_continuous}).

Let $N$ be a non-explosive point process on $\realLine\times\anyCSMS$. We can define the filtration $\internalHistoryCollection{N}= \collection{\internalHistory{N}{t}}{t\in\realLine}$ that corresponds to the internal history of $N$ by
\begin{equation*}
	\internalHistory{N}{t}:=\sigma\left\{N(A\times U) \, : \, A\in\borelSetsOf{\realLine},\,  A\subset(-\infty,t],\, U\in\anyCSMSBorel \right\} \quad \mbox{for all } t\in\realLine.
\end{equation*}
Using Lemma 1.4 in \citet[p.~4]{kallenberg2006foundations} along with the characterisation of $\boundedlyFiniteMeasuresBorel{\realLine\times\anyCSMS}$ given in Theorem A2.6.III in \citet[p.~404]{daleyVereJonesVolume1}, one can check that
	$\internalHistory{N}{t} = \sigma(\shift{t}\restrictionToNonPositiveOfMeasure{N})$.

In the following, we call a \emph{history} any filtration that contains the internal history of $N$, that is any filtration  $\aHistoryCollection = \collection{\aHistory{t}}{t\in\realLine}$ such that
	$\internalHistory{N}{t} \subset\aHistory{t}$, $t\in\realLine$.
Equivalently, one says that $N$ is $\aHistoryCollection$-adapted. The notation $\aHistoryCollection= \collection{\aHistory{t}}{t\in\realLine}$ will always be used to refer to a history.
We also need to define the predictable $\sigma$-algebra $\predictableSigmaAlgebra$ on $\realisationsSpace\times\realLine\times\anyCSMS$ corresponding to a history $\aHistoryCollection$. The $\sigma$-algebra $\predictableSigmaAlgebra$ is the one which is generated by all the sets of the form
\begin{equation*}
	A\times (s,t] \times U, \quad  s,t\in\realLine,\, s< t, \, U\in\anyCSMSBorel ,\, A\in\aHistory{s}.
\end{equation*}
Any mapping $H:\realisationsSpace\times\realLine\times\anyCSMS\rightarrow\realLine$ that is $\predictableSigmaAlgebra$-measurable is called an $\aHistoryCollection$-predictable process. Any mapping $H:\realisationsSpace\times\realLinePositive\times\anyCSMS\rightarrow\realLine$ that is $(\realisationsSpace\times\realLinePositive\times\anyCSMS)\cap \predictableSigmaAlgebra$-measurable is also called an $\aHistoryCollection$-predictable process. Given an $\aHistoryCollection$-stopping time $\tau$, the strict past $\aHistory{\tau-}$ is defined as the $\sigma$-algebra generated by the all the classes $\{t<\tau\}\cap \aHistory{t}$, $t\in\realLine$.

\subsubsection{Intensity process and functional.} \label{subsec:intensity}

We equip the mark space $(\markSpace,\markSpaceBorel)$ with a reference measure $\markSpaceMeasure$, allowing us to define the concept of intensity rigorously.
Let $N$ be a marked point process on $\theCSMS$ and $\aHistoryCollection= \collection{\aHistory{t}}{t\in\realLine}$ a history. Let $\intensity:\realisationsSpace\times\realLinePositive\times\markSpace\rightarrow\realLineNonNegative$ be a non-negative $\aHistoryCollection$-predictable process. We say that $\intensity$ is the $\aHistoryCollection$-intensity of $N$ relative to $\markSpaceMeasure$ if for every non-negative $\aHistoryCollection$-predictable process $H:\realisationsSpace\times\realLinePositive\times\markSpace\rightarrow\realLineNonNegative$,
\begin{equation} \label{eq:intensity_definition}
	\E{\iint_{\realLinePositive\times\markSpace}H(t,m)N(dt,dm)} = \E{\iint_{\realLinePositive\times\markSpace}H(t,m)\intensityAtTimeAtMarkWithIndex{t}{m}{}\markSpaceMeasure(dm)dt}.
\end{equation}
Note that if an intensity exists, it is then unique up to $\prob(d\omega)\markSpaceMeasure(dm)dt$-null sets thanks to the predictability requirement, see \citet[Section II.4]{bremaud:1981:MartingaleDynamics} and \citet[p.~391]{daleyVereJonesVolume2}.
In this paper, we will be particularly interested in intensities that are expressed in terms of a functional applied to the point process.
\newtheorem{intensity_functional}[marked_point_process]{Definition}
\begin{intensity_functional}[Intensity functional] \label{def:intensity_functional}
	Let $\intensityFunctional:\markSpace\times\boundedlyFiniteMeasures{\theCSMS}\rightarrow\realLineNonNegative\cup\{\infty\}$ be a measurable functional. We say that a non-explosive marked point process $N:\realisationsSpace\rightarrow\mppSpace$ admits $\intensityFunctional$ as its \textit{intensity functiona}l if $N$ admits an $\internalHistoryCollection{N}$-intensity $\intensity:\realisationsSpace\times\realLinePositive\times\markSpace\rightarrow\realLineNonNegative$ relative to $\markSpaceMeasure$ such that
	\begin{equation} \label{eq:intensity_functional_definition}
		\intensityAtTimeAtMarkWithIndex{\omega,t}{m}{} = \intensityFunctionalAtMarkGivenPastWithIndex{m}{\pointProcessAtTime{N(\omega)}{t}}{},\quad \prob(d\omega)dt\markSpaceMeasure(dm)\aeverywhere
	\end{equation}
\end{intensity_functional}

\subsubsection{Initial condition.} \label{subsec:initial_condition}

Let $\probSpaceInitialCondition$ be a given probability space and $N_{\leq 0}$ be a given marked point process on $\theCSMS$ such that $\restrictionToNonPositiveOfMeasure{N_{\leq 0}(\omega_{\leq 0})}=N_{\leq 0}(\omega_{\leq 0})$ for all $\omega_{\leq 0}\in\realisationsSpaceInitialCondition$ (i.e., there are no events on $\realLinePositive$). We will reserve the notation $N_{\leq 0}$ to refer to an initial condition.

Let $\probSpaceFuture$ be another probability space that will correspond to the driving Poisson process in the SDE introduced in Subsection \ref{subsec:existence_uniqueness} below.
In the context of strong existence, we will work with the probability space $\probSpace$ which we define as the completion \citep[p.~13]{kallenberg2006foundations} of the product probability space given by
\begin{equation} \label{eq:probability_space}
	\realisationsSpace := \realisationsSpaceInitialCondition\times\realisationsSpaceFuture, \quad \tilde{\sigmaAlgebra}:=\sigmaAlgebraInitialCondition\otimes\sigmaAlgebraFuture, \quad \tilde{\prob}:=\probInitialCondition\times\probFuture.
\end{equation}
Such a structure of the probability space is motivated by the fact that the driving noise and the initial condition are independent.
\newtheorem{satisfying_initial_condition_strong}[marked_point_process]{Definition}
\begin{satisfying_initial_condition_strong}[Strong initial condition] \label{def:satisfying_initial_condition_strong}
	Let $N:\realisationsSpace\rightarrow\mppSpace$ be a non-explosive marked point process on $\theCSMS$. We say that $N$ \textit{satisfies a strong initial condition} $N_{\leq 0}$ if $\restrictionToNonPositiveOfMeasure{N(\omega)}=N_{\leq 0}(\omega_{\leq 0})$  $\as$,
	where $\omega=(\omega_{\leq 0},\omega_{>0})\in\realisationsSpace$.
\end{satisfying_initial_condition_strong}

In the context of weak uniqueness, one needs a another concept of initial condition. Let $(\Omega',\mathcal{F}',\mathbb{P}')$ be another probability space potentially different from $\probSpace$.
\newtheorem{satisfying_initial_condition_weak}[marked_point_process]{Definition}
\begin{satisfying_initial_condition_weak}[Weak initial condition] \label{def:satisfying_initial_condition_weak}
	Let $N':\Omega'\rightarrow\mppSpace$ be a non-explosive marked point process on $\theCSMS$. We say that $N'$ \textit{satisfies a weak initial condition} $N_{\leq 0}$ if the induced probability $\inducedProb{\restrictionToNonPositiveOfMeasure{N'}}$ coincides with $\inducedProb{N_{\leq 0}}$.
\end{satisfying_initial_condition_weak}

\subsection{Hybrid marked point processes: an event--state viewpoint} \label{subsec:hybrid_marked_point_processes}

\subsubsection{Mark space and state process.}

Let $\eventSpaceMeasurable$ and $\stateSpaceMeasurable$ be two measure spaces where both $\eventSpace$ and $\stateSpace$ are complete separable metric spaces and both $\eventSpaceMeasure$ and $\stateSpaceMeasure$ are a boundedly finite Borel measures. Each $e\in\eventSpace$ represents a type of event and we call $\eventSpace$ the event space. Each $x\in\stateSpace$ represents a possible state of a system and we call $\stateSpace$ the state space. Motivated by the need to jointly model events and the state of the system (see \hyperref[sec:intro]{Introduction}), we consider a mark space $\markSpaceMeasurable$ of the form
\begin{equation} \label{eq:product_mark_space}
	\markSpace:=\eventSpace\times\stateSpace, \quad \markSpaceBorel=\borelSetsOf{\eventSpace}\otimes\borelSetsOf{\stateSpace}, \quad \markSpaceMeasure:=\eventSpaceMeasure\times\stateSpaceMeasure.
\end{equation}
Such a decomposition of the mark space admits the following interpretation. Let $\xi\in\mppSpace$ be a realisation of a marked point process on $\theCSMS$.  A point $t\in\realLine$ and a point $m=(e,x)\in\markSpace$ such that $\xi(\{t,m\})=1$ can now be interpreted as an event of type $e$ occurring at time $t$ and moving the state of the system to $x$. To formalise this viewpoint, we define the state functional and the state process as follows.
\newtheorem{state_process}[marked_point_process]{Definition}
\begin{state_process}[State functional and state process] \label{def:state_process}
We define the measurable state functional \\ $F:\boundedlyFiniteMeasures{\theCSMS}\rightarrow\stateSpace$ by
\begin{equation*}
F(\xi):=
\begin{cases}
	\iint_{\{\kappa(\xi)\}\times\markSpace} x \xi(dt,de,dx), \quad &\mbox{if } \xi\in\mppSpace\mbox{ and } \kappa(\xi)>-\infty, \\
	x_0, \quad &\mbox{otherwise,}\\
\end{cases}
\end{equation*}
where $\kappa(\xi):=\inf\{t<0\,:\,\xi((t,0)\times\markSpace)=0\}$ and $x_0\in\stateSpace$ is an arbitrary initial state. Given a non-explosive marked point process $N$ on $\markSpace$, we define the state process $\collection{\stateProcess_t}{t\in\realLine}$ by
\begin{equation*}
	\stateProcess_t := F(\pointProcessAtTime{N}{t}),\quad t\in\realLine.
\end{equation*}
\end{state_process}
Note that $\kappa(\pointProcessAtTime{N}{t})$ is the time of the last event up to time $t$ and, thus, $\stateProcess_{t}$ is the coordinate $x\in\stateSpace$ of the mark $m=(e,x)\in\markSpace$ of the most recent event. As a consequence, we indeed have that a point $t\in\realLine$ and a point $m=(e,x)\in\markSpace$ such that $N(\{t,m\})=1$ can be interpreted as an event of type $e$ occurring at time $t$ and moving the state of the system to $x$.
With this viewpoint, a marked point process on $\theCSMS$ allows to jointly model the evolution of a system with state process $\collection{\stateProcess_t}{t\in\realLine}$ in $\stateSpace$ and the arrival in time of the event types $\eventSpace$. To check that the state functional $F$ is indeed measurable, one can adapt the proof of Lemma \ref{prop:mpp_is_enumeration} in the \hyperref[sec:appendix]{Appendix}.

\subsubsection{Implicit definition through the intensity, implied dynamics and characterisation.}

We can now introduce the class of hybrid marked point processes, which provides a unified treatment of various types of processes, including continuous-time Markov chains and Hawkes processes, whence the qualifier \textit{hybrid}. Moreover, this class contains new types of processes, such as state-dependent Hawkes processes (Example \ref{ex:generalised_hawkes_processes} below), which are applicable to the joint modelling of events and systems. This class is specified implicitly through a specific form of the intensity. We were inspired by \citet[p.~13]{cartea2015:enhancing} who, in the context of a continuous-time Markov chain model, propose a decomposition of the intensity that is similar in spirit. We should also mention the connection to the decomposition of the rate kernel $\alpha$ of a continuous-time Markov chain into a rate function $c$ and transition kernel $\mu$ (i.e., $\alpha=\mu c$), see \citet[p.~238-239]{kallenberg2006foundations}, even though this is slightly different as $c$ is the total intensity and does not depend on the event variable $e\in\eventSpace$.
\newtheorem{hybrid_marked_pp}[marked_point_process]{Definition}
\begin{hybrid_marked_pp}[Hybrid marked point processes] \label{def:hybrid_marked_pp}
Let $\stateFunctional:\stateSpace\times\eventSpace\times\stateSpace\rightarrow\realLineNonNegative$ be a measurable non-negative function such that $\stateFunctionalAtStateGivenEventGivenPast{\cdot}{e}{x}$ is a probability density over $\stateSpaceMeasurable$ for all $e\in\eventSpace$, $x\in\stateSpace$. Let $\eventFunctional:\eventSpace\times\boundedlyFiniteMeasures{\theCSMS} \rightarrow \realLineNonNegative\cup\{\infty\}$ be a measurable non-negative functional. Define the measurable intensity functional $\intensityFunctional:\markSpace\times\boundedlyFiniteMeasures{\theCSMS}\rightarrow\realLineNonNegative\cup\{\infty\}$ by $\intensityFunctionalAtMarkGivenPastWithIndex{m}{\xi}{}:=\stateFunctionalAtStateGivenEventGivenPast{x}{e}{F(\xi)}\eventFunctionalAtEventGivenPast{e}{\xi}$ for all $m=(e,x)\in\markSpace$, $\xi\in\boundedlyFiniteMeasures{\theCSMS}$. A \emph{hybrid marked point process} with transition function $\stateFunctional$ and event functional $\eventFunctional$ is a non-explosive marked point process $N:\realisationsSpace\rightarrow\mppSpace$ that admits $\intensityFunctional$ as its intensity functional. In other words, $N$ admits an $\internalHistoryCollection{N}$-intensity $\intensity: \realisationsSpace\times\realLinePositive\times\markSpace\rightarrow\realLineNonNegative$ relative to $\markSpaceMeasure$ that satisfies
\begin{equation} \label{eq:intensity_form_canonical_example}
	\intensityAtTimeAtMarkWithIndex{\omega,t}{e,x}{}=\stateFunctionalAtStateGivenEventGivenPast{x}{e}{\stateProcess_{t}(\omega)} \eventFunctionalAtEventGivenPast{e}{\pointProcessAtTime{N(\omega)}{t}}, \quad \prob(d\omega)dt\markSpaceMeasure(de,dx)\aeverywhere
\end{equation}
\end{hybrid_marked_pp}

To demonstrate the generality and flexibility of hybrid marked point processes, we give three examples of well-known processes that belong to the class. While at first these examples might be understood only at an intuitive level, the reader should become fully convinced of their validity once Theorem \ref{thm:implied_dynamics_characterisation} is introduced.
\newtheorem{compound_poisson}[marked_point_process]{Example}
\begin{compound_poisson}[Compound Poisson process] \label{ex:compound_poisson}
  Let $\eventSpace=\{0\}$ (i.e., just one type of event), $\eventSpaceMeasure=\delta_{0}$, $\stateSpace=\realLine$, $\stateSpaceMeasure(dx)=dx$ (i.e., the Lebesgue measure) and suppose that $f:\realLine\rightarrow\realLineNonNegative$ is a probability density function.
  Consider a hybrid marked point process $N$ with constant event functional $\eventFunctional\equiv \nu \in\realLinePositive$ and transition function given by $\stateFunctional(x'\,|\,x)= f(x'-x)$, $x',x\in\realLine$. Then, the $\internalHistoryCollection{N}$-intensity of $N$ satisfies $\intensity(t,x)=f(x-\stateProcess_t)\nu$ and the state process $\collection{X_t}{t\in\realLineNonNegative}$ is a compound Poisson process with rate $\nu$ and jump size distribution $f(x)dx$.
\end{compound_poisson}
\newtheorem{markov_chain}[marked_point_process]{Example}
\begin{markov_chain}[Continuous-time Markov chain] \label{ex:markov_chain}
  Let $\eventSpace$, $\eventSpaceMeasure$, $\stateSpace$ and $\stateSpaceMeasure$ be as in Example \ref{ex:compound_poisson}. Suppose that a hybrid marked point process $N$ has event functional of the form $\eventFunctional(\xi)=c(F(\xi))$, $\xi\in\boundedlyFiniteMeasures{\realLine^2}$, where $c$ is a positive function. The $\internalHistoryCollection{N}$-intensity of $N$ is then given by $\intensity(t,x)=\stateFunctional(x\,|\,X_t)c(X_t)$ and the state process $\collection{X_t}{t\in\realLineNonNegative}$ is a continuous-time Markov chain with rate function $c$ and transition kernel $\mu(x,B)=\int_B\stateFunctional(y\,|\,x)dy$, $x\in\realLine$, $B\in\borelSetsOf{\realLine}$.
\end{markov_chain}
\newtheorem{multivariate_Hawkes_process}[marked_point_process]{Example}
\begin{multivariate_Hawkes_process}[Multivariate Hawkes process] \label{ex:multivariate_Hawkes_process}
  Let $\eventSpace=\{1,\ldots,d\}$, $d\in\integers$, $\eventSpaceMeasure=\sum_{n=1}^{d}\delta_n$, $\stateSpace=\{0\}$ (i.e., only one possible state), $\stateSpaceMeasure=\delta_0$, $\nu=(\nu_1,\ldots,\nu_d)\in\realLinePositive^d$ and $k:\realLinePositive\times\eventSpace^2\rightarrow\realLineNonNegative$. Consider a hybrid marked point process $N$ with event functional
\begin{equation*}
  \eventFunctionalAtEventGivenPast{e}{\xi}= \nu_e + \iint_{(-\infty,0)\times\eventSpace}k(-t',e',e)\xi(dt',de'),\quad e=1,\ldots,d,\, \xi\in\boundedlyFiniteMeasures{\realLine\times\eventSpace},
\end{equation*}
  and note that the transition function must satisfy $\stateFunctional\equiv 1$. Then $N$ is a multivariate Hawkes process with base rate $\nu$ and kernel $k$.
\end{multivariate_Hawkes_process}

\theoremstyle{plain}

Let us next explain the dynamics that the intensity \eqref{eq:intensity_form_canonical_example} implies. If one integrates out the state variable $x$ by using the fact that $\stateFunctionalAtStateGivenEventGivenPast{\cdot}{e}{\stateProcess_{t}}$ is a probability density, one can see that the $\internalHistoryCollection{N}$-intensity of the marked point process $N(\cdot\times\stateSpace)$ on $\realLine\times\eventSpace$ is exactly $\eventFunctionalAtEventGivenPast{e}{\shift{t	}\restrictionToNegativeOfMeasure{N}}$. In other words, $\eventFunctionalAtEventGivenPast{e}{\shift{t	}\restrictionToNegativeOfMeasure{N}}$ is the intensity of the aggregation of events of type $e$ (irrespectively of how they impact the state process $\stateProcess_{t}$). Then, $\eventFunctionalAtEventGivenPast{e}{\shift{t	}\restrictionToNegativeOfMeasure{N}}$ is distributed in the state space $\stateSpace$ according to $\stateFunctionalAtStateGivenEventGivenPast{x}{e}{\stateProcess_{t}}$, specifying the intensity of events with mark $(e,x)$. This suggests that $\stateFunctionalAtStateGivenEventGivenPast{\cdot}{e}{\stateProcess_{t}}$ is the probability density of the next state of the system given that the next event is of type $e$ and that the current state is $\stateProcess_{t}$. This intuition is confirmed by Theorem \ref{thm:implied_dynamics_characterisation}, which actually goes further and states that these dynamics characterise hybrid marked point processes. The proof is presented in Subsection \ref{subsec:implied_dynamics_characterisation}.
\newtheorem{implied_dynamics_characterisation}[marked_point_process]{Theorem}
\begin{implied_dynamics_characterisation}[Implied dynamics and characterisation] \label{thm:implied_dynamics_characterisation}
Let $\stateFunctional$ and $\eventFunctional$ be as in Definition \ref{def:hybrid_marked_pp}. Moreover, suppose that $N$ is a non-explosive marked point process on $\theCSMS$ with an $\internalHistoryCollection{N}$-intensity relative to $\markSpaceMeasure$. Then, $N$ is a hybrid marked point process with transition function $\stateFunctional$ and event functional $\eventFunctional$ if and only if the following two statements hold.
\begin{enumerate}[label=\textup{(\roman*)}]
			\item $N_\eventSpace(\cdot):=N(\cdot\times\stateSpace)$ is a non-explosive marked point process on $\realLine\times\eventSpace$ that admits an $\internalHistoryCollection{N}$-intensity $\intensity_\eventSpace:\realisationsSpace\times\realLinePositive\times\eventSpace\rightarrow\realLineNonNegative$ relative to  $\eventSpaceMeasure$ such that
			$\intensityAtTimeAtMarkWithIndex{\omega, t}{e}{\eventSpace}=\eventFunctionalAtEventGivenPast{e}{\pointProcessAtTime{N(\omega)}{t}}$ holds  $\prob(d\omega)dt\eventSpaceMeasure(de)\aeverywhere$
			\item Let $t\in\realLineNonNegative$ and define the stopping time $\tau_t:=\sup\{ u > t \,:\, N((t,u)\times\markSpace) = 0 \}$ and the random elements $(E,X):=\iint_{\{\tau_t\}\times\markSpace} (e,x) N(du,de,dx)$ such that $\tau_t$ is the time of the first event after time $t$ and $(E,X)$ is the corresponding mark. We have that
	\begin{equation} \label{eq:probability_of_next_state}
		\prob\left( X\in dx \,|\, \sigma(E)\vee \internalHistory{N}{\tau_t-} \right)\indicator_{\{\tau_t<\infty\}}
		=\stateFunctionalAtStateGivenEventGivenPast{x}{E}{\stateProcess_{t}}\stateSpaceMeasure(dx)\indicator_{\{\tau_t<\infty\}},\quad \as
	\end{equation}
		\end{enumerate}
\end{implied_dynamics_characterisation}
\theoremstyle{definition}
\newtheorem{transition_probability_remark}[marked_point_process]{Remark}
\begin{transition_probability_remark} \label{rem:transition_probability_remark}
As shown in the proof of Theorem \ref{thm:implied_dynamics_characterisation}, Equation \eqref{eq:probability_of_next_state} implies that
\begin{equation*}
	\prob\left( X\in dx \,|\, \sigma(E)\vee \internalHistory{N}{t},\{\tau_t<\infty\} \right)
		=\stateFunctionalAtStateGivenEventGivenPast{x}{E}{\stateProcess_{t}}\stateSpaceMeasure(dx)\quad \as
\end{equation*}

\end{transition_probability_remark}

We now add a fourth example to show that Definition \ref{def:hybrid_marked_pp} contains also new types of processes. This example extends Hawkes processes to what could be called state-dependent Hawkes processes. Together, the four examples demonstrate that hybrid marked point processes provide a common framework to construct and analyse various types of processes.
\theoremstyle{definition}
\newtheorem{generalised_hawkes_processes}[marked_point_process]{Example}
\begin{generalised_hawkes_processes}[State-dependent Hawkes process] \label{ex:generalised_hawkes_processes}
Consider hybrid marked point processes with event functionals $\eventFunctional$ of the form
\begin{equation} \label{eq:event_functional_hawkes}
	\eventFunctionalAtEventGivenPast{e}{\xi}=\nu(e) + \iint_{(-\infty,0)\times\markSpace}k(-t',m',e)\xi(dt',dm'),\quad e\in\eventSpace, \xi\in\boundedlyFiniteMeasures{\theCSMS},
\end{equation}
where $\nu:\eventSpace\rightarrow\realLineNonNegative$ and $k:\realLine\times\markSpace\times\eventSpace\rightarrow\realLineNonNegative$ are non-negative measurable functions.
We show that such functionals are indeed measurable (Proposition \ref{prop:hawkes_functionals}).
By Theorem \ref{thm:implied_dynamics_characterisation}, this gives rise to a marked point process $N_\eventSpace$ with marks in $\eventSpace$ and intensity $\eventFunctional$ that interacts with a state process $\collection{\stateProcess_t}{t\in\realLine}$ on $\stateSpace$ with transition probabilities $\stateFunctional$. On the one hand, events in $N_\eventSpace$ occur like in a Hawkes process except that now the kernel depends also on the state process. For example, an event of type $e'\in\eventSpace$ might precipitate an event of type $e\in\eventSpace$ only if it moves the system to some specific state $x_0\in\stateSpace$, i.e., $k(\cdot, e',x',e)\equiv 0$ as soon as $x'\neq x_0$. On the other hand, the occurence of an event in $N_\eventSpace$ prompts a state change according to the transition probabilities $\stateFunctional$. Consequently, such a marked point process defines a state-dependent Hawkes process where the state process is fully coupled with the Hawkes process. Viewing $N_\eventSpace$ and $\collection{\stateProcess_t}{t\in\realLine}$ as one single marked point process $N$ on $\eventSpace\times\stateSpace$ with intensity $\stateFunctional\eventFunctional$ will allow us to prove the existence of such dynamics, see Corollary \ref{cor:existence_hybrid} and Example \ref{ex:existence_state_hawkes}.

This subclass of hybrid marked point processes extends the regime-switching model of \citet{VinkovskayaEkaterina2014Appm}, where the state process triggering the regime switches is not modelled. Besides, since here the events drive the dynamics of the state process, this subclass is different from the Markov-modulated Hawkes processes considered by \citet{Cohen:2013:regimeSwitchingHawkes} or \citet{Swishchuk:2017:compoundHawkes}, where the state process is a continuous-time Markov chain that jumps independently of the events. Moreover, the intensity in \citet{Cohen:2013:regimeSwitchingHawkes} depends only on the current state whereas, here and in \citet{Swishchuk:2017:compoundHawkes}, it may depend on all past states.
\end{generalised_hawkes_processes}

\subsection{Existence and uniqueness of hybrid marked point processes} \label{subsec:existence_uniqueness}

In this subsection, $\markSpace$ is not required to be a product space as in Subsection \ref{subsec:hybrid_marked_point_processes} but can be again an arbitrary complete separable metric space.

\subsubsection{The existence and uniqueness problem.}


A hybrid marked point processes (Definition \ref{def:hybrid_marked_pp}) is defined implicitly via its intensity process, which, in turn, depends on the history of the hybrid marked point process. Due to the self-referential nature of the definition, it is not clear a priori that such a marked point process exists. More generally, given an initial condition $N_{\leq 0}$ (see Subsection \ref{subsec:initial_condition}) and a measurable intensity functional $\psi:\markSpace\times\boundedlyFiniteMeasures{\theCSMS}\rightarrow\realLineNonNegative\cup\{\infty\}$, one can ask if there exists a unique non-explosive marked point process $N$ that satisfies the initial condition $N_{\leq 0}$ on $\realLineNonPositive$ and admits $\intensityFunctional$ as its intensity functional on $\realLinePositive$.

\citet{Massoulie:1998aa:StabilityResults} tackles this question by reformulating the existence problem as a Poisson-driven SDE, extending the works of \citet{Bremaud:1996aa:StabilityNonLinearHawkes}, \citet{grigelionis1971representation} and \citet{kerstan1964teilprozesse}. \citet{Delattre:2016aa} also employ this Poisson embedding technique in the context of Hawkes processes on infinite directed graphs. However, in these papers, strong existence and uniqueness is obtained by imposing a Lipschitz condition on the intensity functional $\intensityFunctional$. More precisely, it is assumed that there exists a non-negative kernel $\overline{k}:\realLinePositive\times\markSpace\times\markSpace\rightarrow\realLineNonNegative$ such that
\begin{equation} \label{eq:lipschitz_condition}
	| \intensityFunctionalAtMarkGivenPastWithIndex{m}{\xi}{} - \intensityFunctionalAtMarkGivenPastWithIndex{m}{\xi'}{} | \leq \iint_{\realLineNegative\times\markSpace}\overline{k}(-t',m',m)|\xi-\xi'|(dt',dm'),\quad m\in\markSpace, \xi,\xi'\in\boundedlyFiniteMeasures{\theCSMS}.
\end{equation}
Unfortunately, this condition is too restrictive in the context of hybrid marked point processes. A simple, yet natural, example of a hybrid marked point process not satisfying \eqref{eq:lipschitz_condition} is given in Subsection \ref{subsubsec:counterexample}. Hence, our goal is to construct a strong solution to a Poisson-driven SDE without imposing the Lipschitz condition \eqref{eq:lipschitz_condition} on the intensity functional $\intensityFunctional$.
We will in fact extend the existence result in \citet{Massoulie:1998aa:StabilityResults} by imposing only a weaker sublinearity condition on $\intensityFunctional$.
The idea to define a random measure as a strong solution to an SDE driven by another random measure was also studied by \citet{jacod1979calcul}. Similarly, a Lipschitz condition that does not seem to apply to Hawkes processes and hybrid marked point processes is required \citep[Chapter 14, Section 1]{jacod1979calcul}.

Let us also briefly review some weak existence and uniqueness results. \citet{jacod:1975:projection} proved that there exists a unique probability measure on the canonical space of marked point processes such that the canonical marked point process admits a given compensator.
However, this marked point process may be explosive a priori. Still, we will apply this result in the proof of Theorem \ref{thm:uniqueness_weak} below to obtain weak uniqueness. A similar approach is followed by \citet[Proposition 4.3.5, Corollary 4.4.4]{jacobsen2006point} who, furthermore, gives a domination condition on the intensity functional ensuring that the corresponding marked point process is non-explosive. Proposition \ref{prop:construction_is_dominated1} will be the counterpart of this result in the strong setting. These weak existence results are however limited to intensities with respect to the internal history $\internalHistoryCollection{N}$. The advantage of the strong setting is that the results of \citet{Massoulie:1998aa:StabilityResults} and the pathwise construction of this paper also hold when the intensity functional depends additionally on an auxiliary process, meaning that intensities with respect to larger filtrations can be considered. Besides, the Poisson-driven SDE representation of marked point processes considered in the strong setting directly suggests a simulation (thinning) algorithm. In fact, the Poisson embedding lemma (Lemma \ref{lem:poisson_embedding} below), which is a stepping stone to the strong setting, was first given in the simulation literature \citep{Lewis1976:SimulationPoissonProcesses, Ogata:1981aa}.

Finally, there is a third approach to obtain existence, based on a change of measure, see \citet[Theorem 11, p.~242]{bremaud:1981:MartingaleDynamics} and \citet{Hansen:2015:ExponentialMartingales}. While this technique also accommodates filtrations that are larger than the internal history, existence is generally obtained only on finite time intervals.

\subsubsection{The Poisson-driven SDE.}

Let $\probSpaceFuture$ be given and let $\drivingPoisson_{>0}:\realisationsSpaceFuture\rightarrow\boundedlyFiniteMeasures{\theCSMS\times\realLine}$ be a Poisson process on $\theCSMS\times\realLine$ with mean measure $dt\markSpaceMeasure(dm)dz$. As usual, denote by $\collection{\internalHistory{\drivingPoisson_{>0}}{t}}{t\in\realLine}$ the internal history of $\drivingPoisson_{>0}$ on $\realisationsSpaceFuture$.
In this Subsection, we work under the assumption that underlying probability space $\probSpace$ is the completion of the product probability space defined by \eqref{eq:probability_space}. In particular, $\realisationsSpace:=\realisationsSpaceInitialCondition\times\realisationsSpaceFuture$, where $\realisationsSpaceInitialCondition$ corresponds to the probability space of an initial condition $N_{\leq 0}$, see Subsection \ref{subsec:initial_condition}.
We extend $\drivingPoisson_{>0}$ to a mapping $M:\realisationsSpace\rightarrow\boundedlyFiniteMeasures{\theCSMS\times\realLine}$ by simply setting
\begin{equation} \label{eq:driving_Poisson}
	\drivingPoisson(\omega) := \drivingPoisson_{>0}(\omega_{>0}), \quad \omega=(\omega_{\leq 0},\omega_{>0})\in\realisationsSpace.
\end{equation}
Let $\aHistoryCollection=\collection{\aHistory{t}}{t\in \realLine}$ be the filtration on $\realisationsSpace$ such that, for all $t\in\realLine$, $\aHistory{t}$ is the $\prob$-completion of
	$\internalHistory{N_{\leq 0}}{t} \otimes \internalHistory{\drivingPoisson_{>0}}{t}$ in $\sigmaAlgebra$. In particular, the filtration $\aHistoryCollection$ is complete \citep[p.~123]{kallenberg2006foundations}.
Similarly to \citet{Massoulie:1998aa:StabilityResults}, we want to solve the following Poisson-driven SDE.
\newtheorem{the_SDE}[marked_point_process]{Definition}
\begin{the_SDE}[The Poisson-driven SDE] \label{SDE:the_Poisson_driven_SDE}
	Let $\intensityFunctional:\markSpace\times\boundedlyFiniteMeasures{\theCSMS}\rightarrow\realLineNonNegative\cup\{\infty\}$ be a given measurable functional. By a solution to the \emph{Poisson-driven SDE}, we mean an $\aHistoryCollection$-adapted non-explosive marked-point process $N:\realisationsSpace\rightarrow\mppSpace$ that solves
	\begin{equation} \label{eq:the_SDE_driven_by_Poisson}
	\begin{cases}
		N(dt,dm)=\drivingPoisson(dt,dm,(0,\intensityAtTimeAtMarkWithIndex{t}{m}{}]), \quad & t\in\realLinePositive,\, \as, \\
		\intensityAtTimeAtMarkWithIndex{\omega, t}{m}{} = \intensityFunctionalAtMarkGivenPastWithIndex{m}{\pointProcessAtTime{N(\omega)}{t}}{}, \quad & t\in\realLinePositive, m\in\markSpace, \omega\in\realisationsSpace, \\
		\restrictionToNonPositiveOfMeasure{N}(\omega)= N_{\leq 0}(\omega_{\leq 0}), & \omega=(\omega_{\leq 0},\omega_{>0}) \in\realisationsSpace, \as,
	\end{cases}
	\end{equation}
	where $N_{\leq 0}$ is a given initial condition (see Subsection \ref{subsec:initial_condition}).
\end{the_SDE}
\theoremstyle{plain}
Still, notice that our problem differs slightly as we only search for solutions in the space of non-explosive marked point processes, a smaller space than the one considered in \citet{Massoulie:1998aa:StabilityResults}.

\subsubsection{Assumptions.}

The following assumptions are only required for the strong existence result (Theorem \ref{thm:pathwise_existence} below). We first need to assume that the mark space $\markSpace$ has finite total mass.
\newtheorem{markspacebounded}{Assumption}
\renewcommand{\themarkspacebounded}{\Alph{markspacebounded}}
\begin{markspacebounded}[] \label{ass:markspacebounded}
The reference measure $\markSpaceMeasure$ is finite, i.e., $\markSpaceMeasure(\markSpace)<\infty$.
\end{markspacebounded}

Next, we need to control for both the intensity functional $\intensityFunctional$ and the initial condition $N_{\leq 0}$. We will prove Theorem \ref{thm:pathwise_existence} for two different scenarios. In the first scenario, the intensity is dominated by an increasing function of the total number of past events, while the number of events before time $0$ is finite.
\newtheorem{intensity_dominated_1}[markspacebounded]{Assumption}
\begin{intensity_dominated_1}[] \label{ass:intensity_dominated_1}
There exists a non-decreasing function $a:\integers\cup\{\infty\}\rightarrow\integers\cup\{\infty\}$ with $a(n)<\infty$ for all $n\in\integers$ and $a(\infty)=\infty$ such that:
\begin{enumerate}[label=\textup{(\roman*)}]
\item $\intensityFunctionalAtMarkGivenPastWithIndex{m}{\xi}{}\leq a(\xi((-\infty,0)\times\markSpace))$, $m\in\markSpace$, $\xi\in\boundedlyFiniteMeasures{\theCSMS}$ \textup{;}
\item $\sum_{n=0}^{\infty}a(n)^{-1}=\infty$.
\end{enumerate}
\end{intensity_dominated_1}

\newtheorem{initial_condition_1}[markspacebounded]{Assumption}
\begin{initial_condition_1}[] \label{ass:initial_condition_1}
The initial condition satisfies $N_{\leq 0}(\omega_{\leq 0},(-\infty,0]\times\markSpace))<\infty$ for all $\omega_{\leq 0}\in\realisationsSpaceInitialCondition$.
\end{initial_condition_1}

In the second scenario, the intensity functional $\intensityFunctional$ is dominated by a Hawkes functional. Note that this requirement is weaker than the Lipschitz condition \eqref{eq:lipschitz_condition} in \citet{Massoulie:1998aa:StabilityResults}.
\newtheorem{sublinear_intensity}[markspacebounded]{Assumption}
\begin{sublinear_intensity}[] \label{ass:sublinear_intensity}
There exists $\lambda_{0}\in\realLineNonNegative$ and a measurable function $\overline{k}:\realLinePositive\times\markSpace\times\markSpace \rightarrow\realLineNonNegative$ such that:
\begin{enumerate}[label=\textup{(\roman*)}]
	\item $\intensityFunctionalAtMarkGivenPastWithIndex{m}{\xi}{} \leq \lambda_{0} + \iint_{(-\infty,0)\times\markSpace}\overline{k}(-t',m',m)\xi(dt',dm')$, $m\in\markSpace$, $\xi\in\boundedlyFiniteMeasures{\theCSMS}$ \textup{;}
	\item $\rho:=\sup_{m\in\markSpace} \iint_{(0,\infty)\times\markSpace}\overline{k}(t',m',m)\markSpaceMeasure(dm')dt' < 1$ \textup{;}
	\item $\sup_{m\in\markSpace}k(t',m',m)<\infty$ for all $t'\in\realLinePositive$, $m'\in\markSpace$.
\end{enumerate}
\end{sublinear_intensity}
\newtheorem{initial_condition_good_expectation_and_kernel_contracts}[markspacebounded]{Assumption}
\begin{initial_condition_good_expectation_and_kernel_contracts}[] \label{ass:initial_condition_good_expectation_and_kernel_contracts}
The initial condition $N_{\leq 0}$ satisfies:
\begin{enumerate}[label=\textup{(\roman*)}]
	\item $\sup_{t>0,\, m\in\markSpace} \E{\iint_{(-\infty,0]\times\markSpace}\overline{k}(t-t',m',m)N_{\leq 0}(dt',dm')} <\infty$ \textup{;}
	\item $\tilde{\intensity}_{\leq 0}(\omega_{\leq 0},t):=\sup_{m\in\markSpace} \iint_{(-\infty,0]\times\markSpace}\overline{k}(t-t',m',m)N_{\leq 0}(\omega_{\leq 0},dt',dm') < \infty$, $\omega_{\leq 0}\in\realisationsSpaceInitialCondition$, $t\in\realLinePositive$.
\end{enumerate}
\end{initial_condition_good_expectation_and_kernel_contracts}
Note that Assumptions \ref{ass:sublinear_intensity}.(ii) and \ref{ass:initial_condition_good_expectation_and_kernel_contracts}.(i) are needed in order to reuse Theorem 2 in \citet{Massoulie:1998aa:StabilityResults}. It will allow us to dominate the marked point process $N$ by a Hawkes process with kernel $\overline{k}$.

\subsubsection{Existence.}

We construct a solution to the Poisson-driven SDE in two mains steps. First, by taking advantage of the discrete nature of the driving Poisson process, we construct in a pathwise fashion a mapping $N:\realisationsSpace\rightarrow\integerValuedMeasures{\theCSMS}$ that solves \eqref{eq:the_SDE_driven_by_Poisson} up to each event time, generalising the construction in \citet[Chapter 6, p.~302-306]{ccinlar2011probability} and \citet[p.~127]{Lindvall:1988aa}. Second, we dominate $N$ by a non-explosive marked point process to show that $N$ is itself non-explosive, generalising the argument in \citet[Lemma B.1, p.~30]{chevallier:2015:Hawkes_Generalised}. When working under Assumptions \ref{ass:sublinear_intensity} and \ref{ass:initial_condition_good_expectation_and_kernel_contracts}, these two steps must actually be performed concurrently.
Then, it turns out that this constructed $N$ admits $\intensityFunctional$ as its intensity functional on $\realLinePositive$ and, thus, solves the existence problem. The proof of the following theorem, which extends the existence result in \citet{Massoulie:1998aa:StabilityResults}, is given in Subsection \ref{subsec:pathwise_construction}.
\newtheorem{pathwise_existence}[marked_point_process]{Theorem}
\begin{pathwise_existence}[Strong existence] \label{thm:pathwise_existence}
Under either Assumptions \ref{ass:markspacebounded}, \ref{ass:intensity_dominated_1}, \ref{ass:initial_condition_1} or Assumptions  \ref{ass:markspacebounded}, \ref{ass:sublinear_intensity}, \ref{ass:initial_condition_good_expectation_and_kernel_contracts}, there exists a non-explosive marked point process $N:\realisationsSpace\rightarrow\mppSpace$ that solves the Poisson-driven SDE (Definition \ref{SDE:the_Poisson_driven_SDE}). Any such $N$ satisfies the strong initial condition $N_{\leq 0}$ and admits $\intensityFunctional$ as its intensity functional on $\realLinePositive$.
\end{pathwise_existence}
As a corollary, we obtain conditions that ensure the existence of hybrid marked point processes.
\newtheorem{existence_hybrid}[marked_point_process]{Corollary}
\begin{existence_hybrid}[Existence of hybrid marked point processes] \label{cor:existence_hybrid}
	Suppose that Assumption \ref{ass:markspacebounded} holds and $\|\phi\|_\infty<\infty$. Moreover, suppose that either Assumptions \ref{ass:intensity_dominated_1} and \ref{ass:initial_condition_1} or Assumptions \ref{ass:sublinear_intensity} and \ref{ass:initial_condition_good_expectation_and_kernel_contracts} hold with $\intensityFunctionalAtMarkGivenPastWithIndex{m}{\xi}{}$ replaced by $\eventFunctionalAtEventGivenPast{e}{\xi}$, where the dominating kernel $\overline{k}$ is now a function $\overline{k}:\realLinePositive\times\markSpace\times\eventSpace\rightarrow\realLineNonNegative$, and with the constraint $\rho<\|\phi\|_\infty^{-1}$.
		Then, there exists a hybrid marked point process $N:\realisationsSpace\rightarrow\mppSpace$ with transition function $\stateFunctional$ and event functional $\eventFunctional$ that satisfies the strong initial condition $N_{\leq 0}$.
\end{existence_hybrid}
\theoremstyle{definition}
\newtheorem{existence_state_hawkes}[marked_point_process]{Example}
\begin{existence_state_hawkes}[Existence of state-dependent Hawkes processes] \label{ex:existence_state_hawkes}
	When the transition function $\stateFunctional$ is bounded, the above corollary encompasses the case of state-dependent Hawkes processes (Example \ref{ex:generalised_hawkes_processes}) for either bounded kernels with no integrability constraint or unbounded kernels (up to the constraint \ref{ass:sublinear_intensity}.(iii)) with an integrability constraint.
\end{existence_state_hawkes}
 \theoremstyle{plain}

\subsubsection{Uniqueness.}

As \citet{Massoulie:1998aa:StabilityResults} considers point processes on $\theCSMS$ that are not necessarily non-explosive marked point processes, he uses the Lipschitz condition \eqref{eq:lipschitz_condition} to obtain strong uniqueness in a space of regular point processes. Here, since we restrict ourselves to non-explosive marked point processes, the enumeration representation allows us to prove strong uniqueness more easily without any specific assumptions. The proof is deferred until Subsection \ref{subsec:uniqueness}.
\newtheorem{uniqueness_pathwise}[marked_point_process]{Theorem}
\begin{uniqueness_pathwise}[Strong uniqueness] \label{thm:uniqueness_pathwise}
Let $N:\realisationsSpace\rightarrow\mppSpace$ and $N':\realisationsSpace\rightarrow\mppSpace$ be two non-explosive marked point processes solving the Poisson-Driven SDE (Definition \ref{SDE:the_Poisson_driven_SDE}). Then $N=N'$ $\as$
\end{uniqueness_pathwise}

By applying Theorem 3.4 in \citet{jacod:1975:projection}, we can also obtain weak uniqueness. Alternatively, we could also have applied Theorem 14.2.IV in \citet[p.~381]{daleyVereJonesVolume2}. The idea is that the intensity and the conditional distributions $\prob((T_{n+1},M_{n+1})\in \cdot\,|\, \internalHistory{N}{T_n})$ uniquely determine each other, see also \citet{last1995marked} and \citet[Theorem 4.3.2, p.~54]{jacobsen2006point}. Another approach, as suggested by \citet{Massoulie:1998aa:StabilityResults}, could be to use the fact that any marked point process with an intensity functional can be represented as the strong solution to a Poisson-driven SDE like in Definition \ref{SDE:the_Poisson_driven_SDE}, see \citet[Theorem 14.56,  p.~472]{jacod1979calcul},  and use the strong uniqueness result. We prove the following result in Subsection \ref{subsec:uniqueness}.
\newtheorem{uniqueness_weak}[marked_point_process]{Theorem}
\begin{uniqueness_weak}[Weak uniqueness] \label{thm:uniqueness_weak}
Let $N_1$ and $N_2$ be two non-explosive marked point processes (possibly on distinct probability spaces) that admit the same intensity functional $\intensityFunctional$ on $\realLinePositive$.	Assume also that both $N_1$ and $N_2$ satisfy the weak initial condition $N_{\leq 0}$. Then, we have that $\inducedProb{N_1}=\inducedProb{N_2}$, i.e., the induced probabilities measures on $\boundedlyFiniteMeasures{\theCSMS}$ coincide.
\end{uniqueness_weak}
As a corollary, we obtain the weak uniqueness of hybrid marked point processes.
\newtheorem{uniqueness_hybrid}[marked_point_process]{Corollary}
\begin{uniqueness_hybrid}[Uniqueness of hybrid marked point processes] \label{cor:uniqueness_hybrid}
	All hybrid marked point processes with transition function $\stateFunctional$ and event functional $\eventFunctional$ that satisfy the weak initial condition $N_{\leq 0}$ induce the same probability measure on $\boundedlyFiniteMeasures{\theCSMS}$.
\end{uniqueness_hybrid}
\theoremstyle{definition}
\newtheorem{uniqueness_weak_rem}[marked_point_process]{Remark}
\begin{uniqueness_weak_rem} \label{rem:uniqueness_weak_rem}
	Note that weak uniqueness might not hold for a general history $\aHistoryCollection$. Given an $\aHistoryCollection$-predictable process $\intensity$, there could be two marked point processes $N$ and $N'$ that both admit $\intensity$ as their $\aHistoryCollection$-intensity, but such that $\inducedProb{N}\neq \inducedProb{N'}$, see Proposition 9.54.(ii) in \citet{kallenberg2017random} for such an example. The fact the we restrict ourselves to the natural filtration $\internalHistoryCollection{N}$ is crucial here.
\end{uniqueness_weak_rem}
\theoremstyle{plain}

\section{Dynamics of hybrid marked point processes} \label{sec:hybrid_marked_pp}

In this section, we prove Theorem \ref{thm:implied_dynamics_characterisation}, which characterises the dynamics of hybrid marked point processes.

\subsection{Preliminaries}

We first present a lemma that helps us reuse some results in the literature that require a specific form for the filtration. It simply says that the information up to time $u$ is equal to the information up to time $t$ to which we add the information between time $t$ and $u$, where $t< u$.
\newtheorem{history_decomposition}[marked_point_process]{Lemma}
\begin{history_decomposition}[] \label{lem:history_decomposition}
Let $N$ be a non-explosive point process on $\realLine\times\anyCSMS$. Let $t,u\in\realLine$ such that  $u>t$. Then,  we have that $\internalHistory{N}{u} = \internalHistory{N}{t}	\vee \internalHistory{\pointProcessAfterTime{N}{t}}{u}$.
\end{history_decomposition}
\begin{proof}
	Note that
	\begin{equation*}
		\internalHistory{\pointProcessAfterTime{N}{t}}{u} = \sigma\left\{N(A\times U) \, : \, A\in\borelSetsOf{\realLine},\,  A\subset(t,u],\, U\in\anyCSMSBorel \right\} \quad \mbox{for all } u>t.
	\end{equation*}
	Then, clearly $\internalHistory{\pointProcessAfterTime{N}{t}}{u}\subset \internalHistory{N}{u}$. Also, $\internalHistory{N}{t}\subset \internalHistory{N}{u}$ and, thus $ \internalHistory{N}{t}	\vee \internalHistory{\pointProcessAfterTime{N}{t}}{u} \subset \internalHistory{N}{u}$. On the other hand, let $A\in\borelSetsOf{\realLine}$ be such that $A\subset (-\infty, u]$ and let $U\in\anyCSMSBorel$. We have that
	\begin{equation*}
		N(A\times U) = N(A\cap(-\infty,t])\times U) + N(A\cap(t,u]\times U).
	\end{equation*}
	The first term is $\internalHistory{N}{t}$-measurable while the second term is $\internalHistory{\pointProcessAfterTime{N}{t}}{u}$-measurable. Hence, $N(A\times U)$ is $\internalHistory{N}{t}	\vee \internalHistory{\pointProcessAfterTime{N}{t}}{u}$ measurable. Since, by definition, $\internalHistory{N}{u}$ is the smallest $\sigma$-algebra that makes all the $N(A\times U)$ measurable, this implies that $\internalHistory{N}{u} \subset \internalHistory{N}{t}	\vee \internalHistory{\pointProcessAfterTime{N}{t}}{u}$, which concludes the proof.
\end{proof}

As defined in Subsection \ref{subsec:intensity}, an intensity process has to always be finite. We verify that, if one finds a potentially infinite process that satisfies the definition of the intensity, then one can take a finite version of this process and identify it with the intensity.
\newtheorem{intensity_finite_version}[marked_point_process]{Lemma}
\begin{intensity_finite_version}[] \label{lem:intensity_finite_version}
Let $N$ be a non-explosive marked point process on $\theCSMS$ and let $\tilde{\intensity}:\realisationsSpace\times\realLinePositive\times\markSpace\rightarrow\realLineNonNegative\cup\{\infty\}$ be an $\aHistoryCollection$-predictable process that satisfies \eqref{eq:intensity_definition} for all non-negative $\aHistoryCollection$-predictable processes $H$. Then $N$ admits an $\aHistoryCollection$-intensity $\intensity:\realisationsSpace\times\realLinePositive\times\markSpace\mapsto\realLineNonNegative$ relative to $\markSpaceMeasure$ such that	 $\intensity(\omega,t,m)=\tilde{\intensity}(\omega,t,m)$ holds $\prob(d\omega)\markSpaceMeasure(dm)dt\aeverywhere$
\end{intensity_finite_version}
\begin{proof}
	Since the marked point process $N$ is non-explosive, using similar arguments as in Lemma L2 of \citet[p.~24]{bremaud:1981:MartingaleDynamics}, one can show that, for all bounded sets $A\in\borelSetsOf{\realLinePositive}$,
	\begin{equation*}
		\iint_{A\times\markSpace} \tilde{\intensity}(t,m)\markSpaceMeasure(dm)dt <\infty\,,\quad \mbox{a.s.,}
	\end{equation*}
	which implies that $\tilde{\intensity}(\omega,t,m)<\infty$ holds $\prob(d\omega)dt\markSpaceMeasure(dm)\aeverywhere$ By a composition argument (see the beginning of the proof of Lemma \ref{lem:poisson_embedding}), since $\tilde{\intensity}$ is $\aHistoryCollection$-predictable, we have that $(\omega,t,m)\mapsto\indicator_{\{\tilde{\intensity}(\omega,t,m)<\infty\}}$ is also $\aHistoryCollection$-predictable. It is then easy to check that $\intensity(\omega,t,m):=\indicator_{\{\tilde{\intensity}(\omega,t,m)<\infty\}}\tilde{\intensity}(\omega,t,m)$ is the $\aHistoryCollection$-intensity of $N$ where we use the convention $0\times\infty=0$.
\end{proof}

The next lemma says that by integrating the intensity against the state variable $x$, we obtain the intensity of the marked point process that tracks the event types, ignoring the state process.
\newtheorem{intensity_integrated_mpp}[marked_point_process]{Lemma}
\begin{intensity_integrated_mpp}[] \label{lem:intensity_integrated_mpp}
Let $N$ be a marked point process on $\theCSMS$ with $\aHistoryCollection$-intensity $\intensity$ relative to $\markSpaceMeasure$.	Then, $N_\eventSpace(\cdot):=N(\cdot\times\stateSpace)$ is a non-explosive marked point process on $\realLine\times\eventSpace$ with $\aHistoryCollection$-intensity $\intensity_\eventSpace:\realisationsSpace\times\realLinePositive\times\eventSpace\rightarrow\realLineNonNegative$ relative to  $\eventSpaceMeasure$ such that
			$\intensity_\eventSpace(\omega,t,e)=\int_\stateSpace \intensity(\omega,t,e,x)\stateSpaceMeasure(dx)$ holds $\prob(d\omega)dt\eventSpaceMeasure(de)\aeverywhere$
\end{intensity_integrated_mpp}
\begin{proof}
	Let $H:\realisationsSpace\times\realLinePositive\times\eventSpace\rightarrow\realLineNonNegative$ be an $\aHistoryCollection$-predictable non-negative process. Then, by applying the definition of $N_\eventSpace$ and using Tonelli's theorem, we obtain that
	\begin{align*}
		\E{\iint_{\realLinePositive\times\eventSpace}H(t,e)N_\eventSpace(dt,de)}&= \E{\iiint_{\realLinePositive\times\eventSpace\times\stateSpace}H(t,e)N(dt,de,dx)} \\
		&= \E{\iiint_{\realLinePositive\times\eventSpace\times\stateSpace}H(t,e)\intensity(t,e,x)\stateSpaceMeasure(dx)\eventSpaceMeasure(de)dt} \\
		&= \E{\iint_{\realLinePositive\times\eventSpace}H(t,e)\left(\int_\stateSpace \intensity(t,e,x)\stateSpaceMeasure(dx)\right)\eventSpaceMeasure(de)dt}.
	\end{align*}
	The process $\int_\stateSpace\intensity(t,e,x)\stateSpaceMeasure(dx)$, $t\in\realLinePositive$, $e\in\eventSpace$, is $\aHistoryCollection$-predictable, see for example Lemma 25.23 in \citet[p.~503]{kallenberg2006foundations} and we conclude using Lemma \ref{lem:intensity_finite_version}.
\end{proof}

We now check that an intensity functional applied to the history of a point process defines a predictable process.
\newtheorem{predictability_functional_process}[marked_point_process]{Lemma}
\begin{predictability_functional_process}[] \label{lem:predictability_functional_process}
Let $\intensityFunctional : \anyCSMS \times\boundedlyFiniteMeasures{\realLine\times\anyCSMS}\rightarrow\realLineNonNegative\cup\{\infty\}$ be a measurable functional and $N$ be a non-explosive point process on $\realLine\times\anyCSMS$ that is $\aHistoryCollection$-adapted. Then, the process $\intensity:\realisationsSpace\times\realLine\times\anyCSMS\rightarrow\realLineNonNegative\cup\{\infty\}$ defined by $\intensity(\omega,t,u)=\intensityFunctionalAtMarkGivenPastWithIndex{u}{\pointProcessAtTime{N(\omega)}{t}}{	}$, $\omega\in\realisationsSpace, t\in\realLine, u\in\anyCSMS$, is $\aHistoryCollection$-predictable.
\end{predictability_functional_process}
\begin{proof}
	By Lemma \ref{prop:history_left_continuous}, $\pointProcessAtTime{N(\omega)}{t}$ is left-continuous in $t$ and, by assumption, the process $\collection{\pointProcessAtTime{N}{t}}{t\in\realLine}$ is $\aHistoryCollection$-adapted. As a consequence, the mapping $(\omega,t)\mapsto\pointProcessAtTime{N(\omega)}{t}$ is $\aHistoryCollection$-predictable, see for example Lemmas 25.1  and 1.10 in \citet[p.~491, p.~6]{kallenberg2006foundations}. We then obtain that $\intensity$ is $\aHistoryCollection$-predictable by viewing it as the composition $(\omega,t,u)\mapsto(u, \pointProcessAtTime{N(\omega)}{t})\mapsto\intensityFunctionalAtMarkGivenPastWithIndex{u}{\pointProcessAtTime{N(\omega)}{t}}{}$
	and using the measurability of $\intensityFunctional$.
\end{proof}

The next lemma essentially says that if two predictable processes coincide at all event times of a marked point process, then they coincide everywhere under positive intensity. A less general variant of this result and its proof are suggested in \citet[Theorem T12, p.~31]{bremaud:1981:MartingaleDynamics}.
\newtheorem{equality_on_events}[marked_point_process]{Lemma}
\begin{equality_on_events}[] \label{lem:equality_on_events}
Let $N_\eventSpace$ be a non-explosive marked point process on $\realLine\times\eventSpace$ with $\aHistoryCollection$-intensity $\intensity_\eventSpace$ relative to $\eventSpaceMeasure$. Let $H_1:\realisationsSpace\times\realLinePositive\times\markSpace\rightarrow\realLineNonNegative\cup\{\infty\}$ and 	$H_2:\realisationsSpace\times\realLinePositive\times\markSpace\rightarrow\realLineNonNegative\cup\{\infty\}$ be two non-negative $\aHistoryCollection$-predictable processes. Then, $H_1=H_2$ holds $\prob(d\omega)N_\eventSpace(\omega,dt,de)\stateSpaceMeasure(dx)\aeverywhere$ if and only if $H_1 = H_2$ holds $\prob(d\omega)\intensity_\eventSpace(\omega,t,e)dt\markSpaceMeasure(de,dx)\aeverywhere$
\end{equality_on_events}
\begin{proof}
	By a composition argument, since $H_1$ and $H_2$ are $\aHistoryCollection$-predictable, we have that the function $(\omega,t,m)\mapsto\indicator_{\{H_1(\omega,t,m)\neq H_1(\omega,t,m)\}}$ is $\aHistoryCollection$-predictable (see the beginning of the proof of Lemma \ref{lem:poisson_embedding}). By Lemma 25.23 in \citet[p.~503]{kallenberg2006foundations}, we also have that the process $\int_\stateSpace \indicator_{\{H_1(\cdot,\cdot,\cdot, x)\neq H_1(\cdot,\cdot,\cdot, x)\}}\stateSpaceMeasure(dx)$ is $\aHistoryCollection$-predictable. Using the definition of the intensity and Tonelli's theorem, we obtain that
	\begin{align*}
		&\int_{\realisationsSpace}\int_{\realLinePositive\times\markSpace} \indicator_{\{H_1(\omega,t,m)\neq H_1(\omega,t,m)\}} \stateSpaceMeasure(dx)N_\eventSpace(\omega,dt,de)\prob(d\omega)\\
		&= \int_{\realisationsSpace}\int_{\realLinePositive\times\eventSpace}\left(\int_{\stateSpace} \indicator_{\{H_1(\omega,t,m)\neq H_1(\omega,t,m)\}} \stateSpaceMeasure(dx)\right)N_\eventSpace(\omega,dt,de)\prob(d\omega)\\
		&= \int_{\realisationsSpace}\int_{\realLinePositive\times\markSpace} \indicator_{\{H_1(\omega,t,m)\neq H_1(\omega,t,m)\}} \stateSpaceMeasure(dx)\intensity_\eventSpace(\omega,t,e)\eventSpaceMeasure(de)dt \prob(d\omega),
	\end{align*}
	from which the assertion follows.
\end{proof}

Finally, we show that the link between joint densities and conditional densities still holds when we pre-condition on a sub-$\sigma$-algebra. Since $\markSpace$ is a complete separable metric space and, in particular, Borel, random elements in $\markSpace$ always have regular conditional distributions \citep[p.~106, Theorem A1.2, p.~561]{kallenberg2006foundations}.
\newtheorem{conditional_conditional_densities}[marked_point_process]{Lemma}
\begin{conditional_conditional_densities}[] \label{lem:conditional_conditional_densities}
Let $(E,X)$ be a random element in $\markSpace$. Let $\mathcal{G}$ be a sub-$\sigma$-algebra, i.e., $\mathcal{G}\subset \sigmaAlgebra$, and let $A\in\mathcal{G}$ such that $\prob(A)>0$. Moreover, let $f:\realisationsSpace\times\markSpace \rightarrow\realLineNonNegative\cup\{\infty\}$ be a non-negative measurable function that is $\mathcal{G}\otimes\markSpaceBorel$-measurable. If we have
\begin{equation} \label{eq:joint_regular_density_function}
	\prob\left( E\in de, X\in dx \,|\, \mathcal{G} \right)\indicator_A =
		f(e,x)\markSpaceMeasure(de,dx)\indicator_A, \quad \as,
\end{equation}
then
\begin{equation*}
	\prob\left( X\in dx \,|\, \sigma(E)\vee \mathcal{G} \right)\indicator_A =
		\frac{f(E,x)}{\int_\stateSpace f(E,x')\stateSpaceMeasure(dx')}\stateSpaceMeasure(dx)\indicator_A, \quad \as
\end{equation*}
\end{conditional_conditional_densities}
\begin{proof}
 Let $B\in\borelSetsOf{\stateSpace}$, $G\in\mathcal{G}$ and $H\in\sigma(E)$. On the one hand,
	\begin{align} \label{eq:lem_conditional_eq1}
		\E{\indicator_{G}\indicator_{H}\indicator_{\{X\in B\}}\indicator_{A}}&=\E{\indicator_{G}h(E)\indicator_{\{X\in B\}}\indicator_{A}}=\E{\indicator_{G}\E{h(E)\indicator_{\{X\in B\}}\indicator_{A}\,|\,\mathcal{G}}} \nonumber \\
		&= \E{\indicator_{G}\indicator_A \int_{\eventSpace}\int_{B}h(e)f(e,x)\stateSpaceMeasure(dx)\eventSpaceMeasure(de)},
	\end{align}
	where we successively used Lemma 1.13 in \citet[p.~7]{kallenberg2006foundations} to write $\indicator_{H}=h(E)$ using a measurable function $h:\eventSpace\rightarrow\{0,1\}$, the Tower property, the disintegration theorem in \citet[Theorem 6.4, p.~108]{kallenberg2006foundations} with the regular conditional distribution of \eqref{eq:joint_regular_density_function} and, finally, the product form of $\markSpaceMeasure$. Note that, here, the disintegration theorem is applied to the probability measure $\prob(\cdot\cap A)/\prob(A)$ on the measurable space $(A, A\cap\sigmaAlgebra )$.
	On the other hand, observe that \eqref{eq:product_mark_space} and \eqref{eq:joint_regular_density_function} imply that
		\begin{equation*}
			\prob\left( E\in de \,|\, \mathcal{G} \right)\indicator_A =
		\int_\stateSpace f(e,x)\stateSpaceMeasure(dx)\eventSpaceMeasure(de) \indicator_A, \quad \as
		\end{equation*}
	 Then, using similar arguments,
	\begin{align*}
		&\E{\indicator_{G}\indicator_{H}\int_{B}\frac{f(E,x)}{\int_\stateSpace f(E,x')\stateSpaceMeasure(dx')}\stateSpaceMeasure(dx)\indicator_{A}} = \E{\indicator_{G}h(E)\int_{B}\frac{f(E,x)}{\int_\stateSpace f(E,x')\stateSpaceMeasure(dx')}\stateSpaceMeasure(dx)\indicator_{A}} \\
		&=\E{\indicator_{G}\E{h(E)\int_{B}\frac{f(E,x)}{\int_\stateSpace f(E,x')\stateSpaceMeasure(dx')}\stateSpaceMeasure(dx)\indicator_{A}\,\Big|\,\mathcal{G}}} \\
		&= \E{\indicator_{G}\indicator_A \int_{\eventSpace}\left(h(e)\int_{B}\frac{f(e,x)}{\int_\stateSpace f(e,x')\stateSpaceMeasure(dx')}\stateSpaceMeasure(dx)\right)\int_\stateSpace f(e,x')\stateSpaceMeasure(dx')\eventSpaceMeasure(de)}.
	\end{align*}
	Tonelli's theorem and \eqref{eq:lem_conditional_eq1} then imply that
	\begin{equation}\label{eq:proof:transition_prob:conditional_charact}
		\E{\indicator_{G}\indicator_{H}\indicator_{\{X\in B\}}\indicator_{A}} = \E{\indicator_{G}\indicator_{H}\int_{B}\frac{f(E,x)}{\int_\stateSpace f(E,x')\stateSpaceMeasure(dx')}\stateSpaceMeasure(dx)\indicator_{A}}.
	\end{equation}
	Using a monotone class argument, we show below that \eqref{eq:proof:transition_prob:conditional_charact} can be extended to
	\begin{equation} \label{eq:proof:transition_prob:conditional_charact_ext}
		\E{\indicator_{F}\indicator_{\{X\in B\}}\indicator_{A}} = \E{\indicator_{F}\int_{B}\frac{f(E,x)}{\int_\stateSpace f(E,x')\stateSpaceMeasure(dx')}\stateSpaceMeasure(dx)\indicator_{A}}
	\end{equation}
	for all $F\in\sigma(E)\vee\mathcal{G}$, which means exactly that
	\begin{equation*}
		\prob\left(X\in B \,|\,  \sigma(E) \vee\mathcal{G}\right)\indicator_A = \int_{B}\frac{f(E,x)}{\int_\stateSpace f(E,x')\stateSpaceMeasure(dx')}\stateSpaceMeasure(dx)\indicator_{A},\quad \as,
	\end{equation*}
	as asserted.

	To prove \eqref{eq:proof:transition_prob:conditional_charact_ext}, define the functions
	\begin{align*}
	\mu_{1}\,:\quad \sigma(E)\vee\mathcal{G} &\rightarrow [0,1]	\\
	F&\mapsto \mu_{1}(F):=\E{\indicator_{F}\indicator_{\{X\in B\}}\indicator_{A}} \,, \\
	\mu_{2}\,:\quad \sigma(E)\vee\mathcal{G} &\rightarrow [0,1]	\\
	F &\mapsto \mu_{2}(F):= \E{\indicator_{F}\int_{B}\frac{f(E,x)}{\int_\stateSpace f(E,x')\stateSpaceMeasure(dx')}\stateSpaceMeasure(dx)\indicator_{A}}.
	\end{align*}
	One can check that $\mu_{1}$ and $\mu_{2}$ are bounded measures on $(\realisationsSpace, \sigma(E)\vee\mathcal{G})$ (to swap an expectation with an infinite sum, use the monotone convergence theorem, see for example Theorem 1.19 in \citet[p.~11]{kallenberg2006foundations}). Define also the class $\mathcal{C}:= \{ G\cap H \,:\, G\in\mathcal{G},\, H\in\sigma(E)\}$.
	Equation \eqref{eq:proof:transition_prob:conditional_charact} means that $\mu_{1}(C)=\mu_{2}(C)$ for all $C\in\mathcal{C}$. Moreover, $\mathcal{C}$ is a $\pi$-system such that $\realisationsSpace\in\mathcal{C}$. Also, note that $\sigma(E)\cup\mathcal{G} \subset \mathcal{C} \subset \sigma(E)\vee\mathcal{G}$ and, thus, $\sigma(\mathcal{C}) = \sigma(E)\vee\mathcal{G}$. As a consequence, we can apply Lemma 1.17 in \citet[p.~9]{kallenberg2006foundations} to conclude that $\mu_{1}(F)=\mu_{2}(F)$ for all $F\in \sigma(E)\vee\mathcal{G}$, meaning that \eqref{eq:proof:transition_prob:conditional_charact_ext} holds.
\end{proof}

\subsection{Implied dynamics and characterisation} \label{subsec:implied_dynamics_characterisation}
\begin{proof}[Proof of Theorem \ref{thm:implied_dynamics_characterisation}]
	Recall that we denote by $\intensity$ the $\internalHistoryCollection{N}$-intensity of $N$ relative to $\markSpaceMeasure$ and by $\intensity_\eventSpace$ the $\internalHistoryCollection{N}$-intensity of $N_{\eventSpace}:=N(\cdot\times\stateSpace)$ relative to $\eventSpaceMeasure$.

	\textbf{Necessity.} Assume that $N$ is a hybrid marked point process with transition function $\stateFunctional$ and event functional $\eventFunctional$. We first observe that statement (i) holds simply by applying Lemma \ref{lem:intensity_integrated_mpp} and using the fact that $\stateFunctionalAtStateGivenEventGivenPast{\cdot}{e}{x}$ is a probability density for all $e\in\eventSpace$ and $x\in\stateSpace$.

	Next, we show that statement (ii) holds. This is clearly true when $\prob(\tau_t < \infty)=0$ and, thus, we assume that $\prob(\tau_t<\infty)>0$.
	By applying Theorem T6 in \citet[p.~236]{bremaud:1981:MartingaleDynamics}, we obtain that, for all $M\in\markSpaceBorel$,
\begin{equation*}
	\prob\left((E,X)\in M \,|\, \internalHistory{N}{\tau_t-} \right)\indicator_{\{\tau_t<\infty\}} = \frac{\int_{M}\intensity(\tau_t,m)\markSpaceMeasure(dm)}{\int_\markSpace\intensity(\tau_t,m')\markSpaceMeasure(dm')} \indicator_{\{\tau_t<\infty\}},\quad \as
\end{equation*}
	This is allowed since Lemma \ref{lem:history_decomposition} tells us that the filtration $\internalHistoryCollection{N}$ is within the framework of this result.
	Hence, we have identified the unique regular conditional distribution of $(E,X)$ given $\internalHistory{N}{\tau_t -}$ on the measurable space  $(\{\tau_t<\infty\}, \{\tau_t<\infty\}\cap\sigmaAlgebra)$ equipped with the measure $\prob(\cdot\cap \{\tau_t<\infty\})/\prob(\{\tau_t<\infty\})$ \citep[Theorem 6.3, p.~107]{kallenberg2006foundations}.
	Besides, observe that the mapping $(\omega,m)\mapsto\intensityAtTimeAtMarkWithIndex{\omega,\tau_t(\omega)}{m}{}\indicator_{\{\tau_t(\omega)<\infty\}}$ is $\internalHistory{N}{\tau_t-}\otimes\markSpaceBorel$-measurable, see for example Lemma 25.3 in \citet[p.~492]{kallenberg2006foundations}. Using Lemma 1.26 in \citet[p.~14]{kallenberg2006foundations}, we obtain that the function $f$ defined by
	\begin{equation*}
		f(\omega,e,x)=\frac{\intensity(\omega,\tau_t(\omega),e,x)}{\int_{\markSpace}\intensity(\omega,\tau_t(\omega),m')\markSpaceMeasure(dm')} \indicator_{\{\tau_t(\omega)<\infty\}},\quad \omega\in\realisationsSpace, e\in\eventSpace, x\in\stateSpace,
	\end{equation*}
	is $\internalHistory{N}{\tau_t-}\otimes\markSpaceBorel$-measurable.
	We can then apply Lemma \ref{lem:conditional_conditional_densities} with $\mathcal{G}=\internalHistory{N}{\tau_t-}$ and $A=\{\tau_t<\infty\}$. This yields that
	\begin{align} \label{eq:probability_next_state}
	\prob\left( X\in dx \,|\, \sigma(E)\vee \internalHistory{N}{\tau_t-} \right)\indicator_{\{\tau_t<\infty\}} &=
		\frac{\intensityAtTimeAtMarkWithIndex{\tau_t}{E,x}{}}{\int_\stateSpace \intensityAtTimeAtMarkWithIndex{\tau_t}{E,x'}{}\stateSpaceMeasure(dx')}\stateSpaceMeasure(dx)\indicator_{\{\tau_t<\infty\}},\quad \as
	\end{align}
	By viewing the term $\stateFunctionalAtStateGivenEventGivenPast{x}{e}{\stateProcess_t}$  as a measurable function $\varphi$ applied to $(e,x,\pointProcessAtTime{N}{t})$ where $\varphi(x,e\,|\,\xi)=\stateFunctionalAtStateGivenEventGivenPast{x}{e}{F(\xi)}$, and using the measurability of the state functional $F$ and the transition function $\stateFunctional$, we obtain by Lemma \ref{lem:predictability_functional_process} that $\stateFunctionalAtStateGivenEventGivenPast{x}{e}{\pointProcessAtTime{N}{t}}$, $t\in\realLine$, $e\in\eventSpace$, $m\in\markSpace$, is $\internalHistoryCollection{N}$-predictable.
	Similarly, note that $\eventFunctionalAtEventGivenPast{e}{\pointProcessAtTime{N}{t}}$, $t\in\realLine$, $e\in\eventSpace$, is also $\internalHistoryCollection{N}$-predictable (this will be useful when proving sufficiency).
	Besides, thanks to the assumption on $\intensity$,
	\begin{equation*}
		\frac{\intensityAtTimeAtMarkWithIndex{\omega, t}{e, x}{}}{\int_\stateSpace \intensityAtTimeAtMarkWithIndex{\omega, t}{e ,x'}{}\stateSpaceMeasure(dx')} = \stateFunctionalAtStateGivenEventGivenPast{x}{e}{\stateProcess_t(\omega)},\quad \prob(d\omega)dt\markSpaceMeasure(de,dx)\aeverywhere
	\end{equation*}
	Hence, using Lemma \ref{lem:equality_on_events}, \eqref{eq:probability_next_state} becomes
	\begin{align*}
	\prob\left( X\in dx \,|\, \sigma(E)\vee \internalHistory{N}{\tau_t-} \right)\indicator_{\{\tau_t<\infty\}} &= \stateFunctionalAtStateGivenEventGivenPast{x}{E}{\stateProcess_{\tau_t}}\stateSpaceMeasure(dx)\indicator_{\{\tau_t<\infty\}} , \quad \as
	\end{align*}
To obtain \eqref{eq:probability_of_next_state}, it remains to notice that $\stateProcess_{\tau_t}=\stateProcess_{t+}$ on $\{\tau_t<\infty\}$ since there is no event on the time interval $(t,\tau_t)$ by definition of $\tau_t$. Also, since the ground point process $N(\cdot\times\markSpace)$ admits an $\internalHistoryCollection{N}$-intensity, we have that $N(\{t\}\times\markSpace)=0$ $\as$, implying that $X_{t+}=X_t$ $\as$  To show the statement in Remark \ref{rem:transition_probability_remark}, simply use \eqref{eq:probability_of_next_state} and the tower property to obtain that
	\begin{align*}
		\E{\indicator_{F}\indicator_{\{\tau_t<\infty\}}\indicator_{X\in B}} = \E{\indicator_{F}\indicator_{\{\tau_t<\infty\}}\E{\indicator_{X\in B}\,|\,\sigma(E)\vee\internalHistory{N}{\tau_t-}}}=\E{\indicator_{F}\indicator_{\{\tau_t<\infty\}}\int_B\stateFunctionalAtStateGivenEventGivenPast{x}{E}{\stateProcess_{t}}\stateSpaceMeasure(dx)}
	\end{align*}
	for all $F\in\sigma(E)\vee\internalHistory{N}{t}$, $B\in\borelSetsOf{\stateSpace}$ and observe that $\int_B\stateFunctionalAtStateGivenEventGivenPast{x}{E}{\stateProcess_{t}}\stateSpaceMeasure(dx)$ is $\sigma(E)\vee\internalHistory{N}{t}$-measurable.

	\textbf{Sufficiency.} Assume that $N$ is a non-explosive marked point process on $\theCSMS$ such that it admits an $\internalHistoryCollection{N}$-intensity relative to $\markSpaceMeasure$ and such that statements (i) and (ii) hold. We want to show that $\intensityAtTimeAtMarkWithIndex{\omega,t}{e,x}{}=\stateFunctionalAtStateGivenEventGivenPast{x}{e}{\stateProcess_t(\omega)}\eventFunctionalAtEventGivenPast{e}{\pointProcessAtTime{N(\omega)}{t}}$ holds $\prob(d\omega)\markSpaceMeasure(de,dx)dt\aeverywhere$ For all $t\in\realLineNonNegative$, by using statement (ii), \eqref{eq:probability_next_state}, Lemmas \ref{lem:intensity_integrated_mpp} and \ref{lem:equality_on_events}, and statement (i), we obtain that
	\begin{align*}
		\stateFunctionalAtStateGivenEventGivenPast{x}{E}{\stateProcess_{\tau_t}}\stateSpaceMeasure(dx)\indicator_{\{\tau_t<\infty\}}&= \prob\left( X\in dx \,|\, \sigma(E)\vee \internalHistory{N}{\tau_t-} \right)\indicator_{\{\tau_t<\infty\}} \\
		&=\frac{\intensityAtTimeAtMarkWithIndex{\tau_t}{E,x}{}}{\int_\stateSpace \intensityAtTimeAtMarkWithIndex{\tau_t}{E,x'}{}\stateSpaceMeasure(dx')}\stateSpaceMeasure(dx)\indicator_{\{\tau_t<\infty\}} \\
		&=  \frac{\intensityAtTimeAtMarkWithIndex{\tau_t}{E,x}{}}{\intensity_\eventSpace(\tau_t,E)}\stateSpaceMeasure(dx)\indicator_{\{\tau_t<\infty\}}\\
		&= \frac{\intensityAtTimeAtMarkWithIndex{\tau_t}{E,x}{}}{\eventFunctionalAtEventGivenPast{E}{\pointProcessAtTime{N}{\tau_t}}}\stateSpaceMeasure(dx)\indicator_{\{\tau_t<\infty\}}, \quad \as
	\end{align*}
	This means that, for all $t\in\realLineNonNegative$, we have that
	\begin{equation*}
		\intensityAtTimeAtMarkWithIndex{\tau_t}{E,x}{}\indicator_{\{\tau_t<\infty\}} = \stateFunctionalAtStateGivenEventGivenPast{x}{E}{\stateProcess_{\tau_t}}\eventFunctionalAtEventGivenPast{E}{\pointProcessAtTime{N}{\tau_t}}\indicator_{\{\tau_t<\infty\}},\quad \stateSpaceMeasure(dx)\aeverywhere,\, \as
	\end{equation*}
	This holds $\as$ simultaneously for all $t\in\mathbb{Q}\cap\realLineNonNegative$, whence, using that the number of events in $N$ is countable and finite in any bounded time interval,
	\begin{equation*}
			\intensityAtTimeAtMarkWithIndex{\omega,t}{e,x}{} = \stateFunctionalAtStateGivenEventGivenPast{x}{e}{\stateProcess_t(\omega)}\eventFunctionalAtEventGivenPast{e}{\pointProcessAtTime{N(\omega)}{t}}, \quad \prob(d\omega)N_\eventSpace(\omega,dt,de)\stateSpaceMeasure(dx)\aeverywhere
	\end{equation*}
	By Lemma \ref{lem:equality_on_events}, the above equality then implies that
	\begin{equation*}
			\intensityAtTimeAtMarkWithIndex{\omega,t}{e,x}{} = \stateFunctionalAtStateGivenEventGivenPast{x}{e}{\stateProcess_t(\omega)}\eventFunctionalAtEventGivenPast{e}{\pointProcessAtTime{N(\omega)}{t}}, \quad \prob(d\omega)\intensity_\eventSpace(\omega,t,e)dt\markSpaceMeasure(de,dx)\aeverywhere
	\end{equation*}
	By noticing that, on $\intensity_\eventSpace(\omega,t,e)=0$, we have that $\eventFunctionalAtEventGivenPast{e}{\pointProcessAtTime{N(\omega)}{t}}=0$ holds $\prob(d\omega)dt\eventSpaceMeasure(de)\aeverywhere$ and that $\intensity(\omega,t,e,x)=0$ holds $\prob(d\omega)dt\eventSpaceMeasure(de)\stateSpaceMeasure(dx)\aeverywhere$ (using again Lemma \ref{lem:intensity_integrated_mpp}), we conclude that the above equation actually holds $\prob(d\omega)dt\markSpaceMeasure(de,dx)\aeverywhere$
\end{proof}

\section{Existence and uniqueness} \label{sec:pathwise_construction}

In this section, we prove the strong existence result (Theorem \ref{thm:pathwise_existence}) by means of a Poisson embedding lemma given below. Subsequently, we also prove the strong and weak uniqueness results  (Theorems \ref{thm:uniqueness_pathwise} and \ref{thm:uniqueness_weak}).

\subsection{Preliminaries}

\subsubsection{Example violating the Lipschitz condition.} \label{subsubsec:counterexample}

We give here an example of a hybrid marked point process that does not satisfy the Lipschitz condition \eqref{eq:lipschitz_condition}, implying that the existence and uniqueness results in \citet{Massoulie:1998aa:StabilityResults} do not apply.
\theoremstyle{definition}
\newtheorem{non_lipschitz_example}[marked_point_process]{Example}
\begin{non_lipschitz_example}[] \label{ex:non_lipschitz_example}
Set $\eventSpace=\left\{ 0, 1 \right\}$ and $\stateSpace=\left\{ 0, 1 \right\}$ with $\eventSpaceMeasure=\delta_{0} + \delta_{1}$ and $\stateSpaceMeasure=\delta_{0}+\delta_{1}$. Consider an intensity functional $\intensityFunctional$ that corresponds to a hybrid marked point process with transition function $\stateFunctional$ and event functional $\eventFunctional$ (see Definition \ref{def:hybrid_marked_pp}). Take $\eventFunctional$ to be a Hawkes functional of the form
\begin{equation*}
	\eventFunctionalAtEventGivenPast{e}{\xi}=\nu + \iint_{\realLineNegative\times\markSpace}k(-t',x',e)\xi(dt',de',dx'),
\end{equation*}
where $\nu\in\realLinePositive$ and $k:\realLinePositive\times\stateSpace\times\eventSpace \rightarrow\realLinePositive$ is continuous in time and strictly positive.
Let $t_0\in\realLineNegative$ and choose $\xi_{0},\xi_{1}\in\mppSpace$ such that $\xi_{0}$ and $\xi_{1}$ coincide on $(-\infty,t_0]$ (i.e., $\shift{t_0}\restrictionToNonPositiveOfMeasure{\xi_0}=\shift{t_0}\restrictionToNonPositiveOfMeasure{\xi_1}$) but $F(\xi_0)=0$ and $F(\xi_1)=1$ (thus, $\xi_0$ and $\xi_1$ do not coincide on $(t_0,0)$).  Assume also that $\phi(0 \,|\,0, 1) > \phi(0 \,|\,0, 0)$, $\eventFunctionalAtEventGivenPast{0}{\xi_1}<\infty$, and $ \iint_{(t_0,0)\times\markSpace}k(-t',x',0)\xi_1(dt',de',dx')>\iint_{(t_0,0)\times\markSpace}k(-t',x',0)\xi_0(dt',de',dx')$. Then, following some computations that are left to the reader,
\begin{equation*}
	|\intensityFunctionalAtMarkGivenPastWithIndex{0,0}{\xi_{1}}{} - \intensityFunctionalAtMarkGivenPastWithIndex{0,0}{\xi_{0}}{}| \geq \left( \phi(0 \,|\, 0, 1) - \phi(0 \,|\,0, 0) \right)\iint_{(-\infty,t_{0}]\times\markSpace} k(-t,x,0) \xi_{0}(dt,dx).
\end{equation*}
Next, consider any non-negative kernel $\overline{k}:\realLinePositive\times\markSpace\times\markSpace\rightarrow\realLineNonNegative$. We have that
\begin{equation*}
	\iint_{\realLineNegative\times\markSpace}\overline{k}(-t,m,0,0)|\xi_{1}-\xi_{0}|(dt,dm) = \iint_{(t_0,0)\times\markSpace}\overline{k}(-t,m,0,0)|\xi_{1}-\xi_{0}|(dt,dm).
\end{equation*}
We can now add as many points as necessary to $\xi_{0}$ and $\xi_{1}$ on $(-\infty,t_{0}]$ to guarantee that
\begin{equation*}
	|\intensityFunctionalAtMarkGivenPastWithIndex{0,0}{\xi_{1}}{} - \intensityFunctionalAtMarkGivenPastWithIndex{0,0}{\xi_{0}}{}| > \iint_{\realLineNegative\times\markSpace}\overline{k}(-t,m,0,0)|\xi_{1}-\xi_{0}|(dt,dm).
\end{equation*}
Consequently, the intensity functional $\intensityFunctional$ does not satisfy the Lipschitz condition \eqref{eq:lipschitz_condition}.
\end{non_lipschitz_example}
\theoremstyle{plain}

\subsubsection{Integration with respect to Poisson processes.}

We first clarify briefly the link between point processes and random measures. A random measure $M$ on a measurable space $(S,\mathcal{S})$ is a mapping $M: \realisationsSpace\times\mathcal{S}\rightarrow\realLineNonNegative\cup\{\infty\}$ such that $M(\omega,\cdot)$ is a measure on $(S,\mathcal{S})\,$ for all $\omega\in\realisationsSpace$ and $M(\cdot, A)$ is a random variable for all $A\in\mathcal{S}$, see \citet[p.~106]{kallenberg2006foundations} and \citet[Chapter 6, p.~243]{ccinlar2011probability}. Note that the concepts of internal history and adaptedness of Subsection \ref{subsec:shift_histories} can be directly extended to random measures.
Not surprisingly, point processes are exactly the boundedly finite integer-valued random measures.
\newtheorem{pprocesses_are_random_measures}[marked_point_process]{Proposition}
\begin{pprocesses_are_random_measures} \label{prop:pprocesses_are_random_measures}
Let $N$ be a random measure on $(\anyCSMS , \anyCSMSBorel)$	such that $N(\omega,\cdot)\in\boundedlyFiniteMeasures{\anyCSMS}$ for all $\omega\in\realisationsSpace$. Then $N$ is a non-explosive point process on $\anyCSMS$. In return, any non-explosive point process $N$ on $\anyCSMS$ is a random measure on $(\anyCSMS , \anyCSMSBorel)$ such that $N(\omega,\cdot)\in\boundedlyFiniteMeasures{\anyCSMS}$ for all $\omega\in\realisationsSpace$.
\end{pprocesses_are_random_measures}
\begin{proof}
	See Proposition 9.1.VIII in \citet[p.~8]{daleyVereJonesVolume2}.
\end{proof}

One can show that Poisson processes are Poisson random measures in the sense of \citet[Chapter 6, p.~249]{ccinlar2011probability}. This enables us to apply an important result on integration with respect to Poisson random measures. Before stating the result, we need to define what it means for a Poisson process to be Poisson relative to a filtration.
\theoremstyle{definition}
\newtheorem{poisson_relative_filtration}[marked_point_process]{Definition}
\begin{poisson_relative_filtration} \label{def:poisson_relative_filtration}
Let $N$ be a Poisson process on $\realLine\times\anyCSMS$ and $\aHistoryCollection=\collection{\aHistory{t}}{t\in\realLine}$ be a filtration. We say that $N$ is Poisson \textit{relative} to $\aHistoryCollection$ if	 for all $t\in\realLine$, the point process $\shift{t}\restrictionToNonPositiveOfMeasure{N}$ is $\aHistory{t}$-measurable and $\sigma(\pointProcessAfterTime{N}{t})$ is independent of $\aHistory{t}$.
\end{poisson_relative_filtration}
\theoremstyle{plain}
Trivially, a Poisson process $N$ is always Poisson relative to its internal history $\internalHistoryCollection{N}$. The next result plays a crucial role in the Poisson embedding technique, which is later used to construct marked point processes with given intensities.
\newtheorem{poisson_process_integration}[marked_point_process]{Theorem}
\begin{poisson_process_integration} \label{thm:poisson_process_integration}
Let $N$ be a Poisson process on $\realLine\times\anyCSMS$ with parameter measure $\nu$. Let $\aHistoryCollection=\collection{\aHistory{t}}{t\in\realLine}$ be a filtration and suppose that $N$ is Poisson relative to $\aHistoryCollection$. Then, for every non-negative $\aHistoryCollection$-predictable process $H:\realisationsSpace\times\realLine\times\anyCSMS\rightarrow\realLineNonNegative$, we have that
\begin{equation*}
	\E{\iint_{\realLine\times\anyCSMS}H(t,u)N(dt,du)} = \E{\iint_{\realLine\times\anyCSMS}H(t,u)\nu(dt,du)}.
\end{equation*}
\end{poisson_process_integration}
\begin{proof}
	See Theorem 6.2 in \citet[Chapter 6, p.~299]{ccinlar2011probability}.
\end{proof}

\subsubsection{Driving Poisson process.}

We prove that the mapping $M:\realisationsSpace\rightarrow\boundedlyFiniteMeasures{\theCSMS\times\realLine}$ defined by \eqref{eq:driving_Poisson} is still a Poisson process.

\newtheorem{extended_poisson_is_poisson}[marked_point_process]{Lemma}
\begin{extended_poisson_is_poisson} \label{lem:extended_poisson_is_poisson}
	 The mapping $\drivingPoisson :\realisationsSpace\rightarrow\boundedlyFiniteMeasures{\theCSMS\times\realLine}$ is a Poisson process on $\theCSMS\times\realLine$ with parameter measure $dt\markSpaceMeasure(dm)dz$. Moreover, $\drivingPoisson$ is Poisson relative to $\aHistoryCollection$.
\end{extended_poisson_is_poisson}
\begin{proof}
	By composition, using the measurability of $\drivingPoisson_{>0}$, it is easy to check that $M$ is a measurable mapping and, thus, it is a non-explosive point process. To show that $M$ is a Poisson process with parameter measure $dt\markSpaceMeasure(dm)dz$, it is enough notice that, for any $n\in\integers$, for every family of bounded sets $\collection{A_{i}}{i\in\{1,\ldots,n\}}$, for all $k_{1},\ldots,k_{n}\in\integers$,
	\begin{equation*}
		\prob(\drivingPoisson(A_{i})=k_{i},\, i=1,\ldots,n) = \probFuture(\drivingPoisson_{>0}(A_{i})=k_{i},\, i=1,\ldots,n),
	\end{equation*}
	and use the fact that $\drivingPoisson_{>0}$ is a Poisson process with parameter measure $dt\markSpaceMeasure(dm)dz$. To show that $\shift{t}\restrictionToNonPositiveOfMeasure{\drivingPoisson}$ is $\aHistory{t}$-measurable for any $t\in\realLine$, use the fact that $\shift{t}\restrictionToNonPositiveOfMeasure{\drivingPoisson_{>0}}$ is $\internalHistory{\drivingPoisson_{>0}}{t}$-measurable  (since a Poisson process is always Poisson relative to its internal history) along with a composition argument. Similarly, one can show that $\sigma(\shift{t}\restrictionToPositiveOfMeasure{\drivingPoisson} )\subset \{\varnothing,\realisationsSpaceInitialCondition\}\otimes\sigma(\shift{t}\restrictionToPositiveOfMeasure{\drivingPoisson_{>0}} )$ and, thus, to show that $\sigma(\shift{t}\restrictionToPositiveOfMeasure{\drivingPoisson} )$ is independent of $\aHistory{t}$, it is enough to show that $\{\varnothing,\realisationsSpaceInitialCondition\}\otimes\sigma(\shift{t}\restrictionToPositiveOfMeasure{\drivingPoisson_{>0}} )$ is independent of $\aHistory{t}$. For this, let $A_{\leq 0}\in \{\varnothing,\realisationsSpaceInitialCondition\}$, $A_{>0}\in \sigma(\shift{t}\restrictionToPositiveOfMeasure{\drivingPoisson_{>0}} )$, $B_{\leq 0}\in \internalHistory{N_{\leq 0}}{t}$ and $B_{>0}\in \internalHistory{\drivingPoisson_{>0}}{t}$. Then, using the fact that $\drivingPoisson _{>0}$ is Poisson relative to $\internalHistory{\drivingPoisson_{>0}}{}$, we have that
	\begin{align*}
		\prob(A_{\leq 0}\times A_{>0} \cap B_{\leq 0}\times B_{>0}) &= \prob(A_{\leq 0}\cap B_{\leq 0} \times A_{>0}\cap B_{>0}) = \probInitialCondition(A_{\leq 0}\cap B_{\leq 0})\probFuture(A_{>0}\cap B_{>0})\\
		&= \probInitialCondition(A_{\leq 0})\probInitialCondition(B_{\leq 0})\probFuture(A_{>0})\probFuture(B_{>0}) \\
		&= \prob(A_{\leq 0}\times A_{>0})\prob(B_{\leq 0}\times B_{>0}).
	\end{align*}
	This shows that two $\pi$-systems generating $\{\varnothing,\realisationsSpaceInitialCondition\}\otimes\sigma(\shift{t}\restrictionToPositiveOfMeasure{\drivingPoisson_{>0}} )$ and $\internalHistory{N_{\leq 0}}{t}\otimes\internalHistory{\drivingPoisson_{>0}}{t}$, respectively, are independent. We conclude using Lemma 3.6 in \citet[p.~50]{kallenberg2006foundations} that $\{\varnothing,\realisationsSpaceInitialCondition\}\otimes\sigma(\shift{t}\restrictionToPositiveOfMeasure{\drivingPoisson_{>0}} )$ and $\internalHistory{N_{\leq 0}}{t}\otimes\internalHistory{\drivingPoisson_{>0}}{t}$ are independent. We can then verify that $\{\varnothing,\realisationsSpaceInitialCondition\}\otimes\sigma(\shift{t}\restrictionToPositiveOfMeasure{\drivingPoisson_{>0}} )$ remains independent of the completion of $\internalHistory{N_{\leq 0}}{t}\otimes\internalHistory{\drivingPoisson_{>0}}{t}$, which by definition is $\aHistory{t}$. Indeed, remember that $\aHistory{t}:=\sigma(\mathcal{C})$ with $\mathcal{C}:=(\internalHistory{N_{\leq 0}}{t}\otimes\internalHistory{\drivingPoisson_{>0}}{t})\cup \mathcal{A}$ and where $\mathcal{A}$ denotes the class of all subsets of $\prob$-null sets in $\sigmaAlgebra$. It then suffices to notice that $\mathcal{C}$ is a $\pi$-system and that  $\{\varnothing,\realisationsSpaceInitialCondition\}\otimes\sigma(\shift{t}\restrictionToPositiveOfMeasure{\drivingPoisson_{>0}} )$ remains independent of $\mathcal{C}$.
\end{proof}

\subsubsection{Poisson-embedding lemma.}

We are now able to show the following key lemma which demonstrates how the extra-dimension of the Poisson process $M$ allows us to generate a marked point process with a given intensity.
\newtheorem{poisson_embedding}[marked_point_process]{Lemma}
\begin{poisson_embedding}[Poisson embedding] \label{lem:poisson_embedding}
Let $\lambda:\realisationsSpace\times\realLinePositive\times\markSpace\rightarrow\realLineNonNegative$ be an $\aHistoryCollection$-predictable process.
Then, the mapping
\begin{align}  \label{eq:embedding_equation}
	N :\quad \realisationsSpace\times\borelSetsOf{\realLinePositive\times\markSpace} &\rightarrow \realLineNonNegative\cup\{\infty\} \nonumber \\
	(\omega, A) &\mapsto N(\omega, A):=\iint_{A}\int_{(0,\lambda(\omega,t,m)]}\drivingPoisson(\omega,dt,dm,dz)
\end{align}
is an $\aHistoryCollection$-adapted integer-valued random measure on $\realLinePositive\times\markSpace$.
Moreover, for every non-negative $\aHistoryCollection$-predictable process $H:\realisationsSpace\times\realLinePositive\times\markSpace\rightarrow\realLineNonNegative$, we have that
\begin{equation*}
	\E{\iint_{\realLinePositive\times\markSpace}H(t,m)N(dt,dm)} = \E{\iint_{\realLinePositive\times\markSpace}H(t,m)\intensityAtTimeAtMarkWithIndex{t}{m}{}\markSpaceMeasure(dm)dt}.
\end{equation*}
\end{poisson_embedding}
\begin{proof}
First, let $A\in\borelSetsOf{\realLinePositive\times\markSpace}$ and consider the following composition
\begin{equation*}
	(\omega, t,m,z) \mapsto (\lambda(\omega,t,m), z) \mapsto \indicator_{(0,\lambda(\omega,t,m)]}(z)
\end{equation*}
to notice that $\indicator_{(0,\lambda(\omega,t,m)]}(z)$ is $\aHistoryCollection$-predictable by means of Lemma 1.7 and Lemma 1.8 in \citet[p.~5]{kallenberg2006foundations}. Then, the product $\indicator_{A}(t,m)\indicator_{(0,\lambda(\omega,t,m)]}(z)$ of two $\aHistoryCollection$-predictable processes is also $\aHistoryCollection$-predictable by Lemma 1.12 in \citet[p.~7]{kallenberg2006foundations}. This ensures that the integral
\begin{equation*}
	\iint_{A}\int_{(0,\lambda(\omega,t,m)]}\drivingPoisson(\omega,dt,dm,dz) = \iiint_{\realLinePositive\times\markSpace\times\realLine}\indicator_{A}(t,m)\indicator_{(0,\lambda(\omega,t,m)]}(z)\drivingPoisson(\omega,dt,dm,dz)
\end{equation*}
is well defined for all $\omega\in\realisationsSpace$ and that $N(\cdot, A)$ is a random variable (see Subsection \ref{subsec:integration}).

Second, let $\omega\in\realisationsSpace$. For any finite family of disjoint sets $A_{1},\ldots,A_{n}\in\borelSetsOf{\realLinePositive\times\markSpace}$, $n\in\integers$, we clearly have that $N(\omega,\bigcup_{i\leq n} A_{i}) = \sum_{i\leq n} N(\omega, A_{i})$, which means that $N(\omega, \cdot)$ is finitely additive. To prove that $N(\omega, \cdot)$ is countably additive, invoke finite additivity and apply the monotone convergence theorem. These first two steps show that $N$ is indeed a random measure.

Third, to show that $N$ is $\aHistoryCollection$-adapted, first consider processes $\lambda:\realisationsSpace\times\realLinePositive\times\markSpace\rightarrow\realLineNonNegative$ of the form $\lambda(\omega,t,m)=\indicator_F(\omega)\indicator_{(s,u]}(t)\indicator_C(m)$ where $F\in\aHistory{s}$, $s,u\in\realLinePositive$, $s<u$, $C\in\markSpaceBorel$. For any $t\in\realLinePositive$, any $A\in\borelSetsOf{\realLinePositive}$ such that $A\subset(0,t]$ and any $B\in\markSpaceBorel$, we obtain that
\begin{equation*}
	N(\omega, A\times B)=\indicator_F(\omega)M(\omega, A\cap(s,u]\times B\cap C \times (0,1]),
\end{equation*}
	which is $\aHistory{t}$-measurable since $M$ is $\aHistoryCollection$-adapted by Lemma \ref{lem:extended_poisson_is_poisson}. Hence, $N$ is $\aHistoryCollection$-adapted. To extend this result to any $\aHistoryCollection$-predictable process $\lambda$, one can use a monotone class argument like in the proof of Proposition \ref{prop:hawkes_functionals} for example.

Fourth, let $\omega\in\realisationsSpace$.
By the definition of $N$ and by the linearity of the integral, for all simple non-negative functions $f$ on $\realLinePositive\times\markSpace$, we have that
\begin{equation*}
	\iint_{\realLinePositive\times\markSpace}f(t,m)N(\omega,dt,dm) = \iiint_{\realLinePositive\times\markSpace\times\realLine}f(t,m)\indicator_{(0,\intensityAtTimeAtMarkWithIndex{\omega,t}{m}{}]}(z)M(\omega, dt,dm,dz).
\end{equation*}
Then, by Lemma 1.11 in \citet[p.~7]{kallenberg2006foundations} and the monotone convergence theorem, we have that the above equality holds for any $\borelSetsOf{\realLinePositive\times\markSpace}$-measurable non-negative function $f$. In particular, we have that, for all $\omega\in\realisationsSpace$,
\begin{equation*}
		\iint_{\realLinePositive\times\markSpace}H(\omega,t,m)N(\omega,dt,dm) = \iiint_{\realLinePositive\times\markSpace\times\realLine}H(\omega,t,m)\indicator_{(0,\intensityAtTimeAtMarkWithIndex{\omega,t}{m}{}]}(z)M(\omega,dt,dm,dz).
\end{equation*}

Fifth, using Lemma \ref{lem:extended_poisson_is_poisson} and Theorem \ref{thm:poisson_process_integration}, we deduce that
	\begin{align*}
		\E{\iint_{\realLinePositive\times\markSpace}H(t,m)N(dt,dm)} &= \E{\iiint_{\realLinePositive\times\markSpace\times\realLine}H(t,m)\indicator_{(0,\intensityAtTimeAtMarkWithIndex{t}{m}{}]}(z)M(dt,dm,dz)} \\
		&= \E{\iiint_{\realLinePositive\times\markSpace\times\realLine}H(t,m)\indicator_{(0,\intensityAtTimeAtMarkWithIndex{t}{m}{}]}(z)dt\markSpaceMeasure(dm)dz} \\
		&= \E{\iint_{\realLinePositive\times\markSpace}H(t,m)\intensityAtTimeAtMarkWithIndex{t}{m}{}\markSpaceMeasure(dm)dt}. \qedhere
	\end{align*}
\end{proof}
\theoremstyle{definition}
\newtheorem{poisson_embedding_remark}[marked_point_process]{Remark}
\begin{poisson_embedding_remark} \label{rem:poisson_embedding_remark}
	Similar results are given by \citet[Lemma 3, p.~1571]{Bremaud:1996aa:StabilityNonLinearHawkes}, \citet[Lemma 1, p.~3]{Massoulie:1998aa:StabilityResults} and \citet[Lemma 2.1, p.~4]{Torrisi2016:GaussianApprox}. They refer to \citet{Lewis1976:SimulationPoissonProcesses} and \citet{Ogata:1981aa} for proofs. The fifth part of our proof follows \citet[Proposition 14.7.I, p.~427]{daleyVereJonesVolume2}, but we could not find the first four parts anywhere. For the special case $\markSpace=\{0\}$ (i.e., for univariate point processes), a similar proof is given by \citet[Theorem 6.11, p.~303]{ccinlar2011probability} while an alternative proof is given by \citet[Theorem B.11]{chevalier:2015:microscopic}. Besides, our version of this lemma does not impose any local integrability condition on $\intensity$ and, thus, does not say if the obtained random measure $N$ is boundedly finite. Finally, note that \eqref{eq:embedding_equation} can be rewritten using the compact notation of \citet{Massoulie:1998aa:StabilityResults} as $N(dt,dm)=\drivingPoisson(dt,dm,[0,\lambda(t,m)])$, $t\in\realLinePositive$.
\end{poisson_embedding_remark}
\theoremstyle{plain}
We can now prove the final statement in Thoerem \ref{thm:pathwise_existence}, which we restate here as a corollary.
\newtheorem{intensity_functional_indeed}[marked_point_process]{Corollary}
\begin{intensity_functional_indeed} \label{cor:intensity_functional_indeed}
Let $N:\realisationsSpace\rightarrow\mppSpace$ be a solution to the Poisson-driven SDE (Definition \ref{SDE:the_Poisson_driven_SDE}) under either Assumptions \ref{ass:markspacebounded}, \ref{ass:intensity_dominated_1}, \ref{ass:initial_condition_1}, or Assumptions \ref{ass:markspacebounded}, \ref{ass:sublinear_intensity}, \ref{ass:initial_condition_good_expectation_and_kernel_contracts}. Then, $N$ admits $\intensityFunctional$ as its intensity functional on $\realLinePositive$.
\end{intensity_functional_indeed}
\begin{proof}
	Let $G\in\sigmaAlgebra$ be the almost sure event that \eqref{eq:the_SDE_driven_by_Poisson} holds. Consider the following modifications of $N$ and $\intensity$, where $\intensity$ is defined as in  \eqref{eq:the_SDE_driven_by_Poisson}:
	\begin{equation*}
		\tilde{N}(\omega):=N(\omega)\indicator_{G}(\omega),\quad \omega\in\realisationsSpace,\quad \mbox{and} \quad \tilde{\intensity}(\omega,t,m):=\intensityAtTimeAtMarkWithIndex{\omega,t}{m}{} \indicator_{G}(\omega),\quad \omega\in\realisationsSpace,t\in\realLinePositive,m\in\markSpace.
	\end{equation*}
	Then, $\tilde{N}$ and $\tilde{\lambda}$ satisfy \eqref{eq:embedding_equation} and, using either Assumptions \ref{ass:intensity_dominated_1}.(i) and \ref{ass:initial_condition_1} or Assumptions \ref{ass:sublinear_intensity}.(i) and \ref{ass:initial_condition_good_expectation_and_kernel_contracts}.(ii), one can check that $\tilde{\intensity}(\omega,t,m)<\infty$ for all $\omega\in\realisationsSpace$, $t\in\realLinePositive$, $m\in\markSpace$. Moreover, by Lemma \ref{lem:predictability_functional_process}, $\intensity$ is $\internalHistoryCollection{N}$-predictable, and, thus, $\aHistoryCollection$-predictable as $N$ is $\aHistoryCollection$-adapted. Since the filtration $\aHistoryCollection$ is complete, this implies that $\tilde{\intensity}$ is also $\aHistoryCollection$-predictable. Now, consider any non-negative $\internalHistoryCollection{N}$-predictable process $H:\realisationsSpace\times\realLinePositive\times\markSpace\rightarrow\realLineNonNegative$ and apply Lemma \ref{lem:poisson_embedding} to obtain
	\begin{align*}
		\E{\iint_{\realLinePositive\times\markSpace}H(t,m)N(dt,dm)} &= \E{\iint_{\realLinePositive\times\markSpace}H(t,m)\tilde{N}(dt,dm)}\\
		 	&=\E{\iint_{\realLinePositive\times\markSpace}H(t,m)\tilde{\intensity}(t,m)\markSpaceMeasure(dm)dt}\\
		 	&=\E{\iint_{\realLinePositive\times\markSpace}H(t,m)\intensityAtTimeAtMarkWithIndex{t}{m}{}\markSpaceMeasure(dm)dt}.
	\end{align*}
	We conclude that $N$ admits $\intensityFunctional$ as its intensity functional using Lemma \ref{lem:intensity_finite_version}.
\end{proof}

Given a non-explosive point process $N$ on $\theCSMS$ that solves \eqref{eq:the_SDE_driven_by_Poisson} or is defined through a Poisson embedding as in Lemma \ref{lem:poisson_embedding}, one can ask when $N$ is in fact a non-explosive marked point process. To this end, it is useful to define the following random measures induced by the driving Poisson process $M$:
	\begin{equation*}
		L_n(\omega,\cdot):=M(\omega,\cdot\times\markSpace\times(0,n]),\quad \omega\in\realisationsSpace, n\in\integers.
	\end{equation*}
	We are then able to find the following sufficient condition on $\intensity$.
\newtheorem{ground_measure_simple}[marked_point_process]{Lemma}
\begin{ground_measure_simple}[Simple ground measure] \label{lem:ground_measure_simple}
Let $\lambda:\realisationsSpace\times\realLinePositive\times\markSpace\rightarrow\realLineNonNegative$ be an $\aHistoryCollection$-predictable process and let $N$ be the $\aHistoryCollection$-adapted integer-valued random measure on $\realLinePositive\times\markSpace$ defined by \eqref{eq:embedding_equation}. Then, if Assumption \ref{ass:markspacebounded} holds and if $\sup_{m\in\markSpace}\intensityAtTimeAtMarkWithIndex{t}{m}{}<\infty$ for all $t\in\realLinePositive$, $\as$, we have that $N(\{t\}\times\markSpace)\leq 1$ for all $t\in\realLinePositive$, $\as$
\end{ground_measure_simple}
\begin{proof}
	Each $L_n$ is a Poisson random measure on $\realLine$ in the sense of \citet[Chapter 6, p.~249]{ccinlar2011probability} with boundedly finite parameter measure $n\markSpaceMeasure(\markSpace)dt$. Applying Theorem 2.17 in \citet[Chapter 6, p.~256]{ccinlar2011probability} for each $n\in\integers$, there exists a set $B\in\sigmaAlgebra$ such that $\prob(B)=1$ and such that, for all $\omega\in B$ and $n\in\integers$, $L_n(\omega)\in\mathcal{N}_{\realLine}^{\# g}$ (i.e., $L_n$, $n\in\integers$, are simultaneously simple). Next, let $A$ be the almost sure event that $\sup_{m\in\markSpace}\intensityAtTimeAtMarkWithIndex{t}{m}{}<\infty$ for all $t\in\realLinePositive$. Fix $\omega\in A\cap B$ and use the assumption on $\intensity$ to find that
	\begin{align*}
		N(\omega,\{t\}\times\markSpace) &= \int_{\{t\}}\int_{\markSpace}\int_{(0,\intensityAtTimeAtMarkWithIndex{\omega,s}{m}{}]} \drivingPoisson(\omega,ds,dm,dz) \\
		&\leq M\left(\omega,\{t\}\times\markSpace\times \left(0,\sup_{m\in\markSpace}\intensityAtTimeAtMarkWithIndex{\omega,t}{m}{}\right]\right) \\
		&\leq M\left(\omega,\{t\}\times\markSpace\times \left(0,p(\omega,t)\right]\right) = L_{p(\omega,t)}(\omega,\{t\}) \leq 1,
	\end{align*}
	where $p(\omega,t)\in\integers$ is such that $\sup_{m\in\markSpace}\intensityAtTimeAtMarkWithIndex{\omega,t}{m}{} \leq p(\omega,t)$.
\end{proof}

\subsection{Strong existence: pathwise construction via Poisson embedding} \label{subsec:pathwise_construction}

\subsubsection{Existence under Assumptions \ref{ass:markspacebounded}, \ref{ass:intensity_dominated_1}, \ref{ass:initial_condition_1}.}

We begin by proposing a construction of a candidate solution $N:\realisationsSpace\rightarrow\integerValuedMeasures{\theCSMS}$ to \eqref{eq:the_SDE_driven_by_Poisson}. We proceed in a pathwise fashion. Under Assumption \ref{ass:markspacebounded}, using the definition of a Poisson process, it is not difficult to see that, given $n\in\integers$, $L_n\in\boundedlyFiniteMeasures{\realLine}$ $\as$ This implies that
\begin{equation*}
	F_1:=\{\omega\in\realisationsSpace\,|\, L_n(\omega)\in\boundedlyFiniteMeasures{\realLine}, n\in\integers\}\in\sigmaAlgebra
\end{equation*}
is an almost sure event, which plays a key role in our pathwise construction.
\theoremstyle{definition}
\newtheorem{pathwise_construction}[marked_point_process]{Algorithm}
\begin{pathwise_construction}[] \label{constr:pathwise_construction}
Construct the mapping $N:\realisationsSpace\rightarrow\integerValuedMeasures{\theCSMS}$ as follows. For all $\omega=(\omega_{\leq 0}, \omega_{>0})\in F_1$, initialise $N_{0}(\omega):=N_{\leq 0}(\omega_{\leq 0})$, $T_{0}(\omega):=0$, $\mathbb{M}_0(\omega):=\varnothing$, and $\intensityAtTimeAtMarkWithIndex{\omega, t}{m}{0}:=\intensityFunctionalAtMarkGivenPastWithIndex{m}{\pointProcessAtTime{N_{0}(\omega)}{t}}{}$ for all $t\in\realLinePositive$, $m\in\markSpace$. Define recursively the sequences $\collection{N_{n}}{n\in\integers}$, $\collection{T_{n}}{n\in\integers}$, $\collection{\mathbb{M}_n}{n\in\integers}$, and $\collection{\intensity_{n}}{n\in\integers}$ as follows. For all $n\in\integers$,
\begin{itemize}
	\item if $T_{n}(\omega)<\infty$, then
		\begin{equation} \label{eq:event_time_recursion_pathwise}
		T_{n+1}(\omega):=\sup\left\{ u>T_{n}(\omega) \,:\, \iint_{(T_{n}(\omega),u)\times\markSpace}\int_{(0,\intensityAtTimeAtMarkWithIndex{\omega,t}{m}{n}]}\drivingPoisson(\omega,dt,dm,dz)=0 \right\} \,;
		\end{equation}
		\begin{itemize}
		\item if $T_{n+1}(\omega)<\infty$, then
		\begin{align}
		\mathbb{M}_{n+1}(\omega)&:=\left\{ m\in\markSpace \,:\, \drivingPoisson(\omega,\{T_{n+1}(\omega)\}\times\{m\}\times\ (0,\intensityAtTimeAtMarkWithIndex{\omega,T_{n+1}(\omega)}{m}{n}])>0\right\} \,;  \nonumber \\
		N_{n+1}(\omega)&:= \sum_{i=1}^{n+1}\sum_{m\in \mathbb{M}_{i}(\omega)} \drivingPoisson(\omega,\{T_{i}(\omega)\}\times\{m\}\times\ (0,\intensityAtTimeAtMarkWithIndex{\omega,T_{i}(\omega)}{m}{i-1}])\delta_{(T_{i}(\omega),m)} \,; \label{eq:point_process_pathwise_recursion} \\
		\intensityAtTimeAtMarkWithIndex{\omega,t}{m}{n+1}&:=\intensityFunctionalAtMarkGivenPastWithIndex{m}{\pointProcessAtTime{N_{n+1}(\omega)}{t}}{}, \quad t\in\realLinePositive,m\in\markSpace \,; \label{eq:intensity_pathwise_recursion}
	\end{align}
		\item if $T_{n+1}(\omega)=\infty$, then
		\begin{align*}
		\mathbb{M}_{n+1}(\omega)&:=\varnothing \,; \\
		N_{n+1}(\omega)&:= N_{n}(\omega) \,; \\
		\intensityAtTimeAtMarkWithIndex{\omega,t}{m}{n+1}&:=\intensityAtTimeAtMarkWithIndex{\omega,t}{m}{n}, \quad t\in\realLinePositive,m\in\markSpace\,;
	\end{align*}
		\end{itemize}
	\item if $T_{n}(\omega)=\infty$, then
	\begin{align*}
		T_{n+1}(\omega)&:=\infty \,; \\
		\mathbb{M}_{n+1}(\omega)&:=\varnothing \,; \\
		N_{n+1}(\omega)&:= N_{n}(\omega) \,; \\
		\intensityAtTimeAtMarkWithIndex{\omega,t}{m}{n+1}&:=\intensityAtTimeAtMarkWithIndex{\omega,t}{m}{n}, \quad t\in\realLinePositive,m\in\markSpace.
	\end{align*}
\end{itemize}
For all $\omega=(\omega_{\leq 0}, \omega_{>0})\in\realisationsSpace\setminus F_1$, set $N_n(\omega):=N_{\leq 0}(\omega_{\leq 0})$, $T_n(\omega):=\infty$, $\mathbb{M}_{n}(\omega):=\varnothing$ , $\intensityAtTimeAtMarkWithIndex{\omega,t}{m}{n} :=0$, $t\in\realLinePositive$, $m\in\markSpace$, for all $n\in\integers$.
Then, for all $\omega\in\realisationsSpace$, for all $n\in\integers$, define $N(\omega)$ on $(-\infty,T_{n+1}(\omega))$ by $\pointProcessAtTime{N(\omega)}{T_{n+1}(\omega)} :=\pointProcessAtTime{N_n(\omega)}{T_{n+1}(\omega)}$. Define also the explosion time $T_{\infty}(\omega):=\lim_{n\rightarrow\infty} T_{n}(\omega)$. If $T_{\infty}(\omega) <\infty$, extend $N(\omega)$ to $[T_{\infty}(\omega),\infty)$ by $\restrictionToNonNegativeOfMeasure{\shift{T_{\infty}(\omega)}N(\omega)}:=0$. This is equivalent to defining $N(\omega)$ as
\begin{equation*}
	N(\omega):=\lim_{n\rightarrow\infty}N_n(\omega)= \sum_{n=1}^{\infty}\sum_{m\in \mathbb{M}_{n}(\omega)} \drivingPoisson(\omega,\{T_{n}(\omega)\}\times\{m\}\times\ (0,\intensityAtTimeAtMarkWithIndex{\omega,T_{n}(\omega)}{m}{n-1}])\delta_{(T_{n}(\omega),m)}\indicator_{\{T_n(\omega)<\infty\}}.
\end{equation*}
\end{pathwise_construction}
\theoremstyle{plain}
Algorithm \ref{constr:pathwise_construction} would be ill-defined if the set in \eqref{eq:event_time_recursion_pathwise} were empty. This would mean that there are infinitely many events just after the time $T_{n}$. The following proposition shows that this actually never happens and, thus, ensures that Algorithm \ref{constr:pathwise_construction} is well-defined. We also need to prove that the set $\mathbb{M}_n$ is finite and that $N_{n}(\omega)\in\boundedlyFiniteMeasures{\theCSMS}$ for all $n\in\integers$, because otherwise $\intensityAtTimeAtMarkWithIndex{\omega,t}{m}{n}$ might be ill-defined ($\intensityFunctional$ is a functional on $\markSpace\times\boundedlyFiniteMeasures{\theCSMS}$).
\newtheorem{pathwise_construction_well_defined}[marked_point_process]{Proposition}
\begin{pathwise_construction_well_defined}[] \label{prop:pathwise_construction_well_defined}
In Algorithm \ref{constr:pathwise_construction}, under Assumptions  \ref{ass:markspacebounded}, \ref{ass:intensity_dominated_1}.(i) and \ref{ass:initial_condition_1}, we have that, for every $\omega\in F_1$, $\mbox{\textup{card}}(\mathbb{M}_n(\omega))<\infty$, $N_{n}(\omega)\in\boundedlyFiniteMeasures{\theCSMS}$, $\|\intensity_{i}\|(\omega):=\sup_{t>0,\, m\in\markSpace}\intensityAtTimeAtMarkWithIndex{\omega,t}{m}{i}<\infty$, for all $n\in\integers$, and
\begin{align*}
	\left\{ u>T_{n}(\omega) \,:\, \iint_{(T_{n}(\omega),u)\times\markSpace}\int_{(0,\intensityAtTimeAtMarkWithIndex{\omega,t}{m}{n}]}\drivingPoisson(\omega,dt,dm,dz)=0 \right\} \neq \varnothing\quad &\mbox{for all } n\in\integers \mbox{ s.t. } T_{n}(\omega)<\infty.
\end{align*}
Hence, Algorithm \ref{constr:pathwise_construction} is well-defined under these assumptions.
\end{pathwise_construction_well_defined}
\begin{proof}
	We show the desired result by induction.
	Take any $\omega=(\omega_{\leq 0},\omega_{>0})\in F_1$. Let $n\in\integers$ and for all $i\in\integers$ such that $i< n$ and $T_{i}(\omega)<\infty$, assume that
	\begin{equation} \label{eq:non_empty_set_induction_proof}
		\left\{ u>T_{i}(\omega) \,:\, \iint_{(T_{i}(\omega),u)\times\markSpace}\int_{(0,\intensityAtTimeAtMarkWithIndex{\omega,t}{m}{i}]}\drivingPoisson(\omega,dt,dm,dz)=0 \right\} \neq \varnothing.
	\end{equation}
	For all $i\in\integers$ such that $i\leq n$, assume that $N_{i}(\omega)\in\boundedlyFiniteMeasures{\theCSMS}$ and that $\|\intensity_{i}\|(\omega):=\sup_{t>0,\, m\in\markSpace}\intensityAtTimeAtMarkWithIndex{\omega,t}{m}{i}<\infty$. If $T_{n}(\omega)=\infty$, then, by construction, this is also true for $n+1$.

	Now, assume that $T_{n}(\omega)<\infty$. We first show that \eqref{eq:non_empty_set_induction_proof} holds also for $i=n$.
	Take any $\varepsilon>0$. We have that
	\begin{align*}
		 \iint_{(T_{n}(\omega),T_{n}(\omega)+\varepsilon)\times\markSpace}\int_{(0,\intensityAtTimeAtMarkWithIndex{\omega,t}{m}{n}]}\drivingPoisson(\omega,dt,dm,dz) &\leq \drivingPoisson(\omega,(T_{n}(\omega),T_{n}(\omega)+\varepsilon)\times\markSpace\times(0,\|\intensity_{n}\|(\omega)])\\
		 &\leq L_{p_n(\omega)}(\omega,(T_{n}(\omega),T_{n}(\omega)+\varepsilon)) =:U_{n}(\omega,\varepsilon) < \infty,
	\end{align*}
	where $p_n(\omega)\in\integers$ is such that $\|\intensity_{n}\|(\omega)\leq p_n(\omega)$ and we used the fact that $L_{p_n(\omega)}(\omega)\in\boundedlyFiniteMeasures{\realLine}$. If $U_{n}(\omega,\varepsilon)=0$, then clearly \eqref{eq:non_empty_set_induction_proof} is satisfied for $i=n$. If not, $\drivingPoisson(\omega)$ has a finite number of points in $(T_{n}(\omega),T_{n}(\omega)+\varepsilon)\times\markSpace\times(0,\|\intensity_{n}\|(\omega)]$ and there exists $0<\varepsilon'<\varepsilon$ such that $U_{n}(\omega,\varepsilon')=0$, in which case \eqref{eq:non_empty_set_induction_proof} is again satisfied for $i=n$. Note that the integral in \eqref{eq:non_empty_set_induction_proof} is well-defined since $\intensityAtTimeAtMarkWithIndex{\omega, \cdot}{\cdot}{i}$ is a measurable function on $\realLinePositive\times\markSpace$ for all $\omega\in\realisationsSpace$. To see this, consider the composition $(t,m)\mapsto (m,\pointProcessAtTime{N_i(\omega)}{t})\mapsto\intensityFunctionalAtMarkGivenPastWithIndex{m}{\pointProcessAtTime{N_i(\omega)}{t}}{}$ and use Lemma \ref{prop:history_left_continuous}, the measurability of $\intensityFunctional$ and Lemma 1.8 in \citet{kallenberg2006foundations}.

	Second, we show that  $\mbox{\textup{card}}(\mathbb{M}_{n+1})<\infty$ and $N_{n+1}(\omega)\in\boundedlyFiniteMeasures{\theCSMS}$. If $T_{n+1}(\omega)=\infty$, then this is immediate. If not, using again that $L_{p_n(\omega)}(\omega)\in\boundedlyFiniteMeasures{\realLine}$,
	\begin{align*}
		\sum_{m\in \mathbb{M}_{n+1}}\drivingPoisson(\omega,\{T_{n+1}(\omega)\}\times\{m\}\times (0,\intensityAtTimeAtMarkWithIndex{\omega,T_{n+1}(\omega)}{m}{n}])&\leq\drivingPoisson(\omega,\{T_{n+1}(\omega)\}\times\markSpace\times (0,\|\intensity_{n}\|(\omega)] ) \\
		&\leq L_{p_n(\omega)}(\omega,\{T_{n+1}(\omega)\})<\infty,
	\end{align*}
	which implies that the set $\mathbb{M}_{n+1}(\omega)$ is finite and, in view of \eqref{eq:point_process_pathwise_recursion}, that $N_{n+1}(\omega)\in\boundedlyFiniteMeasures{\theCSMS}$. Note that this also proves that $N_{n+1}(\omega,(0,T_{n+1}(\omega))\times\markSpace)<\infty$.

	Third, we show that $\|\intensity_{n+1}\|(\omega)<\infty$. If $T_{n+1}(\omega)=\infty$, then this is immediate. If not, by \eqref{eq:intensity_pathwise_recursion} and using Assumptions \ref{ass:intensity_dominated_1}.(i) and \eqref{eq:point_process_pathwise_recursion}, we
	 have that for all $t>0,\, m=(x,e)\in\markSpace$,
	\begin{align}
		\intensityAtTimeAtMarkWithIndex{\omega,t}{m}{n+1} &\leq a\left(N_{n+1}(\omega,(-\infty,t)\times\markSpace)\right) \nonumber \\
		&= a\left( N_{n+1}(\omega,(-\infty,0]\times\markSpace) + N_{n+1}(\omega,(0,t)\times\markSpace) \right)	\nonumber \\
		&\leq a\left(N_{\leq 0}(\omega_{\leq 0},(-\infty,0]\times\markSpace) + N_{n+1}(\omega,(0,T_{n+1}(\omega))\times\markSpace \label{eq:intermediary_inequality_initialisation})\right).
	\end{align}
	Since $N_{n+1}(\omega,(0,T_{n+1}(\omega))\times\markSpace)<\infty$ and, by Assumption \ref{ass:initial_condition_1}, $N_{\leq 0}(\omega_{\leq 0},(-\infty,0]\times\markSpace)<\infty$, this implies that $\|\intensity_{n+1}\|(\omega)<\infty$.

	Regarding the basis of this induction, it is immediate that $N_0(\omega) = N_{\leq 0}(\omega_{\leq 0})\in\boundedlyFiniteMeasures{\theCSMS}$. To see that $\|\lambda_{0}\|(\omega)<\infty$, simply set $n=-1$ in \eqref{eq:intermediary_inequality_initialisation}.
\end{proof}

We show that the constructed mapping $N:\realisationsSpace\rightarrow\integerValuedMeasures{\theCSMS}$ satisfies indeed \eqref{eq:the_SDE_driven_by_Poisson} up to each event time.
\newtheorem{pathwise_construction_solves_SDE}[marked_point_process]{Proposition}
\begin{pathwise_construction_solves_SDE}[] \label{prop:pathwise_construction_solves_SDE}
Under Assumptions \ref{ass:markspacebounded}, \ref{ass:intensity_dominated_1}.(i) and  \ref{ass:initial_condition_1}, the mapping $N:\realisationsSpace\rightarrow\integerValuedMeasures{\theCSMS}$ given by Algorithm \ref{constr:pathwise_construction} is such that $N(\omega)$ solves \eqref{eq:the_SDE_driven_by_Poisson} on $(-\infty,T_{n}(\omega))$ for all $n\in\integers$, for all $\omega\in F_1$.
\end{pathwise_construction_solves_SDE}
\begin{proof}
	Define the process $\intensityAtTimeAtMarkWithIndex{\omega,t}{m}{}:=\intensityFunctionalAtMarkGivenPastWithIndex{m}{\pointProcessAtTime{N(\omega)}{t}}{}$, $\omega\in\realisationsSpace, t\in(0, T_{\infty}(\omega)),\, m\in\markSpace$.
	Take any $\omega=(\omega_{\leq 0},\omega_{>0})\in F_1$. By construction, $\restrictionToNonPositiveOfMeasure{N(\omega)} = N_{0}(\omega)=N_{\leq 0}(\omega_{\leq 0})$ and, thus, $N$ satisfies the strong initial condition $N_{\leq 0}$. Take any $n\in\integers$ such that $T_{n+1}(\omega)<\infty$ in  Algorithm \ref{constr:pathwise_construction} and consider the time interval $(T_{n}(\omega),T_{n+1}(\omega)]$. By construction, we have that
	\begin{equation} \label{eq:pathwise_SDE_recursive}
		N_{n+1}(\omega,dt,dm)=\drivingPoisson(\omega,dt,dm,(0,\intensityAtTimeAtMarkWithIndex{\omega,t}{m}{n} ]) \quad \mbox{for all } t\in(T_{n}(\omega),T_{n+1}(\omega)].
	\end{equation}
	But, by definition, on $(-\infty,T_{n+1}(\omega)]$, $N(\omega)=N_{n+1}(\omega)$ and thus, for all $t\in (0,T_{n+1}(\omega)]$, $m\in\markSpace$,
	\begin{align*}
		\intensityAtTimeAtMarkWithIndex{\omega,t}{m}{}
		&=\intensityFunctionalAtMarkGivenPastWithIndex{m}{\pointProcessAtTime{N(\omega)}{t}}{}
		=\intensityFunctionalAtMarkGivenPastWithIndex{m}{\pointProcessAtTime{N_{n+1}(\omega)}{t}}{}
		=\intensityFunctionalAtMarkGivenPastWithIndex{m}{\pointProcessAtTime{N_{n}(\omega)}{t}}{}
		= \intensityAtTimeAtMarkWithIndex{\omega,t}{m}{n},
	\end{align*}
	by the definition \eqref{eq:intensity_pathwise_recursion} of $\intensity_{n}$, since $N_{n+1}(\omega)$ and $N_{n}(\omega)$ can only differ by a mass at time $T_{n+1}(\omega)$.
	Consequently, \eqref{eq:pathwise_SDE_recursive} can be rewritten on $(T_{n}(\omega),T_{n+1}(\omega)]$ as
	\begin{equation} \label{eq:SDE_on_positive_real_line}
		N(\omega,dt,dm)=\drivingPoisson(\omega,dt,dm,(0,\intensityAtTimeAtMarkWithIndex{\omega,t}{m}{} ]).
	\end{equation}
	This shows that the constructed $N(\omega)$ solves \eqref{eq:the_SDE_driven_by_Poisson} on $(-\infty,T_{n}(\omega)]$ for all $n\in\integers$ such that $T_{n}(\omega)<\infty$. Now, if there is $n\in\integers$ such that $T_{n}(\omega)<\infty$ and $T_{n+1}(\omega)=\infty$, then clearly the constructed $N(\omega)$ is null on $(T_{n}(\omega),\infty)$ and by similar arguments, \eqref{eq:SDE_on_positive_real_line} holds on $(T_{n}(\omega),\infty)$. This now allows us to conclude that $N(\omega)$ solves \eqref{eq:the_SDE_driven_by_Poisson} on $(-\infty,T_{n}(\omega))$ for all $n\in\integers$ in both cases $T_{n}(\omega)<\infty$ and $T_{n}(\omega)=\infty$.
\end{proof}

It will also be crucial for the strong existence proof to show that, for all $n\in\integers$, $N_n$ is adapted to the filtration $\aHistoryCollection$ and $\intensity_n$ is $\aHistoryCollection$-predictable.
\newtheorem{construction_adapted_predictable}[marked_point_process]{Proposition}
\begin{construction_adapted_predictable}[] \label{prop:construction_adapted_predictable}
In Algorithm \ref{constr:pathwise_construction}, for all $n\in\integers$, $\intensity_n$ is $\aHistoryCollection$-predictable, $N_n$ is  an $\aHistoryCollection$-adapted non-explosive point process and $T_n$ is an $\aHistoryCollection$-stopping time.
\end{construction_adapted_predictable}
\begin{proof}
	We proceed by induction. Regarding the basis, as the filtration $\aHistoryCollection$ is complete, clearly $N_0$ is $\aHistoryCollection$-adapted and $T_0$ is an $\aHistoryCollection$-stopping time. Now assume that $N_n$ is $\aHistoryCollection$-adapted and $T_n$ is an $\aHistoryCollection$-stopping time for some $n\in\integers$. First, observe that this implies that $\intensity_n$ is $\aHistoryCollection$-predictable by simply using the identity $\intensityAtTimeAtMarkWithIndex{\omega, t}{m}{n}=\intensityFunctionalAtMarkGivenPastWithIndex{m}{\pointProcessAtTime{N_n(\omega)}{t}}{}\indicator_{F_1}(\omega)$ and invoking Lemma \ref{lem:predictability_functional_process}, the fact that $\internalHistory{N_n}{t}\subset \aHistory{t}$, $t\in\realLine$, and the assumption that $\aHistoryCollection$ is complete.
	Second, let $t\in\realLine$ and notice that
	\begin{equation} \label{eq:proof:stopping_time}
		\{T_{n+1}\leq t\} =\left\{\iiint_{\theCSMS\times\realLine}\indicator_{(T_n,t]}(s)\indicator_{(0,\intensity_n(s,m)]}(z)M(ds,dm,dz) > 0\right\}\cap\{T_n<\infty\}\cap F_1.
	\end{equation}
	Because $T_n$ is an $\aHistoryCollection$-stopping time, we have that $\collection{\indicator_{(T_n,t]}(s)}{s\in\realLine}$ is $\aHistoryCollection$-adapted and left-continuous, implying that it is $\aHistoryCollection$-predictable, see for example Lemma 25.1 in \citet[p.~491]{kallenberg2006foundations}. Adapting the arguments of the third part of the proof of Lemma \ref{lem:poisson_embedding}, we deduce that the first event on the right-hand side of \eqref{eq:proof:stopping_time} belongs to $\aHistory{t}$ and so $T_{n+1}$ is an $\aHistoryCollection$-stopping time.
	Third, using Proposition \ref{prop:pathwise_construction_solves_SDE} and looking at Algorithm \ref{constr:pathwise_construction}, notice that $N_{n+1}$ satisfies
	\begin{equation*}
	\begin{cases}
		N_{n+1}(\omega,dt,dm)=\drivingPoisson(\omega,dt,dm,(0,\intensityAtTimeAtMarkWithIndex{\omega,t}{m}{n}\indicator_{\{t\leq T_{n+1}(\omega)\}}]), \quad & \omega\in\realisationsSpace, t\in\realLinePositive, \\
		\restrictionToNonPositiveOfMeasure{N}(\omega)= N_{\leq 0}(\omega_{\leq 0}), & \omega=(\omega_{\leq 0},\omega_{>0}) \in\realisationsSpace,
	\end{cases}
	\end{equation*}
	where $\intensityAtTimeAtMarkWithIndex{t}{m}{n}\indicator_{\{t\leq T_{n+1}\}}$ is $\aHistoryCollection$-predictable as a product of $\aHistoryCollection$-predictable processes, note that $\indicator_{\{t\leq T_{n+1}\}}$ is $\aHistoryCollection$-adapted and left-continuous since $T_{n+1}$ is an $\aHistoryCollection$-stopping time. Now, applying Lemma \ref{lem:poisson_embedding}, it follows that $N_{n+1}$ is indeed $\aHistoryCollection$-adapted.
\end{proof}

We are now in a position to prove Theorem \ref{thm:pathwise_existence} under Assumption \ref{ass:markspacebounded}, \ref{ass:intensity_dominated_1} and \ref{ass:initial_condition_1} for the following intensity functional:
\begin{align*}
	\intensityFunctional' : \markSpace\times\boundedlyFiniteMeasures{\theCSMS} &\rightarrow \realLineNonNegative\cup\{\infty\}  \\
	(m,\xi) &\mapsto \intensityFunctional'(m\,|\,\xi):= a\left(\xi((-\infty,0)\times\markSpace)\right).
\end{align*}
Still, note that the first step of the following proof remains true for general intensity functionals $\intensityFunctional$ that satisfy Assumption \ref{ass:intensity_dominated_1} and will be reused in other parts of the proof of Theorem \ref{thm:pathwise_existence}.
\begin{proof}[Proof of Theorem \ref{thm:pathwise_existence}, Part 1]
	Let $N:\realisationsSpace\rightarrow\integerValuedMeasures{\theCSMS}$ be given by Algorithm \ref{constr:pathwise_construction} under Assumptions \ref{ass:markspacebounded}, \ref{ass:intensity_dominated_1} and \ref{ass:initial_condition_1}, which is well-defined by Proposition \ref{prop:pathwise_construction_well_defined}, and consider here the special case $\intensityFunctional=\intensityFunctional'$. We will prove that $N$ admits a version that solves the Poisson-driven SDE. We proceed in four steps.

	First, notice that for all $\omega\in\realisationsSpace$, $t<T_\infty(\omega)$, there exists $n\in\integers$ such that $\pointProcessAtTime{N(\omega)}{t}=\pointProcessAtTime{N_n(\omega)}{t}$, which implies by Proposition \ref{prop:construction_is_dominated1} that the process
	\begin{equation*}
		\intensityAtTimeAtMarkWithIndex{\omega,t}{m}{}:=\intensityFunctionalAtMarkGivenPastWithIndex{m}{\pointProcessAtTime{N(\omega)}{t}}{}\indicator_{F_1}(\omega)\indicator_{\{t<T_\infty(\omega)\}},\quad \omega\in\realisationsSpace, t\in\realLinePositive, m\in\markSpace,
	\end{equation*}
	is well-defined and finite, and that, for all $\omega\in\realisationsSpace$, $t\in\realLinePositive$, $m\in\markSpace$,
	\begin{equation*}
		\intensityAtTimeAtMarkWithIndex{\omega,t}{m}{}=\lim_{n\rightarrow\infty}\intensityFunctionalAtMarkGivenPastWithIndex{m}{\pointProcessAtTime{N_n(\omega)}{t}}{}\indicator_{F_1}(\omega)\indicator_{\{t<T_\infty(\omega)\}}= \lim_{n\rightarrow\infty}\intensityAtTimeAtMarkWithIndex{\omega,t}{m}{n}\indicator_{\{t<T_\infty(\omega)\}}.
	\end{equation*}
	By Proposition \ref{prop:pathwise_construction_solves_SDE}, and because of the way we constructed $N$, we have that $N$ and $\intensity$ satisfy \eqref{eq:embedding_equation}. By Proposition \ref{prop:construction_adapted_predictable}, for all $n\in\integers$, $\intensity_n$ is $\aHistoryCollection$-predictable and $T_n$ is an $\aHistoryCollection$-stopping time. Since $T_{\infty}=\lim_{n\rightarrow\infty}T_n$, we have that $T_{\infty}$ is an $\aHistoryCollection$-predictable time, which implies by Lemma 25.3.(ii) in \citet[p.~492]{kallenberg2006foundations} that $\indicator_{\{t<T_\infty\}}$ is $\aHistoryCollection$-predictable. As $\intensity$ is a limit of $\aHistoryCollection$-predictable processes, we have that $\intensity$ is also $\aHistoryCollection$-predictable by Lemma 1.9 in \citet[p.~6]{kallenberg2006foundations}. Consequently, we can apply Lemma \ref{lem:poisson_embedding} to obtain that $N$ is an $\aHistoryCollection$-adapted integer-valued random measure. The main goal of the next steps is to show that $T_\infty=\infty$ $\as$

	Second, following Proposition \ref{prop:pathwise_construction_well_defined}, we can see that $\sup_{m\in\markSpace}\intensityAtTimeAtMarkWithIndex{\omega,t}{m}{}<\infty$ for all $t\in\realLinePositive$, $\omega\in\realisationsSpace$. Hence, by Lemma \ref{lem:ground_measure_simple}, there exists and almost sure event $G\in\sigmaAlgebra$ on which $N(\{t\}\times\markSpace\})\leq 1$ for all $t\in\realLine$. Let $\tilde{N}$, $\collection{\tilde{N}_n}{n\in\integers}$, $\collection{\tilde{T}_n}{n\in\integers}$ and $\tilde{T}_{\infty}$ coincide with $N$, $\collection{N_n}{n\in\integers}$, $\collection{T_n}{n\in\integers}$ and $T_{\infty}$ on $G$. Outside $G$, set $\tilde{N}:=0$, $\tilde{N}_n:=0$, $\tilde{T}_n:=\infty$, for all $n\in\integers$, and $\tilde{T}_\infty:=\infty$. Define the random measures on $\realLine$
	\begin{equation*}
		\tilde{N}_{\markSpace}(\cdot):=\tilde{N}(\cdot\times\markSpace),\quad \tilde{N}_{\markSpace,n}(\cdot):=\tilde{N}_n(\cdot\times\markSpace), \quad n\in\integers,
	\end{equation*}
and define the process
	\begin{equation*}
		\tilde{\intensity}(\omega,t):=\lim_{n\rightarrow\infty}a(\tilde{N}_{\markSpace,n}(\omega,(-\infty,t)))\indicator_{\{t<\tilde{T}_{\infty}(\omega)\}}=a(\tilde{N}_{\markSpace}(\omega,(-\infty,t)))\indicator_{\{t<\tilde{T}_{\infty}(\omega)\}},\quad \omega\in\realisationsSpace, t\in\realLinePositive.
	\end{equation*}
	Since $\{\tilde{T}_n\leq t\}=\{\tilde{N}_\markSpace((0,t])\geq n\}$, $\tilde{T}_n$ is in fact an $\internalHistoryCollection{\tilde{N}_\markSpace}$-stopping time and, thus, reusing the argument in the first step, we have that $\indicator_{\{t<\tilde{T}_{\infty}\}}$ is $\internalHistoryCollection{\tilde{N}_\markSpace}$-predictable. Moreover, by Lemma \ref{lem:predictability_functional_process}, we have that $\collection{a(\tilde{N}_{\markSpace,n}((-\infty,t)))}{t>0}$ is $\internalHistoryCollection{\tilde{N}_{\markSpace,n}}$-predictable and, thus, $\internalHistoryCollection{\tilde{N}_{\markSpace}}$-predictable. Hence, using again Lemma 1.9 in \citet[p.~6]{kallenberg2006foundations}, we have that $\tilde{\intensity}$ is also $\internalHistoryCollection{\tilde{N}_\markSpace}$-predictable. Next, because $\tilde{N}=N$ $\as$ and $\tilde{\intensity}(t)=\intensity(t,m)$ for all $t\in\realLinePositive$, $m\in\markSpace$, $\as$, and because Lemma \ref{lem:poisson_embedding} applies to $N$ and $\intensity$, we have that, for any non-negative $\internalHistoryCollection{\tilde{N}_\markSpace}$-predictable process $H:\realisationsSpace\times\realLinePositive\rightarrow\realLineNonNegative$,
	\begin{align*}
	\E{\int_{\realLinePositive}H(t)\tilde{N}_\markSpace(dt)} &= \E{\iint_{\realLinePositive\times\markSpace}H(t)N(dt,dm)}\\
		 	&=\E{\iint_{\realLinePositive\times\markSpace}H(t)\intensity(t,m)\markSpaceMeasure(dm)dt}\\
		 	&=\E{\int_{\realLinePositive}H(t)\tilde{\intensity}(t)\markSpaceMeasure(\markSpace)dt}.
	\end{align*}
	Consequently, $\tilde{N}_\markSpace$, or equivalently $\collection{\tilde{T}_n}{n\in\integers}$, defines a simple point process on $\realLinePositive$ with $\internalHistoryCollection{\tilde{N}_\markSpace}$-predictable projection $\collection{\markSpaceMeasure(\markSpace) \int_0^t\tilde{\intensity}(s)ds}{t>0}$ in the sense of \citet{jacod:1975:projection}.

	Third, by Lemma \ref{lem:history_decomposition}, $\internalHistory{\tilde{N}_\markSpace}{t}=\internalHistory{\tilde{N}_\markSpace}{0}\vee\internalHistory{\restrictionToPositiveOfMeasure{\tilde{N}_\markSpace}}{t}$, and, thus, Assumption A.1 of \citet{jacod:1975:projection} holds, see also the proof of Theorem \ref{thm:uniqueness_weak} and Remark \ref{rem:jacod_notations}. Then, by Proposition 3.1 in \citet{jacod:1975:projection}, we have that, conditional on $N_{\leq 0}((-\infty,0])=n_0\in\integers$, $S_n:=\tilde{T}_{n+1}-\tilde{T}_n$, $n\in\integers$, follows an exponential distribution with parameter $a(n+n_0)\markSpaceMeasure(\markSpace)$ and the $\collection{S_n}{n\in\integers}$ are independent. Thanks to Assumption \ref{ass:intensity_dominated_1}.(ii), by Example 3.1.4 in \citet[p.20]{jacobsen2006point}, we deduce that, conditional on $N_{\leq 0}((-\infty,0])=n_0$, $\tilde{T}_\infty=\lim_{n\rightarrow\infty}\tilde{T}_n = \infty$ $\as$, see also Proposition 12.19 in \citet[p.240]{kallenberg2006foundations}. Consequently, it holds that $\tilde{T}_\infty = \infty$ $\as$ unconditionally.

	Fourth, following these first three steps, we have proved that there exists a version of $N$ such that $N\in\mppSpace$, i.e., this version of $N$ is a non-explosive marked point process, see Proposition \ref{prop:pprocesses_are_random_measures}, and such that $N$ solves the Poisson-driven SDE. We conclude by Corollary \ref{cor:intensity_functional_indeed}.
\end{proof}

To prove Theorem \ref{thm:pathwise_existence} under Assumptions \ref{ass:markspacebounded}, \ref{ass:intensity_dominated_1} and \ref{ass:initial_condition_1} in the general case, we will use a solution to the special case $\intensityFunctional = \intensityFunctional'$ to show that the constructed mapping $N:\realisationsSpace\rightarrow\integerValuedMeasures{\theCSMS}$ actually takes values in $\mppSpace$. First, we need to define what we mean for a marked point process $N$ to be dominated by another marked point process $\dominatingMarkedPointProcess$.
\theoremstyle{definition}
\newtheorem{marked_point_process_dominated}[marked_point_process]{Definition}
\begin{marked_point_process_dominated}[] \label{def:marked_point_process_dominated}
Let $\xi,\overline{\xi}\in\integerValuedMeasures{\theCSMS}$. We say that $\xi$ is \textit{dominated} by $\overline{\xi}$ and write $\xi\prec\overline{\xi}$ if, for all $A\in\csmsBorel$, $\xi(A) \leq \overline{\xi} (A)$.
Let $T\in\realLine$. We say that $\xi$ is dominated by $\overline{\xi}$ on $(-\infty,T]$ if $\restrictionToNonPositiveOfMeasure{\shift{T}\xi}\prec \restrictionToNonPositiveOfMeasure{\shift{T}\overline{\xi}}$.
Consider two mappings $N:\realisationsSpace\rightarrow\integerValuedMeasures{\theCSMS}$ and $\dominatingMarkedPointProcess: \realisationsSpace\rightarrow\integerValuedMeasures{\theCSMS}$.	We say that $N$ is \textit{dominated} by $\dominatingMarkedPointProcess$ if $N\prec \dominatingMarkedPointProcess$ $\as$
\end{marked_point_process_dominated}
\theoremstyle{plain}
When $N\prec \dominatingMarkedPointProcess$ $\as$, one could also say that $N$ is a thinning of $\dominatingMarkedPointProcess$.
Indeed, notice that $\xi\prec\overline{\xi}$ implies that all the atoms of $\xi$ are also atoms of $\overline{\xi}$.

We will now show that the constructed mapping $N:\realisationsSpace\rightarrow\integerValuedMeasures{\theCSMS}$ is dominated by any solution to the special case $\intensityFunctional=\intensityFunctional'$.
\newtheorem{construction_is_dominated1}[marked_point_process]{Proposition}
\begin{construction_is_dominated1}[] \label{prop:construction_is_dominated1}
Let $N':\realisationsSpace\rightarrow\mppSpace$ be a solution to the Poisson-driven SDE  with intensity functional $\intensityFunctional'$. Then, under Assumptions \ref{ass:markspacebounded}, \ref{ass:intensity_dominated_1}.(i) and  \ref{ass:initial_condition_1}, the mapping $N:\realisationsSpace\rightarrow\integerValuedMeasures{\theCSMS}$ obtained from Algorithm \ref{constr:pathwise_construction} satisfies $N\prec N'$ $\as$
\end{construction_is_dominated1}
\begin{proof}
	Fix $\omega=(\omega_{\leq 0},\omega_{>0})\in A\cap F_1$, where $A\in\sigmaAlgebra$ is the almost sure event that $N'$ solves \eqref{eq:the_SDE_driven_by_Poisson}, where $\intensityFunctional$ is replaced by $\intensityFunctional'$.  Clealry, we have that $N(\omega)\prec \dominatingMarkedPointProcess(\omega)$ on $(-\infty,0]$. Now take any $n\in\integers$ such that $T_{n}(\omega)<\infty$ and assume that $N(\omega)\prec N'(\omega)$ on $(-\infty,T_{n}(\omega)]$. If $T_{n+1}(\omega)=\infty$, then $N(\omega)$ is null on $(T_{n},\infty)$ and we have $N(\omega)\prec N'(\omega)$. If $T_{n+1}(\omega)<\infty$, we have that for all $t\in(T_{n}(\omega),T_{n+1}(\omega)]$, $m\in\markSpace$,
	\begin{align*}
		\intensityAtTimeAtMarkWithIndex{\omega,t}{m}{}&=\intensityFunctionalAtMarkGivenPastWithIndex{m}{\pointProcessAtTime{N(\omega)}{t}}{} \\
		\mbox{(by construction)} \quad &= \intensityFunctionalAtMarkGivenPastWithIndex{m}{\pointProcessAtTime{N_{n}(\omega)}{t}}{} \\
		\mbox{(by Assumption \ref{ass:intensity_dominated_1}.(i))} \quad &\leq \intensityFunctional'(m\,|\,\pointProcessAtTime{N_n(\omega)}{t}) \\
		\mbox{(by the definition of $\intensityFunctional'$, Assumption \ref{ass:intensity_dominated_1} and since $N_{n}(\omega)\prec N'(\omega)$)} \quad &\leq \intensityFunctional'(m\,|\,\pointProcessAtTime{N'(\omega)}{t}) =: \intensity'(\omega,t,m).
	\end{align*}
	By Proposition \ref{prop:pathwise_construction_solves_SDE} for $N(\omega)$ and by assumption for $N'(\omega)$, we have that $N(\omega)$ and $N'(\omega)$ both satisfy \eqref{eq:the_SDE_driven_by_Poisson} on $(-\infty, T_{n+1}(\omega)]$, where $\intensity$ is replaced by $\intensity'$ for $N'(\omega)$. Consequently, we must have $N(\omega)\prec N'(\omega)$ on $(-\infty,T_{n+1}(\omega)]$. As, by construction, $N(\omega)$ has mass on $\realLinePositive$ only at the times $T_1(\omega)<T_2(\omega)<\ldots<\infty$, we have shown that  $N(\omega)\prec N'(\omega)$. This implies that $N\prec N'$ $\as$
\end{proof}
This allows us to conclude the proof of Theorem \ref{thm:pathwise_existence} under Assumptions \ref{ass:markspacebounded}, \ref{ass:intensity_dominated_1} and \ref{ass:initial_condition_1}.
\begin{proof}[Proof of Theorem \ref{thm:pathwise_existence}, Part 2]
Repeat the first step of Part 1. Then, by Proposition \ref{prop:construction_is_dominated1}, we deduce that $N\in\mppSpace$ and $T_{\infty}=\infty$ $\as$ We then conclude by repeating the fourth step of Part 1.
\end{proof}

\subsubsection{Existence under Assumptions \ref{ass:markspacebounded}, \ref{ass:sublinear_intensity}, \ref{ass:initial_condition_good_expectation_and_kernel_contracts}.}
To prove Theorem \ref{thm:pathwise_existence} under Assumptions \ref{ass:markspacebounded}, \ref{ass:sublinear_intensity} and \ref{ass:initial_condition_good_expectation_and_kernel_contracts}, we will also use Algorithm \ref{constr:pathwise_construction} to construct a candidate solution, but the almost sure event $F_1$ needs to be replaced by another almost sure event $F_2$ that guarantees that the algorithm is well-defined under these new assumptions. Whereas under Assumptions \ref{ass:intensity_dominated_1} and \ref{ass:initial_condition_1}, we were able to first construct the candidate solution and then dominate it by a solution to the special case $\intensityFunctional=\intensityFunctional'$, here we will dominate the candidate solution while constructing it. The dominating non-explosive marked point process is nothing else than a solution to the Poisson-driven SDE with the Hawkes intensity functional
\begin{align*}
	\dominatingIntensityFunctional : \markSpace\times\boundedlyFiniteMeasures{\theCSMS} &\rightarrow \realLinePositive\cup\{\infty\} \\
	(m,\xi) &\mapsto \dominatingIntensityFunctionalAtMarkGivenPastWithIndex{m}{\xi}{}:= \lambda_{0} + \iint_{(-\infty,0)\times\markSpace}\overline{k}(-t',m',m)\xi(dt',dm'),
\end{align*}
where $\lambda_{0}$ and $\overline{k}$ are as in Assumption \ref{ass:sublinear_intensity}. Indeed, by applying the results of \citet{Massoulie:1998aa:StabilityResults} and Lemma \ref{lem:ground_measure_simple}, we can prove Theorem \ref{thm:pathwise_existence} under  Assumptions \ref{ass:markspacebounded}, \ref{ass:sublinear_intensity} and \ref{ass:initial_condition_good_expectation_and_kernel_contracts} for the special case $\intensityFunctional = \dominatingIntensityFunctional$.
\begin{proof}[Proof of Theorem \ref{thm:pathwise_existence}, Part 3]
Clearly, $\dominatingIntensityFunctional$ satisfies the Lipschitz condition \eqref{eq:lipschitz_condition} with the kernel $\overline{k}$. Under Assumptions \ref{ass:markspacebounded}, \ref{ass:sublinear_intensity}.(ii) and \ref{ass:initial_condition_good_expectation_and_kernel_contracts}.(i), by Theorem 2 of
\citet{Massoulie:1998aa:StabilityResults}, we know that there exists a non-explosive point process $\dominatingMarkedPointProcess:\realisationsSpace\rightarrow\boundedlyFiniteMeasures{\theCSMS}$ that solves \eqref{eq:the_SDE_driven_by_Poisson}, where $\intensityFunctional$ is replaced by $\dominatingIntensityFunctional$, and such that $\dominatingMarkedPointProcess(\cdot\times\markSpace)\in\boundedlyFiniteMeasures{\realLine}$. Moreover, applying Assumptions \ref{ass:sublinear_intensity}.(iii) and  \ref{ass:initial_condition_good_expectation_and_kernel_contracts}.(ii), we obtain that
\begin{align*}
	\dominatingIntensityAtTimeAtMarkWithIndex{\omega,t}{m}{} :&= \dominatingIntensityFunctionalAtMarkGivenPastWithIndex{m}{\pointProcessAtTime{\dominatingMarkedPointProcess(\omega)}{t}}{} =\lambda_{0} + \iint_{(-\infty,t)\times\markSpace}\overline{k}(t-t',m',m)\dominatingMarkedPointProcess(\omega,dt',dm') \\
	&\leq \lambda_0 + \tilde{\intensity}_{\leq 0}(\omega_{\leq 0},t) + \iint_{(0,t)\times\markSpace}\sup_{m''\in\markSpace}\overline{k}(t-t',m',m'')\dominatingMarkedPointProcess(\omega,dt',dm') \\
	&<\infty,\quad \omega\in\realisationsSpace, t\in\realLinePositive, m\in\markSpace,
\end{align*}
which proves that $\sup_{m\in\markSpace}\dominatingIntensityAtTimeAtMarkWithIndex{t}{m}{}<\infty$, for all $t\in\realLinePositive$, $\as$ Hence, by Lemma \ref{lem:ground_measure_simple}, we conclude that $\dominatingMarkedPointProcess$ admits a version such that $\dominatingMarkedPointProcess(\omega)\in\mppSpace$, $\omega\in\realisationsSpace$, meaning that this version solves the Poisson-driven SDE. Conclude by Corollary \ref{cor:intensity_functional_indeed}.
\end{proof}

From now on, denote by $\dominatingMarkedPointProcess$ a solution to the Poisson-driven SDE in the special case $\intensityFunctional=\dominatingIntensityFunctional$ and by $F_2\in\sigmaAlgebra$ the almost sure event that $\dominatingMarkedPointProcess$ solves \eqref{eq:the_SDE_driven_by_Poisson}, where $\intensityFunctional$ is replaced by $\dominatingIntensityFunctional$. The following statement is the analogue of Proposition \ref{prop:pathwise_construction_well_defined} and ensures that Algorithm \ref{constr:pathwise_construction} is well-defined under this different set of assumptions.
\newtheorem{pathwise_construction_well_defined2}[marked_point_process]{Proposition}
\begin{pathwise_construction_well_defined2}[] \label{prop:pathwise_construction_well_defined2}
In Algorithm \ref{constr:pathwise_construction}, where $F_1$ is replaced by $F_2$, under Assumptions  \ref{ass:markspacebounded}, \ref{ass:sublinear_intensity} and \ref{ass:initial_condition_good_expectation_and_kernel_contracts}, we have that, for every $\omega\in F_2$, $\mbox{\textup{card}}(\mathbb{M}_n(\omega))\leq 1$, $N_{n}(\omega)\prec \dominatingMarkedPointProcess(\omega)$, for all $n\in\integers$, and
\begin{align*}
	\left\{ u>T_{n}(\omega) \,:\, \iint_{(T_{n}(\omega),u)\times\markSpace}\int_{(0,\intensityAtTimeAtMarkWithIndex{\omega,t}{m}{n}]}\drivingPoisson(\omega,dt,dm,dz)=0 \right\} \neq \varnothing\quad &\mbox{for all } n\in\integers \mbox{ s.t. } T_{n}(\omega)<\infty.
\end{align*}
Hence, Algorithm \ref{constr:pathwise_construction} is well-defined under these assumptions.
\end{pathwise_construction_well_defined2}
\begin{proof}
	We show the assertion by induction.
	Take any $\omega=(\omega_{\leq 0},\omega_{>0})\in F_2$. Let $n\in\integers$ and, for all $i\in\integers$ such that $i< n$ and $T_{i}(\omega)<\infty$, assume that
	\begin{equation} \label{eq:non_empty_set_induction_proof2}
		\left\{ u>T_{i}(\omega) \,:\, \iint_{(T_{i}(\omega),u)\times\markSpace}\int_{(0,\intensityAtTimeAtMarkWithIndex{\omega,t}{m}{i}]}\drivingPoisson(\omega,dt,dm,dz)=0 \right\} \neq \varnothing.
	\end{equation}
	For all $i\in\integers$ such that $i\leq n$, assume that $N_{i}(\omega)\prec \dominatingMarkedPointProcess$. If $T_{n}(\omega)=\infty$, then, by construction, this is also true for $n+1$.

	Now, assume that $T_{n}(\omega)<\infty$. We first show that \eqref{eq:non_empty_set_induction_proof2} holds also for $i=n$. By adapting the proof of Proposition \ref{prop:construction_is_dominated1} and using Assumption \ref{ass:sublinear_intensity}.(i), we get that $\intensityAtTimeAtMarkWithIndex{\omega,t}{m}{n} \leq \dominatingIntensityAtTimeAtMarkWithIndex{\omega,t}{m}{}$, for all $t>T_n(\omega)$, $m\in\markSpace$. Hence, for any $\varepsilon>0$, we have that
	\begin{align*}
		 \iint_{(T_{n}(\omega),T_{n}(\omega)+\varepsilon)\times\markSpace}\int_{(0,\intensityAtTimeAtMarkWithIndex{\omega,t}{m}{n}]}\drivingPoisson(\omega,dt,dm,dz) &\leq \iint_{(T_{n}(\omega),T_{n}(\omega)+\varepsilon)\times\markSpace}\int_{(0,\dominatingIntensityAtTimeAtMarkWithIndex{\omega,t}{m}{}]}\drivingPoisson(\omega,dt,dm,dz)\\
		 &= \dominatingMarkedPointProcess(\omega,(T_n(\omega),T_n(\omega)+\varepsilon)\times\markSpace) =:\overline{U}_{n}(\omega,\varepsilon) < \infty,
	\end{align*}
	since $\dominatingMarkedPointProcess(\omega,\cdot\times\markSpace)\in\boundedlyFiniteMeasures{\realLine}$. If $\overline{U}_{n}(\omega,\varepsilon)=0$, then clearly \eqref{eq:non_empty_set_induction_proof2} is satisfied for $i=n$. If not, $\dominatingMarkedPointProcess(\omega,\cdot\times\markSpace)$ has a finite number of points in $(T_{n}(\omega),T_{n}(\omega)+\varepsilon)$ and there exists $0<\varepsilon'<\varepsilon$ such that $\overline{U}_{n}(\omega,\varepsilon')=0$, in which case \eqref{eq:non_empty_set_induction_proof2} is again satisfied for $i=n$.

	Second, we show that $\mbox{\textup{card}}(\mathbb{M}_{n+1})\leq 1$ and $N_{n+1}(\omega)\prec \dominatingMarkedPointProcess(\omega)$. If $T_{n+1}(\omega)=\infty$, then this is immediate. If not, as $\intensityAtTimeAtMarkWithIndex{\omega,t}{m}{n} \leq \dominatingIntensityAtTimeAtMarkWithIndex{\omega,t}{m}{}$, for all $t>T_n(\omega)$, it is enough to notice that
	\begin{align*}
		\mathbb{M}_{n+1}(\omega):&=\left\{ m\in\markSpace \,:\, \drivingPoisson(\omega,\{T_{n+1}(\omega)\}\times\{m\}\times\ (0,\intensityAtTimeAtMarkWithIndex{\omega,T_{n+1}(\omega)}{m}{n}])>0\right\} \\
		&\subset \left\{ m\in\markSpace \,:\, \drivingPoisson(\omega,\{T_{n+1}(\omega)\}\times\{m\}\times\ (0,\dominatingIntensityAtTimeAtMarkWithIndex{\omega,T_{n+1}(\omega)}{m}{}])>0\right\} \\
		&= \left\{ m\in\markSpace \,:\, \dominatingMarkedPointProcess(\omega,\{T_{n+1}(\omega)\}\times\{m\})>0\right\} \leq 1,
	\end{align*}
	since $\dominatingMarkedPointProcess(\omega)\in\mppSpace$. As we already know that $N_{n}(\omega)\prec\dominatingMarkedPointProcess(\omega)$, looking at \eqref{eq:point_process_pathwise_recursion} and observing that
	$N_{n+1}(\omega)$ and $N_{n}(\omega)$ only differ by a mass at time $T_{n+1}(\omega)$, we further deduce that $N_{n+1}(\omega)\prec\dominatingMarkedPointProcess(\omega)$.

	Regarding the basis of this induction, it is immediate that $N_0(\omega) \prec \dominatingMarkedPointProcess(\omega)$ since $N_0(\omega) = N_{\leq 0}(\omega_{\leq 0})=\restrictionToNonPositiveOfMeasure{\dominatingMarkedPointProcess}(\omega)$.
\end{proof}

We are now in a position to finish the proof of Theorem \ref{thm:pathwise_existence} under Assumption \ref{ass:markspacebounded}, \ref{ass:sublinear_intensity} and \ref{ass:initial_condition_good_expectation_and_kernel_contracts}.
\begin{proof}[Proof of Theorem \ref{thm:pathwise_existence}, Part 4]
	Let $N:\realisationsSpace\rightarrow\integerValuedMeasures{\theCSMS}$ be given by Algorithm \ref{constr:pathwise_construction} under Assumptions \ref{ass:markspacebounded}, \ref{ass:sublinear_intensity} and \ref{ass:initial_condition_good_expectation_and_kernel_contracts}, where $F_1$ is replaced by $F_2$. By Proposition \ref{prop:pathwise_construction_well_defined2}, this mapping is well defined. Moreover, we notice that Propositions \ref{prop:pathwise_construction_solves_SDE} and \ref{prop:construction_adapted_predictable} still hold under the present assumptions. Hence, we can repeat the first step of Part 1 of the proof. By Proposition \ref{prop:pathwise_construction_well_defined2}, we know that $N_n(\omega)\prec \dominatingMarkedPointProcess(\omega)$ for all $\omega\in F_2$, which implies by construction that $N(\omega)\prec\dominatingMarkedPointProcess(\omega)$ on $\realLinePositive$ for all $\omega\in\realisationsSpace$. Consequently, we have that $N(\omega)\in\mppSpace$ and $T_\infty(\omega)=\infty$ for all $\omega\in\realisationsSpace$, which implies that $N$ solves the Poisson-driven SDE. We conclude again by Corollary \ref{cor:intensity_functional_indeed}.
\end{proof}

\subsection{Strong and weak uniqueness} \label{subsec:uniqueness}

\begin{proof}[Proof of Theorem \ref{thm:uniqueness_pathwise}]
	Let $\tilde{\Omega}\in\sigmaAlgebra$ be the almost sure event that both $N$ and $N'$ solve \eqref{eq:the_SDE_driven_by_Poisson}. Let $\collection{T_n,M_n}{n\in\integers}$ and $\collection{T'_n,M'_n}{n\in\integers}$ be the enumerations in $(0,\infty]\times\markSpace$ to which $N$ and $N'$ are respectively equivalent. Now fix arbitrary $\omega\in\tilde{\Omega}$. We show by strong induction that $T_{n}(\omega)=T'_n(\omega)$ and $M_n(\omega) = M'_n(\omega)$ for all $n\in\integers$.

	Let $n\in\integers$ and assume that  $T_{i}(\omega)=T'_i(\omega)$ and $M_i(\omega) = M'_i(\omega)$ for all $i=1,\ldots,n-1$. By contradiction, assume that $T_n(\omega)\neq T'_n(\omega)$ and, moreover, without loss of generality, that $T_n(\omega)<T'_n(\omega)$. Then, this implies that
	\begin{align*}
		N(\omega,(0,T_n(\omega)]\times\markSpace) &= \int_{(0,T_n(\omega)]}\int_{\markSpace}\int_{(0,\intensityAtTimeAtMarkWithIndex{\omega,t}{m}{}]}M(\omega,dt,dm,dz) = n, \\
		N'(\omega,(0,T_n(\omega)]\times\markSpace) &= \int_{(0,T_n(\omega)]}\int_{\markSpace}\int_{(0,\lambda'(\omega,t,m) ]}M(\omega,dt,dm,dz) = n-1,
	\end{align*}
	where $\intensityAtTimeAtMarkWithIndex{\omega,t}{m}{} = \intensityFunctionalAtMarkGivenPastWithIndex{m}{\pointProcessAtTime{N(\omega)}{t}}{}$ and $\lambda'(\omega,t,m)=\intensityFunctionalAtMarkGivenPastWithIndex{m}{\pointProcessAtTime{N'(\omega)}{t}}{}$. But since $\restrictionToNonPositiveOfMeasure{N(\omega)}=\restrictionToNonPositiveOfMeasure{N'(\omega)}$ and also $T_{i}(\omega)=T'_i(\omega)$ and $M_i(\omega) = M'_i(\omega)$ for all $i=1,\ldots,n-1$, we have that $\pointProcessAtTime{N(\omega)}{t}=\pointProcessAtTime{N'(\omega)}{t}$ for all $t\leq T_{n}(\omega)$ and, thus, $\intensityAtTimeAtMarkWithIndex{\omega,t}{m}{} = \lambda'(\omega,t,m)$ for all $t\leq T_{n}(\omega), m\in\markSpace$. This implies that $n=n-1$ which is a contradiction and, thus, necessarily, $T_n(\omega)=T'_n(\omega)$.

	Similarly, if we assume that $M_n(\omega)\neq M'_n(\omega)$, then this implies that
	\begin{align*}
		N(\omega,\{T_n(\omega)\}\times\{M_n(\omega)\}) &= \int_{\{T_n(\omega)\}}\int_{\{M_n(\omega)\}}\int_{(0,\intensityAtTimeAtMarkWithIndex{\omega,t}{m}{}]}M(\omega,dt,dm,dz) = 1, \\
		N'(\omega,\{T_n(\omega)\}\times\{M_n(\omega)\}) &= \int_{\{T_n(\omega)\}}\int_{\{M_n(\omega)\}}\int_{(0,\lambda'(\omega,t,m) ]}M(\omega,dt,dm,dz) = 0.
	\end{align*}
	But again, since $\intensityAtTimeAtMarkWithIndex{\omega,t}{m}{} = \lambda'(\omega,t,m)$ for all $t\leq T_{n}(\omega), m\in\markSpace$, this leads to the contradiction $1=0$ and, thus, it follows that $M_n(\omega)=M'_n(\omega)$. The same reasoning allows us to prove the basis of the strong induction (i.e., to show that $T_1(\omega)=T'_1(\omega)$ and $M_1(\omega)=M'_1(\omega)$).
\end{proof}

\begin{proof}[Proof of Theorem \ref{thm:uniqueness_weak}]
	Consider the canonical measurable space $\mppMeasureSpace$, where $\mppSpaceBorel=\mppSpace\cap\boundedlyFiniteMeasuresBorel{\theCSMS}$, and the canonical non-explosive marked point process $N$ defined by $N(\omega)=\omega$ for all $\omega\in\mppSpace$. Under both $\inducedProb{N_1}$ and $\inducedProb{N_2}$ (they only charge $\mppSpace$), $N$ satisfies the weak initial condition $N_{\leq 0}$ and admits an intensity given by \eqref{eq:intensity_functional_definition}. We will now apply Theorem 3.4 in \citet[p.~242]{jacod:1975:projection} to show that $\inducedProb{N_1}=\inducedProb{N_2}$.
	By Lemma \ref{lem:history_decomposition}, we have that
	\begin{equation*}
		\internalHistory{N}{t}=\internalHistory{N}{0}\vee \internalHistory{\pointProcessAfterTime{N}{0}}{t}=\internalHistory{\restrictionToNonPositiveOfMeasure{N}}{0}\vee \internalHistory{\restrictionToPositiveOfMeasure{N}}{t},\quad t\in\realLineNonNegative,
	\end{equation*}
	and thus, Assumption A.1 of \citet{jacod:1975:projection} is satisfied (see Remark \ref{rem:jacod_notations}).

	To apply Theorem 3.4 in \citet{jacod:1975:projection}, it remains to verify that the restrictions of $\inducedProb{N_1}$ and $\inducedProb{N_2}$ coincide on $\internalHistory{N}{0}$. Note that $\internalHistory{N}{0}$ is generated by the $\pi$-system $\mathcal{C}$ of sets of the form
	\begin{equation*}
		\{N\in\mppSpace\,:\,N(A_1\times M_1)\geq n_1,\ldots, N(A_k\times M_k)\geq n_k\},\quad n_1,\ldots,n_k\in\integers,\, k\in\integers,
	\end{equation*}
	where $A_1,\ldots,A_k,\in\borelSetsOf{\realLineNonPositive}$ and $M_1,\ldots,M_k\in\markSpaceBorel$. For any such set $F\in\mathcal{C}$, setting $B_i:=A_i\times M_i$ for $i=1,\ldots,k$ and invoking the fact that both $N_1$ and $N_2$ satisfy the weak initial condition $N_{\leq 0}$, we deduce that
	\begin{align*}
	\inducedProb{N_1}(F) &=\inducedProb{N_1}(N(B_1)\geq n_1,\ldots, N(B_k)\geq n_k) = \inducedProb{N_1}(\restrictionToNonPositiveOfMeasure{N}(B_1)\geq n_1,\ldots, \restrictionToNonPositiveOfMeasure{N}(B_k)\geq n_k)\\
		&=\inducedProb{N_{\leq 0}}(\restrictionToNonPositiveOfMeasure{N}(B_1)\geq n_1,\ldots, \restrictionToNonPositiveOfMeasure{N}(B_k)\geq n_k)
		= \inducedProb{N_2}(\restrictionToNonPositiveOfMeasure{N}(B_1)\geq n_1,\ldots, \restrictionToNonPositiveOfMeasure{N}(B_k)\geq n_k)\\
		&= \inducedProb{N_2}(N(B_1)\geq n_1,\ldots, N(B_k)\geq n_k) = \inducedProb{N_2}(F).
	\end{align*}
	Hence, $\inducedProb{N_1}$ and $\inducedProb{N_2}$ coincide on $\mathcal{C}$, a $\pi$-system that contains $\mppSpace$. As a consequence, $\inducedProb{N_1}=\inducedProb{N_2}$ on $\internalHistory{N}{0}$, see for example Lemma 1.17 in \citet[p.~9]{kallenberg2006foundations}, and we can apply Theorem 3.4 in \citet[p.~242]{jacod:1975:projection} to deduce that $\inducedProb{N_1}=\inducedProb{N_2}$ on $\mppMeasureSpace$.
\end{proof}

\theoremstyle{definition}
\newtheorem{jacod_notations}[marked_point_process]{Remark}
\begin{jacod_notations} \label{rem:jacod_notations}
	Let us clarify the relationship between our notations and those in \citet{jacod:1975:projection}. Our canonical measurable space $\mppMeasureSpace$ plays the role of his measurable space $(\Omega, \mathcal{F}_{\infty})$. Our marked point process $\restrictionToPositiveOfMeasure{N}$ corresponds to his marked point process $\mu$. Our probability measures $\inducedProb{N_1}$ and $\inducedProb{N_2}$ are the counterparts of $P$ and $P'$, respectively. Our filtrations $\internalHistoryCollection{N}$ and $\internalHistoryCollection{\restrictionToPositiveOfMeasure{N}}$ correspond to his filtrations $(\mathcal{F}_t)_{t\geq 0}$ and $(\mathcal{G}_t)_{t\geq 0}$, respectively.
\end{jacod_notations}
\theoremstyle{plain}

\appendix
\section{Appendix}\label{sec:appendix}

\subsection{The subspace $\mppSpace$ is Borel}
The following result is unlikely to be original, but we could not find it in \citet{daleyVereJonesVolume2}.
\newtheorem{mppSpace_is_measurable}[marked_point_process]{Lemma}
\begin{mppSpace_is_measurable} \label{prop:mppSpace_is_measurable}
	 The set $\mppSpace$ is a Borel subset of $\boundedlyFiniteMeasures{\theCSMS}$, that is $\mppSpace\in\boundedlyFiniteMeasuresBorel{\theCSMS}$.
\end{mppSpace_is_measurable}
\begin{proof}
	For all $n\in\integers$, define the sets
	\begin{equation*}
		F_{n}:= \left\{ \xi\in\boundedlyFiniteMeasures{\theCSMS} \,:\, \xi([-n,n]\times\markSpace)<\infty,\, \xi(\{t\}\times\markSpace)\leq 1 \mbox{ for all } t\in[-n,n] \right\}
	\end{equation*}
	and notice that $\mppSpace = \bigcap_{n\in\integers}F_{n}$. Next, let $n\in\integers$. By Proposition A2.1.IV in \citet[p.~385]{daleyVereJonesVolume1}, the interval $[-n,n]$ contains a dissecting system $((A_{ij})_{j\in\{1,\dots,j_{i}\}})_{i\in\integers}$ where $A_{ij}\in\borelSetsOf{\realLine}$ for any $j\in\{1,\ldots,j_i\}$ and $i\in\integers$, see Definition A1.6.1 in \citet[p.~382]{daleyVereJonesVolume1}. We show that
	\begin{equation*}
		F_{n} = \left\{ \xi\in\boundedlyFiniteMeasures{\theCSMS} \,:\, \xi([-n,n]\times\markSpace)<\infty,\, \limsup_{i\rightarrow\infty}\sup_{j\in\{1,\ldots,j_{i}\}}\xi(A_{ij}\times\markSpace)\leq 1 \right\} =: G_{n}.
	\end{equation*}
	Let $\xi\in F_{n}$. Then $\xi(\cdot\times\markSpace)$ has finitely many atoms $t_1,\ldots,t_p$ in $[-n,n]$ for some $p\in\integers$ number $p\in\integers$ and their mass cannot exceed one. A key property of the dissecting system is that, for each pair of distinct atoms $t_{q_{1}}$ and $t_{q_{2}}$ with $q_{1}\neq q_{2}$, there exists $n(q_{1},q_{2})\in\integers$ such that, for all $i>n(q_{1},q_{2})$, $t_{q_{1}}\in A_{ij}$ implies $t_{q_{2}}\notin A_{ij}$. Thus, define
	\begin{equation*}
		i^{*}:=\max_{q_{1},q_{2}\in\{1,\ldots,p\},\, q_{1}\neq q_{2}} n(q_{1},q_{2})
	\end{equation*}
	and then, $\xi(A_{ij}\times\markSpace)\leq 1$ for all $j\in\{1,\ldots,j_{i}\}$ and $i>i^{*}$, which implies that
	\begin{equation*}
		\limsup_{i\rightarrow\infty}\sup_{j\in\{1,\ldots,j_{i}\}}\xi(A_{ij}\times\markSpace)\leq 1,
	\end{equation*}
	which in turn indicates that $\xi\in G_{n}$. Now, let $\xi\in G_{n}$ and $t\in[-n,n]$. Another salient property of the dissecting system is that there exists a sequence $(j_{i})_{i\in\integers}$ such that $\xi(\{t\}\times\markSpace) = \lim_{i\rightarrow\infty}\xi(A_{ij_{i}}\times\markSpace)$. But since $\xi\in G_{n}$, we have that
	\begin{equation*}
		\xi(\{t\}\times\markSpace) = \lim_{i\rightarrow\infty}\xi(A_{ij_{i}}\times\markSpace) \leq \limsup_{i\rightarrow\infty}\sup_{j\in\{1,\ldots,j_{i}\}}\xi(A_{ij}\times\markSpace)\leq 1,
	\end{equation*}
	which means that $\xi\in F_{n}$. Now that we have shown that $F_{n}=G_{n}$, we invoke Theorem A2.6.III in \citet[p.~404]{daleyVereJonesVolume1} to deduce that $\xi\mapsto \xi([-n,n]\times\markSpace)$ and $\xi\mapsto\xi(A_{ij}\times\markSpace)$, for any $j\in\{1,\ldots,j_i\}$ and $i\in\integers$, are measurable and use Lemma 1.9 in \citet[p.~6]{kallenberg2006foundations} to conclude that $\xi\mapsto \limsup_{i\rightarrow\infty}\sup_{j\in\{1,\ldots,j_{i}\}}\xi(A_{ij}\times\markSpace)$ is measurable. It then follows that $F_{n}\in\boundedlyFiniteMeasuresBorel{\theCSMS}$, whence $\mppSpace = \bigcap_{n\in\integers}F_{n}\in\boundedlyFiniteMeasuresBorel{\theCSMS}$.
\end{proof}

\subsection{Measurability and continuity properties of shifts and resitrictions}

From \citet[p.~178, Lemma 12.1.I]{daleyVereJonesVolume2}, we know the shift operators are continuous under the $\weakHash$-topology. We are able to go further and show that $\shift{t}\xi$ is actually jointly continuous in $\xi$ and $t$. We also prove that  taking the restriction to the positive or negative real line of a boundedly finite measure is a measurable operation. Moreover, we show that $\pointProcessAtTime{\xi}{t}$ is left-continuous as a function of $t\in\realLine$ for any $\xi\in\boundedlyFiniteMeasures{\anyCSMS}$, which is crucial in our proof that an intensity functional applied to the history of a point process generates a predictable process (Lemma \ref{lem:predictability_functional_process}). Before giving the formal statements, we define some notations. Recall that $\weakHashDistancePlain$ is the $\weakHash$-distance \citep[p.~403]{daleyVereJonesVolume1} on the space $\boundedlyFiniteMeasures{\anyCSMS}$. The open ball with centre $u\in\anyCSMS$ and radius $r$ is denoted by $B_r(u)$.
For any subset $A\subset \anyCSMS$ and $\varepsilon>0$, the $\varepsilon$-neighbourhood of $A$ is defined by $A^\varepsilon:=\bigcup_{a\in A} B_\varepsilon(a)$ and the boundary of $A$ is denoted by $\partial A$.
\newtheorem{shift_jointly_continuous}[marked_point_process]{Lemma}
\begin{shift_jointly_continuous} \label{prop:shift_jointly_continuous}
	When $\boundedlyFiniteMeasures{\realLine\times\anyCSMS}$ is equipped with the $\weakHash$-distance $\weakHashDistancePlain$ and $\boundedlyFiniteMeasures{\realLine\times\anyCSMS}\times\realLine$ is equipped with the product metric, the mapping
	\begin{align*}
		\boundedlyFiniteMeasures{\realLine\times\anyCSMS}\times\realLine &\rightarrow \boundedlyFiniteMeasures{\realLine\times\anyCSMS} \\
		(\xi,t) &\mapsto \shift{t}\xi
	\end{align*}
	is continuous.
\end{shift_jointly_continuous}
\begin{proof}
Let $\xi\in\boundedlyFiniteMeasures{\realLine\times\anyCSMS}$, $t\in\realLine$ and let $(\xi_{n},t_{n})_{n\in\integers}$ be a sequence in $\boundedlyFiniteMeasures{\realLine\times\anyCSMS}\times\realLine$ such that $\weakHashDistance{\xi_{n}}{\xi}\rightarrow 0$ and $t_{n}\rightarrow t$ as $n\rightarrow\infty$. By Proposition A2.6.II in \citet[p.~403]{daleyVereJonesVolume1}, it is enough to show that $\xi_n(A+t_n)\rightarrow\xi(A+t)$ as $n\rightarrow\infty$ for any bounded $A\in\borelSetsOf{\realLine\times\anyCSMS}$ such that $\xi(\partial(A+t))=0$. For such a set $A$, which we can assume without loss of generality to be non-empty, there exists $\delta>0$ such that
	$B_{2\delta}(a_n)\subset (A+t)$ and $\xi(\partial B_{\delta}(a_n))=0$, $n=1,\ldots,N$, where $a_1,\ldots,a_N$ are the atoms of $\xi$ in $A+t$, and such that $\xi((A+t)^\delta)=\xi(A+t)$ with $\xi(\partial((A+t)^\delta))=0$.
Introduce the two bounded sets
\begin{equation*}
	S_{1}:=(A+t)^\delta \setminus (A+t)\quad\mbox{and}\quad S_2:= (A+t)\setminus\left(\bigcup_{n=1}^{N}B_{\delta}(a_n)\right)
\end{equation*}
and notice that $\xi(S_1)=\xi(S_2)=\xi(\partial S_1)=\xi(\partial S_2)=0$. Since $\weakHashDistance{\xi_{n}}{\xi}\rightarrow 0$ as $n\rightarrow\infty$, we have that $\xi_n(S_1)=\xi_n(S_2)=\xi_n(\partial S_1)=\xi_n(\partial S_2)=0$ for $n$ large enough. This implies that, for $n$ large enough, all the atoms of $\xi_n$ in $(A+t)^\delta$ actually lie in $(A+t)$ and their distance to the boundary of $(A+t)$ is bigger than $\delta$ (all the atoms are in the balls $B_{\delta}(a_n)$). This means that, for $n$ large enough, $\xi_n((A+t)\setminus(A+s))=\xi_n((A+s)\setminus(A+t))=0$ for all $s\in\realLine$ such that $|t-s|<\delta$, implying that $\xi_n(A+t)=\xi_n(A+s)$ for all such $n$ and $s$. But for $n$ large enough, we also have that $|t_n - t|<\delta$ and $\xi_n(A+t)=\xi(A+t)$, which finally gives that, for such large enough $n$,
	\begin{equation*}
		\xi_n(A+t_n) = \xi_n(A+t + (t_n-t)) = \xi_n(A+t) = \xi(A+t).\qedhere
	\end{equation*}
\end{proof}

\newtheorem{restrictions_are_measurable}[marked_point_process]{Lemma}
\begin{restrictions_are_measurable}[] \label{prop:restrictions_are_measurable}
	The restrictions $\restrictionToNegativeOfMeasure{\xi}$, $\restrictionToNonPositiveOfMeasure{\xi}$, $\restrictionToPositiveOfMeasure{\xi}$ and $\restrictionToNonNegativeOfMeasure{\xi}$ are measurable mappings from $\boundedlyFiniteMeasures{\realLine\times\anyCSMS}$ into itself.
\end{restrictions_are_measurable}
\begin{proof}
	We prove the assertion for $\restrictionToNegativeOfMeasure{\xi}$, the other three restrictions can be treated similarly. Consider the function $f:\,\boundedlyFiniteMeasures{\realLine\times\anyCSMS}\ni\xi\mapsto\restrictionToNegativeOfMeasure{\xi}\in\boundedlyFiniteMeasures{\realLine\times\anyCSMS}$. Remember that, by Theorem A2.6.III in \citet[p.~404]{daleyVereJonesVolume1}, the Borel $\sigma$-algebra  $\boundedlyFiniteMeasuresBorel{\realLine\times\anyCSMS}$ is generated by the sets
	\begin{equation*}
		F_{A,n}:=\{\xi\in\boundedlyFiniteMeasures{\realLine\times\anyCSMS}\,:\, \xi(A)\in[n,\infty]\}, \quad A\in\boundedlyFiniteMeasuresBorel{\realLine\times\anyCSMS} ,\, n\in\realLine.
	\end{equation*}
	Since
	\begin{equation*}
		f^{-1}(F_{A,n})=\{\xi\in\boundedlyFiniteMeasures{\realLine\times\anyCSMS}\,:\, \xi(A\cap \realLineNegative\times\anyCSMS )\in[n,\infty]\}\in\boundedlyFiniteMeasuresBorel{\realLine\times\anyCSMS},
	\end{equation*}
	we conclude that $f$ is measurable by Lemma 1.4 in \citet[p.~4]{kallenberg2006foundations}.
\end{proof}
\newtheorem{history_left_continuous}[marked_point_process]{Lemma}
\begin{history_left_continuous}[] \label{prop:history_left_continuous}
Let $\xi\in\boundedlyFiniteMeasures{\realLine\times\anyCSMS}$. Then the mapping
\begin{align*}
	\realLine &\rightarrow\boundedlyFiniteMeasures{\realLine\times\anyCSMS} \\
	t &\mapsto \restrictionToNegativeOfMeasure{(\shift{t}\xi )}
\end{align*}
is left continuous when $\boundedlyFiniteMeasures{\realLine\times\anyCSMS}$ is equipped with the $\weakHash$-distance $\weakHashDistancePlain$.
\end{history_left_continuous}
\begin{proof}
	Fix $t\in\realLine$ and take any non-decreasing sequence $\collection{t_n}{n\in\integers}$ in $\realLine$ such that $t_n \uparrow t$ as $n\rightarrow\infty$. By Proposition A2.6.II in \citet[p.~403]{daleyVereJonesVolume1}, it is enough to show that $\restrictionToNegativeOfMeasure{(\shift{t_n}\xi)}(A)\rightarrow\restrictionToNegativeOfMeasure{(\shift{t}\xi)}(A)$ as $n\rightarrow\infty$ for all bounded $A\in\borelSetsOf{\realLine\times\anyCSMS}$ such that $\restrictionToNegativeOfMeasure{(\shift{t}\xi)}(\partial A)=0$. Clearly, it suffices to consider bounded Borel sets $A$ such that $A\subset \realLineNegative\times\anyCSMS$. First, consider the case where $\partial A \subset \realLineNegative\times\anyCSMS$. This implies that
	\begin{equation*}
		(\shift{t}\xi)(\partial A)=\restrictionToNegativeOfMeasure{(\shift{t}\xi)}(\partial A)=0.
	\end{equation*}
	By Lemma \ref{prop:shift_jointly_continuous}, and using again the characterisation of Proposition A2.6.II in \citet[p.~403]{daleyVereJonesVolume1}, this implies that
	\begin{equation*}
		\restrictionToNegativeOfMeasure{(\shift{t_n}\xi)}(A) = (\shift{t_n}\xi)(A) \rightarrow (\shift{t}\xi)(A) = \restrictionToNegativeOfMeasure{(\shift{t}\xi)}(A),\quad n\rightarrow\infty.
	\end{equation*}
	Second, consider the remaining case where $\partial A \cap \{0\}\times\anyCSMS \neq \varnothing$. Then, $\restrictionToNegativeOfMeasure{(\shift{t}\xi)}(\partial A)=0$ does not imply anymore that $(\shift{t}\xi)(\partial A)=0$. However, let $\xi_-$ be the measure $\xi$ that omits all atoms with time coordinate $t$. Then, for this measure $\xi_-$, we have again that
	\begin{equation*}
		(\shift{t}\xi_-)(\partial A)=\restrictionToNegativeOfMeasure{(\shift{t}\xi)}(\partial A)=0.
	\end{equation*}
	Since $t_n \leq t$, and adapting the preceding argument for $\xi_-$, we finally find that
	\begin{equation*}
		\restrictionToNegativeOfMeasure{(\shift{t_n}\xi)}(A) = (\shift{t_n}\xi)(A) = (\shift{t_n}\xi_-)(A) \rightarrow (\shift{t}\xi_-)(A) = \restrictionToNegativeOfMeasure{(\shift{t}\xi)}(A),\quad n\rightarrow\infty. \qedhere
	\end{equation*}
\end{proof}

\subsection{Enumeration representation of marked point processes}  \label{subsec:enumeration_representation}

The following result confirms that a non-explosive enumeration in $\realLinePositive\times\markSpace$ corresponds indeed to a non-explosive marked point process.
\newtheorem{enumeration_is_mpp}[marked_point_process]{Lemma}
\begin{enumeration_is_mpp} \label{prop:enumeration_is_mpp}
Let $(T_{n},M_{n})_{n\in\integers}$ be an enumeration in $\realLinePositive\times\markSpace$ such that $\lim_{n\rightarrow\infty}T_{n}=\infty$ $\as$ Let $F\in\sigmaAlgebra$ be the almost sure event that $\lim_{n\rightarrow\infty}T_{n}=\infty$ and define
\begin{equation*}
	N(\omega):=
	\begin{cases}
		\sum_{n}\delta_{(T_{n}(\omega),M_{n}(\omega))}\indicator_{\{T_{n}(\omega)<\infty\}}, \quad &\mbox{if } \omega\in F, \\
		0, \quad &\mbox{if } \omega\notin F.
	\end{cases}
\end{equation*}
Then, $N$ defines a non-explosive marked point process on $\realLineNonNegative\times\markSpace$.
\end{enumeration_is_mpp}
\begin{proof}
	By Proposition 9.1.X in \citet[p.~13]{daleyVereJonesVolume2}, $N$ defines a non-explosive point process on $\realLineNonNegative\times\markSpace$. Moreover, the monotonicity of the sequence $(T_{n})_{n\in\integers}$ implies that $N(\{t\}\times\markSpace)=0$ or $1$ for all $\omega\in\realisationsSpace$. Also, using that $\lim_{n\rightarrow\infty}T_{n}=\infty$ on $F$, notice that $N(\omega, A\times\markSpace)<\infty$, for every bounded set $A\in\borelSetsOf{\realLineNonNegative}$, for all $\omega\in\realisationsSpace$. This means that $N\in\mathcal{N}_{\realLineNonNegative\times\markSpace}^{\# g}$ and, thus, $N$ defines a non-explosive marked point process.
\end{proof}
Conversely, every non-explosive marked point process generates an enumeration.
\newtheorem{mpp_is_enumeration}[marked_point_process]{Lemma}
\begin{mpp_is_enumeration} \label{prop:mpp_is_enumeration}
Let $N$ be a non-explosive marked point process on $\realLineNonNegative\times\markSpace$ such that $N(\{0\}\times\markSpace)=0$ $\as$ Define the sequence $(T_{n})_{n\in\integers}$ by
$T_{n}:=\sup\{t>0\,:\, N((0,t)\times\markSpace)\leq n\}$.
Then $(T_{n})_{n\in\integers}$ is a non-decreasing sequence of random variables in $(0,\infty]$. Moreover, for each $n\in\integers$, on $\{T_{n}<\infty\}$, one can define $M_{n}$ as the unique element in $\markSpace$ such that $N(\{T_{n}\}\times\{M_{n}\})>0$. On $\{T_{n}=\infty\}$, simply set $M_{n}=m_{\infty}$ for some fixed $m_{\infty}\in\markSpace$. Then, $(T_{n},M_{n})_{n\in\integers}$ is an enumeration in $\realLineNonNegative\times\markSpace$ such that
	\begin{equation} \label{eq:mpp_equals_enumeration}
		N = \sum_{n\in\integers} \delta_{(T_{n},M_{n})}\indicator_{\{T_{n}<\infty\}}
	\end{equation}
	and $\lim_{n\rightarrow\infty}T_{n}(\omega)=\infty$ for all $\omega\in\realisationsSpace$.
\end{mpp_is_enumeration}
\begin{proof}
	We proceed in several steps.
	\begin{enumerate}[label=\textup{(\roman*)}]
		\item For each $n\in\integers$, the mapping $\omega\mapsto T_{n}(\omega)$ is measurable. Indeed, notice that $\{T_{n}< t\} = \{N((0,t)\times\markSpace)>n\}$.  Then recall that, by Theorem A2.6.III in \citet[p.~404]{daleyVereJonesVolume1}, $\boundedlyFiniteMeasures{\realLineNonNegative\times\markSpace}\ni \xi \mapsto \xi((0,t)\times\markSpace)$ is measurable and, thus, as a composition \citep[Lemma 1.7, p.~5]{kallenberg2006foundations},
		the mapping $\omega\mapsto N(\omega,(0,t)\times\markSpace)$ is measurable. Consequently, $\{T_{n}<t\}\in\sigmaAlgebra$. We conclude using Lemma 1.4 in \citet[p.~4]{kallenberg2006foundations} that the mapping $\omega\mapsto T_{n}(\omega)$ is measurable.
		\item Using the fact that the ground measure is simple, $T_{n}<\infty$ implies that $T_{n}<T_{n+1}$. Also, it is easy to check that when $T_{n}=\infty$, then $T_{n+1}=\infty$. Hence, $(T_{n})_{n\in\integers}$ is sequence of random variables in $(0,\infty]$ satisfying the monotonicity of an enumeration.
		\item For each $n\in\integers$, the mapping $\omega\mapsto M_{n}(\omega)$ is well defined (using again the fact that the ground measure is simple). Also, this mapping is measurable. Indeed, let $A\in\markSpaceBorel$ and consider the most delicate case where $m_{\infty}\in A$. Notice that
			\begin{equation*}
				\{M_{n}\in A\} = \left(\{T_{n}<\infty\}\cap\{N(\{T_{n}\}\times A)>0\}\right)\cup \{T_{n}=\infty\}.
			\end{equation*}
			Based on what we have seen so far, we know that $\{T_{n}=\infty\}\in\sigmaAlgebra$. Therefore, it suffices to show that the set $\{T_{n}<\infty\}\cap\{N(\{T_{n}\}\times A)>0\}$ is measurable. To this end, notice that
			\begin{equation*}
				\{T_{n}<\infty\}\cap\{N(\{T_{n}\}\times A)>0\} = \{T_{n}<\infty\}\cap \{ \shift{T_{n}} N(\{0\}\times A)>0\},
			\end{equation*}
			where $\shift{T_{n}}$ is the shift operator defined in Subsection \ref{subsec:shift_histories}. Then, by Lemma \ref{prop:shift_jointly_continuous}, we know that the mapping $\boundedlyFiniteMeasures{\realLineNonNegative\times\markSpace}\times\realLineNonNegative \ni (\xi,t)\mapsto \shift{t}\xi\in\boundedlyFiniteMeasures{\realLineNonNegative\times\markSpace}$ is continuous and thus, by Lemma 1.5 in \citet[p.~4]{kallenberg2006foundations}, measurable. Also, by Lemma 1.8 in \citet[p.~5]{kallenberg2006foundations}, the mapping $\omega\mapsto (N(\omega),T_{n}(\omega))$ is measurable, and thus, as a composition \citep[Lemma 1.7, p.~5]{kallenberg2006foundations}, the mapping $\omega\mapsto\shift{T_{n}(\omega)}N(\omega)$ is measurable. Using again Theorem A2.6.III in \citet[p.~404]{daleyVereJonesVolume1}, we conclude that $\{T_{n}<\infty\}\cap \{ \shift{T_{n}} N(\{0\}\times A)>0\}\in\sigmaAlgebra$ and, thus, the mapping $\omega\mapsto M_{n}(\omega)$ is measurable. So far, these first three steps establish that $(T_{n},M_{n})_{n\in\integers}$ is an enumeration. Moreover, \eqref{eq:mpp_equals_enumeration} holds by construction.
			\item Since $N(\cdot\times\markSpace)\in\boundedlyFiniteMeasures{\realLineNonNegative}$, we have that $\lim_{n\rightarrow\infty}T_{n}(\omega) = \infty$ for all $\omega\in\realisationsSpace$. \qedhere
	\end{enumerate}
\end{proof}
\theoremstyle{definition}
\newtheorem{random_measure_equiv_to_enum}[marked_point_process]{Remark}
\begin{random_measure_equiv_to_enum} \label{rem:random_measure_equiv_to_enum}
	On the one hand, Lemma \ref{prop:enumeration_is_mpp} gives us a mapping  that generates a non-explosive marked point process out of a non-explosive enumeration. One can see that if two non-explosive enumerations are not almost surely equal, then the corresponding non-explosive marked point processes cannot be almost surely equal either. In other words, the mapping of Lemma \ref{prop:enumeration_is_mpp} is injective. On the other hand, Lemma \ref{prop:mpp_is_enumeration} tells us that this mapping is surjective. As a consequence, the above two lemmas tell us that non-explosive enumerations and non-explosive marked point processes are two equivalent ways of looking at the same object.
\end{random_measure_equiv_to_enum}
\theoremstyle{plain}


\subsection{Hawkes functionals} \label{subsec:hawkes_processes}

One can generalise the multivariate linear Hawkes processes presented in the introduction by defining Hawkes functionals as intensity functionals $\intensityFunctional :\markSpace\times\boundedlyFiniteMeasures{\theCSMS}\rightarrow\realLineNonNegative\cup\{\infty\}$ of the form
\begin{equation} \label{eq:typical_event_functional}
	\intensityFunctionalAtMarkGivenPastWithIndex{m}{\xi}{}=\nu(m) + \iint_{(-\infty,0)\times\markSpace}k(-t',m',m)\xi(dt',dm'),\quad m\in\markSpace,\xi\in\boundedlyFiniteMeasures{\theCSMS},
\end{equation}
where $\nu:\markSpace\rightarrow\realLineNonNegative$ and $k:\realLine\times\markSpace\times\markSpace\rightarrow\realLineNonNegative$ are non-negative measurable functions. We show that such event functionals are measurable, so that they are admissible in our framework.
\newtheorem{hawkes_functionals}[marked_point_process]{Proposition}
\begin{hawkes_functionals}[] \label{prop:hawkes_functionals}
Hawkes functionals of the form \eqref{eq:typical_event_functional} are jointly measurable in $m\in\markSpace$ and $\xi\in\boundedlyFiniteMeasures{\anyCSMS}$.
\end{hawkes_functionals}
\begin{proof}
	It will be enough to show that the integral term in \eqref{eq:typical_event_functional}, now denoted by $I(m,\xi)$, is measurable as a function of $m\in\eventSpace$ and $\xi\in\boundedlyFiniteMeasures{\theCSMS}$. First, consider the functions $k$ of the form $k(t',m',m)=\indicator_{S}(t',m',m)$ where $S\in\borelSetsOf{\realLine\times\markSpace\times\markSpace}$ and let $\mathcal{C}$ be the class of sets $S\in\borelSetsOf{\realLine\times\markSpace\times\markSpace}$ such that $(m,\xi)\mapsto I(m,\xi)$ is measurable. By monotone convergence, the class $\mathcal{C}$ is a monotone class (i.e., it is closed under monotonically increasing sequences). Denote by $\mathcal{R}$ the class of sets of the form $\bigcup_{i=1}^{n} A_i\times M'_i\times M_i$ where $A_i\in\borelSetsOf{\realLine}$, $M'_i\in\markSpaceBorel$, $M_i\in\borelSetsOf{\markSpace}, n\in\integers$. This class forms a ring (i.e., it is closed under finite intersections and symmetric differences). Indeed, the difference of unions of Cartesian products is a union of Cartesian product. Moreover, since any union of Cartesian products can be decomposed as a union of disjoint Cartesian products, we have that $\mathcal{R}\subset\mathcal{C}$. Indeed, by Theorem A2.6.III in \citet[p.~404]{daleyVereJonesVolume1}, for any $A\in\borelSetsOf{\realLine}$, $M'\in\markSpaceBorel$, $M\in\borelSetsOf{\markSpace}$, the function
	\begin{equation*}
		(m,\xi)\mapsto \iint_{(-\infty,0)\times\markSpace}\indicator_A(-t')\indicator_{M'}(m')\indicator_M(m)\xi(dt',dm') = \indicator_M(m)\xi((-A)\cap(-\infty,0)\times M')
	\end{equation*}
	is measurable. Then, by the monotone class theorem \citep[p.~369]{daleyVereJonesVolume1}, we have that $\borelSetsOf{\realLine\times\markSpace\times\markSpace}=\sigma(\mathcal{R})\subset \mathcal{C}$. The linearity of the integral implies that $(m,\xi)\mapsto I(m,\xi)$ is measurable for all simple functions $k$ and, by monotone convergence, for all non-negative measurable functions $k$ \citep[Lemma 1.11, p.~7]{kallenberg2006foundations}.
\end{proof}

\section*{Acknowledgements}
We would like to thank Rama Cont, Fabrizio Lillo, Charles-Albert Lehalle, Jean-Philippe Bouchaud and Olav Kallenberg for interesting discussions.
Maxime Morariu-Patrichi gratefully acknowledges the Mini-DTC scholarship awarded by the Mathematics Department of Imperial College London. Mikko S. Pakkanen acknowledges partial support from CREATES (DNRF78), funded by the Danish National Research Foundation.


\end{document}